\documentclass[a4paper,10pt]{article}

\usepackage{amsthm,amssymb,amsmath,amscd,empheq,mathrsfs}

\usepackage{todonotes}
\usepackage[margin=1.0in]{geometry}
\usepackage[shortlabels]{enumitem}
\usepackage{xcolor}
\usepackage[utf8]{inputenc}
\usepackage[T2A]{fontenc}
\usepackage{graphicx}
\usepackage{graphicx}
\usepackage{tikz}
\usepackage{tikz-cd}
\usepackage{titlesec}\titleformat*{\section}{\LARGE\bfseries}
\usepackage{hyperref}
\hypersetup{
  colorlinks   = true,    
  urlcolor     = blue,    
  linkcolor    = blue,    
  citecolor    = red      
}

\usepackage[backend=biber,style=alphabetic]{biblatex}
\makeatletter
\newtheorem*{rep@theorem}{\rep@title}
\newcommand{\newreptheorem}[2]{%
\newenvironment{rep#1}[1]{%
 \def\rep@title{#2 \ref{##1}}%
 \begin{rep@theorem}}%
 {\end{rep@theorem}}}
\makeatother

\newreptheorem{theorem}{Theorem}

\newtheorem{theorem}{Theorem}[section]
\newtheorem{lemma}[theorem]{Lemma}
\newtheorem{proposition}[theorem]{Proposition}
\newtheorem{corollary}{Corollary}[theorem]

\theoremstyle{definition}\newtheorem{definition}[theorem]{Definition}
\theoremstyle{definition}\newtheorem{examples}[theorem]{Example}
\theoremstyle{remark}\newtheorem*{remark}{Remark}

\numberwithin{equation}{theorem}

\DeclareMathOperator{\colim}{colim}

\DeclareMathOperator{\spec}{Spec}
\DeclareMathOperator{\cone}{cone}

\author{Grigory Andreychev}
\title{Pseudocoherent and Perfect Complexes and \protect\\ Vector Bundles on Analytic Adic Spaces}

\date{\vspace{-5ex}}

\begin{document}
\maketitle

\begin{abstract}
Using the new approach to analytic geometry developed by Clausen and Scholze by means of condensed mathematics, we prove that for every affinoid analytic adic space $X$, pseudocoherent complexes, perfect complexes, and finite projective modules over $\mathcal{O}_X(X)$ form a stack with respect to the analytic topology on $X$; in particular, we prove that the category of vector bundles on $X$ is equivalent to the category of finite projective modules over $\mathcal{O}_X(X)$. To that end, we construct a fully faithful functor from the category of complete Huber pairs to the category of analytic rings in the sense of Clausen and Scholze and study its basic properties. Specifically, we give an explicit description of the functor of measures of the analytic ring associated to a complete Huber pair. We also introduce, following the ideas of Kedlaya-Liu, notions of \textit{pseudocoherent module} and \textit{pseudocoherent sheaf} in the context of adic spaces (where the latter can be regarded as an analogue of the notion of coherent sheaf in algebraic geometry) and show that there is a similar equivalence of the corresponding categories.
\end{abstract}

\tableofcontents
\thispagestyle{empty}
\clearpage

\reversemarginpar
\section{Introduction}\label{intro}

In the recent works \cite{Condensed} and \cite{Analytic}, Clausen and Scholze develop a new approach to analytic geometry. Roughly speaking, they first construct an appropriate category of ``analytic rings'' that, on the one hand, is apt for studying topological algebra and, on the other hand, has pleasant algebraic properties, and then glue ``analytic spaces'' out of ``analytic spectra'' in a manner similar to the one used in algebraic geometry. The main objective of this paper is to use their approach to prove certain descent results in the context of analytic adic spaces, some of which have already been known and some of which are new. Before we discuss the subject and the main theorems of this paper in detail, let us first recall the definition of an \textit{analytic Huber ring} due to Kedlaya:

\begin{definition}[\cite{Kdl}, Definition 1.1.2]
Let $A$ be a complete Huber ring. It is called \textit{analytic} if its topologically nilpotent elements generate the unit ideal. 

A complete Huber pair $(A,A^+)$ is said to be \textit{analytic} if $A$ is.
\end{definition}

Analytic complete Huber rings admit a number of equivalent characterizations, which we record here for the convenience of the reader:

\begin{proposition}[\cite{Kdl}, Lemma 1.1.3]
Let $(A,A^+)$ be a complete Huber pair. The following are equivalent:
\begin{itemize}
    \item The ring $A$ is analytic.
    \item Any ideal of definition in any ring of definition of $A$ generates the unit ideal in $A$.
    \item Every open ideal of $A$ is equal to the unit ideal of $A$.
    \item For every nontrivial ideal $I$ of $A$, the quotient topology on $A/I$ is not discrete.
    \item The only discrete topological $A$-module is the zero module.
    \item The set $\mathrm{Spa}(A,A^+)$ contains no point on whose residue field the induced valuation is trivial, i.e., the pre-adic space $\mathrm{Spa}(A,A^+)$ is analytic.
    \item The set $\mathrm{Spa}(A,A^+)$ admits a rational open cover whose elements are the adic spectra of Tate Huber pairs.
\end{itemize}
\end{proposition}

Therefore, analytic Huber rings are simply a natural generalization of Tate rings that encompasses all Huber rings corresponding to arbitrary affinoid subspaces of analytic adic spaces\footnote{Recall that an adic space is said to be analytic if all of its points are analytic.}.

Before we proceed, we should note that there is an unfortunate clash of terminology. The term ``analytic ring'' may be slightly confusing as it can mean either of the two things: \textit{analytic Huber ring} as defined above or \textit{analytic condensed ring} in the sense of Clausen and Scholze. A priori, these two notions have nothing to do with each other (although, one of the goals of this paper is to construct an analytic condensed ring for every analytic Huber ring), but as both of them are used in this paper, it may unclear be at times what is exactly meant. Therefore, we do not use the term ``analytic ring'' unless it is clear from the context which word (``Huber'' or ``condensed'') is omitted.

In what follows, we assume that $(A,A^+)$ is a sheafy analytic complete Huber pair unless otherwise specified, and denote the adic space $\mathrm{Spa}(A,A^+)$ by $X$; additionally, for every affinoid open subspace $U\subset X$, we denote the rings $\mathcal{O}_{X}(U)$ and $\mathcal{O}_{X}^+(U)$ by $A_U$ and $A_U^+$, respectively. Before we state the main theorem of the paper, we need to recall some relevant definitions.

\begin{definition}
An object $M$ of $\mathcal{D}(A)$ is said to be a \textit{pseudocoherent complex} if it is quasi-isomorphic to a complex of the form
\begin{equation}
\label{form}
\ldots\rightarrow P_{n+m}\rightarrow \ldots\rightarrow P_{n+1}\rightarrow P_{n}\rightarrow 0,~\tag{$\dagger$}
\end{equation}
where $P_{i}$ are finite projective modules over $A$. 

An object $M$ is said to be a \textit{perfect complex} if it is quasi-isomorphic to a finite complex of the form~(\ref{form}). A perfect complex $M$ is said to have \textit{tor-amplitude in $[a,b]$}, where $a\leq b$ is a pair of integers, if it is quasi-isomorphic to a complex of the form~(\ref{form}) all of whose terms in cohomological degrees outside the interval $[a,b]$ are trivial.
\end{definition}

We denote the full $\infty$-subcategory of pseudocoherent complexes in $\mathcal{D}(A)$ by $\mathrm{PCoh}_A$, the $\infty$-\hspace*{0pt}subcategory of perfect complexes by $\mathrm{Perf}_A$ and the subcategory of finite projective modules by $\mathrm{FP}_A$. In addition, for every integer $n\in \mathbb{Z}$, we denote the full $\infty$-subcategory of pseudocoherent complexes in $\mathcal{D}(A)^{\leq n}$ by $\mathrm{PCoh}_A^{\leq n}$, and for every pair of integers $a\leq b$, we denote the $\infty$-subcategory of perfect complexes with tor-amplitude in $[a,b]$ by $\mathrm{Perf}_A^{[a,b]}$ (note that $\mathrm{FP}_A$ is clearly equal to $\mathrm{Perf}_A^{[0,0]}$). The main result of this paper is the following theorem (see Theorem~\ref{descent of everything} below):

\begin{theorem}
\label{to refer}
The functors sending every affinoid subspace $U\subset X$ to the $\infty$-categories $\mathrm{PCoh}_{A_U}$, $\mathrm{PCoh}_{A_U}^{\leq n}$, $\mathrm{Perf}_{A_U}$, $\mathrm{Perf}_{A_U}^{[a,b]}$, and $\mathrm{FP}_{A_U}$, respectively, define sheaves of $\infty$-categories on $X$.
\end{theorem}

For finite projective modules, this statement has already been known and the first proof is due to Kedlaya and Liu \cite[Theorem 1.4.2]{Kdl} (it was originally proved in \cite[Theorem 2.7.7]{KLI} in the case where $A$ is Tate). However, the descent results for pseudocoherent and perfect complexes are new. 

In \cite{KLII} (see also \cite[§1]{Kdl}), Kedlaya and Liu also introduce the notion of a \textit{pseudocoherent module} (which generalizes the notion of a finite projective module). These modules give rise to a special class of sheaves on analytic adic spaces that can be thought of as an analogue of coherent sheaves in algebraic geometry. As for finite projective modules, they prove in \cite[Theorem 1.4.14]{Kdl} that these modules satisfy descent with respect to the analytic topology (see \cite[Theorem 2.4.15]{KLII} for the original proof when $A$ is Tate). In Section~\ref{vb}, we discuss in detail the notion of pseudocoherence in the context of analytic adic spaces and give a new proof\footnote{We should note, however, that our definition of pseudocoherent modules is slightly different, and, therefore, we actually establish a slightly different result; see the remark after \ref{condensed ps}.} of \cite[Theorem 1.4.14]{Kdl} (see Theorem~\ref{pseudo stack} below).

The proof of Theorem~\ref{to refer} proceeds roughly as follows. Using techniques of condensed mathematics, we first construct a larger category $\mathcal{D}(X)$ of ``quasi-coherent sheaves'' on $X$ and then prove that these ``sheaves'' satisfy descent with respect to the analytic topology. After that, we prove Theorem~\ref{to refer} by isolating the categories $\mathrm{PCoh}_{A_U}$, $\mathrm{PCoh}_{A_U}^{\leq n}$, $\mathrm{Perf}_{A_U}$, $\mathrm{Perf}_{A_U}^{[a,b]}$, and $\mathrm{FP}_{A_U}$ inside $\mathcal{D}(X)$. 

Therefore, we need to produce that category of ``quasi-coherent sheaves'' first. To this end, for every complete Huber pair $(A,A^+)$, we construct a condensed analytic ring, denoted $(A,A^+)_{\blacksquare}$ or, if $A^+$ is equal to $A^{\circ}$, just $A_{\blacksquare}$ for brevity. The construction of $(A,A^+)_{\blacksquare}$ extends the one for discrete Huber pairs used in \cite{Condensed} and is outlined in \cite{Analytic} (see Proposition 13.16 there); we give a complete proof that such a condensed analytic ring is indeed analytic presupposing this only for $R_{\blacksquare}$, where $R$ is a finitely generated $\mathbb{Z}$-algebra (see \cite[Theorem 8.13]{Condensed} for a proof of the latter). Furthermore, we also introduce the notion of a \textit{quasi-finitely generated submodule} of $A$ in Section~\ref{general1}, which we use to prove that $(A,A^+)_{\blacksquare}$ is always $0$-truncated in the sense of \cite{Analytic} and obtain the following explicit description of the constructed functor of measures (see Theorem~\ref{freeobjects} below):

\begin{theorem}
For every profinite set $S$, we have \[(A,A^+)_{\blacksquare}[S]\cong \underset{B\rightarrow A^+,\ M}{\colim} \underset{J}{\prod} \underline{M}\subset \underset{J}{\prod} \underline{A},\] where $J$ is a set such that $\mathbb{Z}_{\blacksquare}[S]\cong \underset{J}{\prod} {\mathbb{Z}}$ and the colimit is taken over all finitely generated subrings $B\subset A^+$ and all quasi-finitely generated $B$-submodules $M$ of $A$.
\end{theorem}

After that, for an affinoid analytic space $X=\mathrm{Spa}(A,A^+)$, we define the category of ``quasi-coherent sheaves'' $\mathcal{D}(X)$ to be the category $\mathcal{D}((A,A^+)_{\blacksquare})$. The main property of $\mathcal{D}(X)$ that facilitates our proof of Theorem~\ref{to refer}, is that it satisfies descent with respect to the analytic topology, which is made precise by the following theorem (see Theorem~\ref{maindescent} below):

\begin{theorem}
\label{main stuff}
Let $X$ be an analytic adic space and let $U$ denote an arbitrary affinoid subspace of $X$. Then the association $U\mapsto \mathcal{D}((\mathcal{O}_X(U),\mathcal{O}_X^+(U))_{\blacksquare})$ defines a sheaf of $\infty$-categories on $X$. Explicitly, let $\{U_i\}_{i\in I}$ be a family of open affinoid subspaces of $X$ covering $U$. Let $\mathcal{C}_I$ denote the category of finite subsets of $I$, and consider the functor mapping $J\in \mathcal{C}_I$ to the open subspace $U_J\overset{\mathrm{def}}{=}\underset{i\in J}{\cap} U_i$. Then the natural functor
\[\mathcal{D}((\mathcal{O}_X(U),\mathcal{O}_X^+(U))_{\blacksquare})\rightarrow \underset{J\in \mathcal{C}_I}{\lim}\mathcal{D}((\mathcal{O}_X(U_J),\mathcal{O}_X^+(U_J))_{\blacksquare})\]
is an equivalence of $\infty$-categories.
\end{theorem}

As explained in Section~\ref{vb}, the derived $\infty$-category of $A$-modules $\mathcal{D}(A)$ can be embedded fully faithful into $\mathcal{D}(X)$ by means of the \textit{condensification functor}; the objects in its essential image are called \textit{(relatively) discrete}. Therefore, the $\infty$-categories $\mathrm{PCoh}_{A_U}$, $\mathrm{PCoh}_{A_U}^{\leq n}$, $\mathrm{Perf}_{A_U}$, $\mathrm{Perf}_{A_U}^{[a,b]}$, and $\mathrm{FP}_{A_U}$ can be viewed as full $\infty$-subcategories of $\mathcal{D}(X)$. Let us briefly explain how to isolate these categories in $\mathcal{D}(X)$ and prove that they also satisfy descent. The main step is to handle the case of $\mathrm{PCoh}_{A}$. One of the reasons necessitating the use of condensed mathematics is that, in general, discrete object do not form a sheaf of $\infty$-categories. Therefore, we need to work around that problem in order to prove Theorem~\ref{to refer}. In Subsection~\ref{vb3}, we introduce and discuss in detail the notions of \textit{pseudocoherent}, \textit{dualizable}, \textit{nuclear} and \textit{compact} objects in $\mathcal{D}(X)$. In particular, we prove that these classes satisfy descent. We also show that every pseudocoherent complex in $\mathcal{D}(A)$ is pseudocoherent and nuclear as an object of $\mathcal{D}(X)$. Furthermore, we prove that a discrete object in $\mathcal{D}(X)$ is pseudocoherent if and only if it comes from a pseudocoherent complex in $\mathcal{D}(A)$. The crux of the argument is to show that every nuclear pseudocoherent object in $\mathcal{D}(X)$ is in fact discrete. After that, we easily deduce Theorem~\ref{to refer} for $\mathrm{PCoh}_{A}$ from the descent results for nuclear and pseudocoherent object.

The paper is organized as follows. First, we begin with brief recollections on condensed mathematics. In particular, we record some basic facts about analytic rings and fix notation. For a thorough introduction to these concepts, see \cite{Condensed} and \cite{Analytic}. We record the main notational and terminological conventions of the paper below, so a reader familiar with the language of condensed mathematics may safely skip Section~\ref{analytic}. The third section is concerned with basic questions about Huber pairs and analytic rings associated to them. We develop there the technical tools used later to prove the aforementioned theorems, and describe the structure of the analytic ring associated to an arbitrary complete Huber pair; in particular, we prove that it is always $0$-truncated. Section~\ref{descentchapter} is entirely devoted to a proof of Theorem~\ref{main stuff}. Finally, we study finite projective modules and pseudocoherent and perfect complexes over analytic complete Huber rings and prove Theorem~\ref{to refer} in the last section.

In what follows, we adopt the following notational and terminological conventions:  

\begin{itemize}
    \item All rings, including condensed ones, are assumed commutative.
    \item In this article, we abuse terminology by using the term ``the derived category of an abelian category $\mathcal{A}$'' to mean ``the derived $\infty$-category of an abelian category $\mathcal{A}$'' in the sense of \cite{HA}, which we denote $\mathcal{D}(\mathcal{A})$
    \item We use the term ``extremally disconnected set'' to refer to an \textit{extremally disconnected compact Hausdorff space}.
    \item In what follows, whenever we use the term ``analytic ring'' outside of Section~\ref{analytic}, we use it in the sense of \cite[Definition 12.1]{Analytic}, i.e., we actually mean ``animated analytic ring''. Therefore, we also work, as has been already mentioned, with $\infty$-enhancements of derived categories and use the notation $(\mathcal{A},\mathcal{M})$ for analytic rings, where $\mathcal{A}$ stands for the underlying animated condensed ring and $\mathcal{M}$ for the functor of measures; accordingly, the derived category of complete modules is denoted by $\mathcal{D}(\mathcal{A},\mathcal{M})$. However, we slightly abuse terminology by generally omitting the word  ``animated''.
    \item As was already transparent in the introduction above, we will generally abuse notation by denoting both the analytic ring associated to a complete Huber pair $(A,A^+)$ and its functor of measures by $(A,A^+)_{\blacksquare}$. Similarly, we denote both simply $A_{\blacksquare}$, if the ring of integral elements $A^+$ is equal to $A$.
    \item Let $(\mathcal{A},\mathcal{M})$ be an analytic ring. An object of $\mathrm{Mod}_{\mathcal{A}}^{\mathrm{cond}}$ (resp. $\mathcal{D}({\mathcal{A}})$) is said to be \textit{$(\mathcal{A},\mathcal{M})$-complete} if and only if it lies in $\mathrm{Mod}_{(\mathcal{A},\mathcal{M})}^{\mathrm{cond}}$ (resp. $\mathcal{D}(\mathcal{A},\mathcal{M})$).
\end{itemize}

{\large\textbf{Acknowledgments}}

This paper is an expanded version of my Master's thesis at the University of Bonn. I would like to express my deepest gratitude to my advisor, Peter Scholze, who suggested its topic and shared with me so much invaluable mathematical insight with great generosity. I am also indebted to him for careful readings of several ugly drafts of this article. Without his help and encouragement, I would have never published this paper. I thank Arthur-César Le Bras and Alberto Vezzani for their questions, which helped me better understand my own work. I am especially grateful to Dmitry Kubrak for many helpful comments on a previous version of this paper. During the time spent at the University of Bonn as a master student, I was supported by the Hausdorff Center for Mathematics. I thank both institutions for their care and support.
\section{Reminders on analytic rings}\label{analytic}
The sole purpose of this section is to review some basic definitions from condensed mathematics and fix notation. The main concept that we recall is the notion of (animated) analytic ring in the sense of \cite{Analytic}, and a reader familiar with it may skip ahead to Section~\ref{general}. As mentioned in the introduction, all rings, including condensed ones, are assumed commutative. We begin with the following fundamental definition:

\begin{definition}[\cite{Condensed}, Definition 1.2]
Consider the pro-étale site $\ast_{\mathrm{pro\acute{e}t}}$ of the point $\ast$, i.e., the category $\mathrm{ProFin}\cong \mathrm{Pro}(\mathrm{Fin})$ of profinite sets with covers given by finite families of jointly surjective maps. A \textit{condensed set} is a sheaf on $\ast_{\mathrm{pro\acute{e}t}}$. Similarly, a \textit{condensed abelian group, ring, etc.,} is a sheaf of abelian groups, rings, etc., on $\ast_{\mathrm{pro\acute{e}t}}$.\footnote{This definition is not quite the ``official'' definition. As explained in \cite{Condensed}, this version has minor set-theoretic issues; see \cite{Condensed} for a proper discussion and the correct definition. In this paper, we will completely ignore this kind of problems. } The category of condensed sets is denoted by $\mathrm{Cond}$ and the category of condensed modules over a condensed ring $R$ by $\mathrm{Mod}^{\mathrm{cond}}_{R}$.
\end{definition}

The main advantage of condensed mathematics is that it provides a framework which is well-suited for dealing with problems that arise in topological algebra. Namely, when one considers algebraic objects, i.e., groups, modules, vector spaces, etc., with some non-trivial topology (that is, not discrete) and continuous morphisms between them, the resultant category is usually poorly behaved: for instance, it is almost never abelian. One particularly easy and enlightening example is given by the identity map $\mathbb{R}_{\mathrm{disc}}\rightarrow \mathbb{R}$, where the left hand side is endowed with the discrete topology and the right hand side with the classical one. It is readily verified that both kernel and cokernel of this map in the category of topological abelian groups are trivial. However, this map is evidently not an isomorphism. The problem is that in the present context, the notions of kernel and cokernel detect only the difference between the sets of continuous maps from a one-point set to $\mathbb{R}_{\mathrm{disc}}$ and $\mathbb{R}$, respectively, and know nothing about maps from arbitrary compact Hausdorff spaces. That is, the fact that both kernel and cokernel of the map $\mathbb{R}_{\mathrm{disc}}\rightarrow \mathbb{R}$ are trivial reflects precisely the fact that this map induces a bijection on the sets of continuous maps from $\{\ast\}$ to $\mathbb{R}_{\mathrm{disc}}$ and $\mathbb{R}$, respectively. It suggests that the key to the problem is to keep track of \textit{all} continuous maps from \textit{all} compact Hausdorff spaces. However, as every compact Hausdorff space admits a surjection from a profinite set, it can be easily proven that it is enough to know only how profinite sets map to a given (compactly generated weakly Hausdorff) topological abelian group (or space, for that matter) in order to recover that group. Therefore, we may translate the problem at hand into the language of condensed mathematics as follows. Consider the functor from the category of compactly generated weakly Hausdorff topological spaces\footnote{One of the reasons for this restriction, natural to most topologists, stems from the fact that the definition above is actually incorrect and the category of condensed sets is not the category of sheaves on anything. For more on that see \cite[Appendix to Lecture II]{Condensed}.} to the category of condensed sets that sends such a topological space $T$ to the condensed set $\underline{T}$ given by $\underline{T}(S)=\mathrm{Cont}(S,T)$ for every profinite set $S$. Applying this functor, we may view the map $\mathbb{R}_{\mathrm{disc}}\rightarrow \mathbb{R}$ as a morphism of condensed abelian groups with nontrivial cokernel. This translation not only obviates the problem at hand, but also provides a systematic way to deal with such topological issues.  The next two theorems explain the properties of this translation and, hopefully, indicate its importance. For the convenience of the reader, we first recall the following basic definition:

\begin{definition}
Let $\mathcal{A}$ be an abelian category that admits all filtered colimits. An object $M$ is called \textit{compact}, if the functor $\mathrm{Hom}_{\mathcal{A}}(M,-)$ commutes with all filtered colimits. An object $X$ of $\mathcal{D}(\mathcal{A})$ is called \textit{compact}, if the functor $\mathrm{R}{\mathrm{Hom}}_{\mathcal{A}}(X, -)$ commutes with all filtered colimits or, equivalently, all direct sums.

Let $\mathcal{C}$ be a category that admits all coproducts. A family of objects $S$ in $\mathcal{C}$ is said to \textit{generate $\mathcal{C}$} if the following equivalent conditions hold:
\begin{itemize}
    \item If $g_1,g_2:c\rightarrow d$ is a pair of morphisms in $\mathcal{C}$ such that $g_1\circ f=g_2\circ f$ for all morphisms $f:s\rightarrow c$, where $s$ is an object of $S$, then $g_1=g_2$.
    \item For every object $c$ of $\mathcal{C}$, there exists an epimorphism of the form $\coprod_{i\in I} s_i \rightarrow c$, where $I$ is an index set and $s_i$ is an object of $\mathcal{C}$ for every $i\in I$.
\end{itemize}

\end{definition}
\begin{theorem}[\cite{Analytic}, Proposition 1.2]
Consider the functor from  the category of compactly generated weakly Hausdorff topological spaces to the category of condensed sets given by the formula $T\mapsto \underline{T}$. 
\begin{enumerate}
    \item It has a left adjoint $X\mapsto X(\ast)_{\mathrm{top}}$, which sends a condensed set $X$ to its underlying set $X(\ast)$ (i.e., the evaluation of the functor $X$ at the point $*$) endowed with the quotient topology arising from the map
    \[\underset{S,a\in X(S)}{\coprod}S\rightarrow X(\ast),\]
    where $S$ ranges over all profinite sets.
    \item It is fully faithful and induces an equivalence between the category of compact Hausdorff spaces and quasi-compact quasi-separated\footnote{Here the notions of quasi-compactness and quasi-separatedness are general notions applicable to sheaves on any (coherent) site. In our case, a condensed set $X$ is quasi-compact if there is some profinite set $S$ with a surjective map $\underline{S}\rightarrow X$. A condensed set $X$ is quasi-separated if for every pair of profinite sets $S_1, S_2$ with maps to $X$, the fiber product $S_1 \underset{X}{\times} S_2$ is quasi-compact.} condensed sets.
    \item It induces an equivalence between the category of ind-compact Hausdorff spaces of the form “$\varinjlim_i$”$T_i$ where all transition maps $T_i\rightarrow T_j$ are closed immersions, and the category of quasi-separated condensed sets.
\end{enumerate}
\end{theorem}

\begin{theorem}[\cite{Condensed}, Theorem 1.10]\label{main}
The category of condensed abelian groups is an abelian category, which satisfies Grothendieck’s axioms (AB3), (AB4), (AB5), (AB6), (AB3*), and (AB4*). Furthermore, it is generated by compact projective objects.
\end{theorem}

Explicitly, there is a functor from the category of profinite sets (in fact, all compact Hausdorff spaces) to the category of condensed abelian groups which sends a profinite set $S$ to the \textit{free condensed abelian group $\mathbb{Z}[S]$}\footnote{In this paper, we employ underscores to denote condensed abelian groups coming from topological ones; however, they are often omitted if the abelian group is discrete and the context is clear.} \textit{on $S$}. The collection of condensed abelian groups $\mathbb{Z}[S]$ for all extremally disconnected spaces $S$ forms a family of compact projective generators of the category of condensed abelian groups. Recall that a compact Hausdorff topological space $T$ is called \textit{extremally disconnected} if every surjection $T'\rightarrow T$ from a compact Hausdorff space admits a section. In what follows, we will slightly abuse terminology by referring to such a space as an \textit{extremally disconnected set}. 

The significance of this notion is explained by the fact that extremally disconnected sets are precisely the projective objects of the category of compact Hausdorff spaces. In fact, they form a family of compact projective generators of that category as every compact Hausdorff space admits a surjection from an extremally disconnected set. Moreover, every such set is profinite, and, therefore, we do not go beyond the limits of $\ast_{\mathrm{pro\acute{e}t}}$ to construct the family above.

There is an obvious analogue of Theorem~\ref{main} for the category of condensed modules over a condensed ring $R$, which can be proved by essentially the same argument. Moreover, the category $\mathrm{Mod}^{\mathrm{cond}}_{R}$ is symmetric monoidal and Cartesian closed. Its tensor product is defined as the sheafification of the presheaf tensor product, and for every pair of condensed $R$-modules $M,N$, the condensed $R$-module $\underline{\mathrm{Hom}}_R(M,N)$ is given by the formula
\[\underline{\mathrm{Hom}}_R(M,N)(S)={\mathrm{Hom}}_R(M\underset{R}{\otimes} R[S],N)\]
for every profinite set $S$.

One of the key advantages of condensed mathematics is that it has proven capable of providing the right algebraic tools to study analysis. In particular, there is an especially well-behaved non-archimedean notion of completeness, which will play a vital role in the paper. Let $S=\varprojlim_i S_i$ be a profinite set written as a limit of finite sets $S_i$. Consider the condensed abelian group $\mathbb{Z}_{\blacksquare}[S]\overset{\mathrm{def}}{=}\varprojlim_i{\mathbb{Z}}[S_i]$. Objects of this form should be thought of as certain extensions of free condensed abelian groups that contain not only finite but also ``convergent'' sums. For example, we will later show (see Lemma~\ref{powerseries}) that the cokernel of the morphism ${\mathbb{Z}}\rightarrow {\mathbb{Z}}_{\blacksquare}[\mathbb{N}\cup \{ \infty\}]$ induced by the map $\{\infty\}\rightarrow \mathbb{N}\cup \{\infty\}$ is isomorphic to the condensed ring $\mathbb{Z}[[U]]$ of formal power series in one variable with coefficients in ${\mathbb{Z}}$, which is isomorphic to $\underset{\mathbb{N}}{\prod}{\mathbb{Z}}$ as a condensed abelian ring. To explain in more detail how these condensed modules, named \textit{solid} by its creators, serve as a basis for framing non-archimedean analysis in algebraic terms, we need to introduce the following notion:

\begin{definition}[\cite{Condensed} Definitions 7.1 \& 7.4 and \cite{Analytic} Definition 6.4]\label{analytic ring}
A \textit{pre-analytic ring} $(\mathcal{A},\mathcal{M})$ is a condensed ring $\mathcal{A}$ (called the \textit{underlying condensed ring of $(\mathcal{A},\mathcal{M})$}) equipped with a functor 
\[\{\mathrm{extremally\ disconnected\ sets}\}\rightarrow \mathrm{Mod}^{\mathrm{cond}}_{\mathcal{A}},\ S\mapsto \mathcal{M}[S],\]
called the \textit{functor of measures of $(\mathcal{A},\mathcal{M})$}\footnote{We should make a small remark regarding one slight abuse of notation committed in the paper: we will denote both the analytic ring associated to a complete Huber pair $(A,A^+)$ and its functor of measures by $(A,A^+)_{\blacksquare}$. Similarly, we will denote both simply $A_{\blacksquare}$, if the ring of integral elements $A^+$ is equal to $A^{\circ}$.}, that sends finite disjoint unions to products, and a natural transformation $\Phi_{(\mathcal{A},\mathcal{M})}: \underline{S} \rightarrow \mathcal{M}[S]$ of functors  
\[\{\mathrm{extremally\ disconnected\ sets}\}\rightarrow \mathrm{Cond}.\] A pre-analytic ring $(\mathcal{A},\mathcal{M})$ is said to be \textit{analytic}, if for every complex
\[C:\dots \rightarrow C_i\rightarrow \dots \rightarrow C_1\rightarrow C_0\rightarrow 0 \]
of $\mathcal{A}$-modules, such that each $C_i$ is a direct sum of objects of the form $\mathcal{M}[T]$ for various extremally disconnected sets $T$, the map $\mathrm{R}\underline{\mathrm{Hom}}_{\mathcal{A}}(\mathcal{M}[S], C) \rightarrow \mathrm{R}\underline{\mathrm{Hom}}_{\mathcal{A}}({\mathcal{A}}[S], C)$ is an isomorphism for all extremally disconnected sets $S$. Equivalently, $(\mathcal{A},\mathcal{M})$ is analytic if and only if for any families $\{S_i\}_{i\in I},\{S_j\}_{j\in J}$ of extremally disconnected sets and a map
\[f: \underset{i\in I}{\bigoplus}{\mathcal{M}}[S_i]\rightarrow \underset{j\in J}{\bigoplus}{\mathcal{M}}[S_j]\]
with kernel $K$, the map $\mathrm{R}\underline{\mathrm{Hom}}_{\mathcal{A}}(\mathcal{M}[S], K) \rightarrow \mathrm{R}\underline{\mathrm{Hom}}_{\mathcal{A}}(\mathcal{A}[S], K)$ is an isomorphism for all extremally disconnected sets $S$. 

An analytic ring $(\mathcal{A},\mathcal{M})$ is called \textit{complete}, if the map $\mathcal{A}\rightarrow\mathcal{M}[\ast]$ is an isomorphism.
\end{definition}

Before we explain the importance of this notion, let us supply some examples.

\begin{examples}\label{examples}
\leavevmode
\begin{enumerate}
    \item\label{1} The most important example, at least for the present paper, is $\mathbb{Z}_{\blacksquare}$. Its underlying condensed ring is simply $\mathbb{Z}$ and the functor of measures is given by the formula $S\mapsto \mathbb{Z}_{\blacksquare}[S]$.\footnote{Note that above we defined $\mathbb{Z}_{\blacksquare}[S]$ for every profinite set $S$. However, the functor of measures of an analytic ring is a priory defined only on extremally disconnected sets. We will address this discrepancy below; the same comment applies to Examples~\ref{examples}.\ref{2} \& \ref{examples}.\ref{3}.} This analytic ring serves as a foundation upon which we build the construction of the analytic ring associated to a complete Huber pair in Section~\ref{general}. The reader should note, however, that the proof of the analyticity of $\mathbb{Z}_{\blacksquare}$ is rather non-trivial, see \cite[Lectures V and VI]{Condensed}.
    \item\label{2} The construction of $\mathbb{Z}_{\blacksquare}$ may be generalized to arbitrary (discrete) finitely generated $\mathbb{Z}$-algebras as follows. Let $A$ be such an algebra and $S$ a profinite set. Write $S$ as a limit $\varprojlim_i S_i$ of finite sets and define
    \[ A_{\blacksquare}[S]\overset{\mathrm{def}}=\varprojlim_i A[S_i].\]
    Then the condensed ring $A$ and the functor
    \[\{\mathrm{extremally\ disconnected\ sets}\}\rightarrow \mathrm{Mod}^{\mathrm{cond}}_{A},\ S\mapsto A_{\blacksquare}[S]\]
    form an analytic ring, which we denote by $A_{\blacksquare}$. For a proof, see \cite[Appendix to Lecture VIII]{Condensed}.
    \item\label{3} This construction can also be applied to $\mathbb{Z}_p$. Let $S$ be a profinite set and write it as a limit $\varprojlim_i S_i$ of finite sets. Define
    \[ \mathbb{Z}_{p,\blacksquare}[S]\overset{\mathrm{def}}=\varprojlim_{i}\underline{\mathbb{Z}}_p[S_i].\]
    Then the condensed ring $\underline{\mathbb{Z}}_p$ and the functor
    \[\{\mathrm{extremally\ disconnected\ sets}\}\rightarrow \mathrm{Mod}^{\mathrm{cond}}_{\underline{\mathbb{Z}}_p},\ S\mapsto \mathbb{Z}_{p,\blacksquare}[S]\]
    form an analytic ring, which we denote by $\mathbb{Z}_{p,\blacksquare}$. For a proof, see \cite[Proposition 7.9]{Condensed}.\footnote{In Section~\ref{general}, we construct an analytic ring for every complete Huber pair. It may seem at first that the new $\mathbb{Z}_{p,\blacksquare}$ is somewhat different. However, it is easy to check that the two analytic rings are isomorphic.}
    \item\label{4} The main objective of \cite{Analytic} is to construct certain analytic rings with underlying rings $\underline{\mathbb{R}}$ and $\underline{\mathbb{C}}$ that make it possible to frame real and complex analysis, respectively, in terms of condensed mathematics. In fact, Clausen and Scholze construct such rings for each real $0<p\leq 1$. As explained below, one obtains a class of ``complete'' $\underline{\mathbb{R}}$- and $\underline{\mathbb{C}}$-condensed modules for each $0<p\leq 1$, which are termed \textit{$p$-liquid $\mathbb{R}$- and $\mathbb{C}$-vector spaces}, respectively. However, the verification that they are analytic is \textit{very complicated}. See especially \cite[Lecture VI]{Analytic}.
\end{enumerate}
\end{examples}

The notion of analytic ring enables one to talk about ``complete'' condensed modules. Their principal application is to transfer problems arising in non-archimedean, real, and complex geometry into the setting of condensed mathematics. The following theorem spells out their definition, encapsulates their basic properties, and, hopefully, indicates their importance:

\begin{theorem}[\cite{Condensed}, Proposition 7.5]\label{wichtig}
Let $(\mathcal{A},\mathcal{M})$ be an analytic ring.
\begin{enumerate}
    \item The full subcategory \[\mathrm{Mod}_{(\mathcal{A},\mathcal{M})}^{\mathrm{cond}}\subset \mathrm{Mod}^{\mathrm{cond}}_{{\mathcal{A}}}\]
    of all $\mathcal{A}$-modules $M$, such that for every extremally disconnected set $S$, the map \[{\mathrm{Hom}}_{\mathcal{A}}(\mathcal{M}[S], M)\rightarrow M(S)\] is an isomorphism, is an subabelian category stable under all limits, colimits, and extensions. Objects of the form $\mathcal{M}[S]$, where $S$ is an extremally disconnected profinite set, constitute a family of compact projective generators of $\mathrm{Mod}_{(\mathcal{A},\mathcal{M})}^{\mathrm{cond}}$. The inclusion functor admits a left adjoint
    \[\mathrm{Mod}_{\mathcal{A}}^{\mathrm{cond}}\rightarrow \mathrm{Mod}^{\mathrm{cond}}_{(\mathcal{A},\mathcal{M})},\ M\mapsto M\underset{\mathcal{A}}{\otimes}(\mathcal{A},\mathcal{M}),\]
    which is the unique colimit-preserving extension of the functor given by $\mathcal{A}[S]\mapsto \mathcal{M}[S]$. Finally, there is a unique symmetric monoidal tensor product $-\underset{(\mathcal{A},\mathcal{M})}{\otimes}-$ on $\mathrm{Mod}_{(\mathcal{A},\mathcal{M})}^{\mathrm{cond}}$ making the functor $-\underset{{\mathcal{A}}}{\otimes}(\mathcal{A},\mathcal{M})$ symmetric monoidal.
    
    \item The functor \[ \mathcal{D}( \mathrm{Mod}_{(\mathcal{A},\mathcal{M})}^{\mathrm{cond}} )\rightarrow \mathcal{D}(\mathrm{Mod}^{\mathrm{cond}}_{{\mathcal{A}}})\] 
    is fully faithful, and its essential image is stable under all limits and colimits and given by those $C\in \mathcal{D}(\mathrm{Mod}_\mathcal{A}^{\mathrm{cond}})$ for which the map \[\mathrm{R}{\mathrm{Hom}}_{{\mathcal{A}}}(\mathcal{M}[S], C) \rightarrow \mathrm{R}{\mathrm{Hom}}_{{\mathcal{A}}}({\mathcal{A}}[S], C)\] is an isomorphism for all extremally disconnected sets $S$. In that case, the map 
    \[\mathrm{R}\underline{\mathrm{Hom}}_{{\mathcal{A}}}(\mathcal{M}[S], C) \rightarrow \mathrm{R}\underline{\mathrm{Hom}}_{{\mathcal{A}}}({\mathcal{A}}[S], C)\] 
    is also an isomorphism. An object $C \in \mathcal{D}(\mathrm{Mod}_{{\mathcal{A}}}^{\mathrm{cond}})$ lies in $\mathcal{D}(\mathrm{Mod}_{(\mathcal{A},\mathcal{M})}^{\mathrm{cond}})$ if and only if, for each $i\in \mathbb{Z}$, the $i$-th cohomology group $H^i(C)$ lies in $\mathrm{Mod}_{(\mathcal{A},\mathcal{M})}^{\mathrm{cond}}$. The inclusion functor $\mathcal{D}(\mathrm{Mod}_{(\mathcal{A},\mathcal{M})}^{\mathrm{cond}})\subset \mathcal{D}(\mathrm{Mod}^{\mathrm{cond}}_{{\mathcal{A}}})$ admits a left adjoint
    \[\mathcal{D}(\mathrm{Mod}_{{\mathcal{A}}}^{\mathrm{cond}})\rightarrow \mathcal{D}(\mathrm{Mod}^{\mathrm{cond}}_{(\mathcal{A},\mathcal{M})}),\ C\mapsto C\overset{L}{\underset{\mathcal{A}}{\otimes}}(\mathcal{A},\mathcal{M}),\] which is the left derived functor of $M\mapsto M\underset{\mathcal{A}}{\otimes}(\mathcal{A},\mathcal{M})$. Finally, there is a unique symmetric monoidal tensor product $-\underset{(\mathcal{A},\mathcal{M})}{\overset{L}{\otimes}}-$ on $\mathcal{D}(\mathrm{Mod}_{(\mathcal{A},\mathcal{M})}^{\mathrm{cond}})$ making the functor $-\overset{L}{\underset{{\mathcal{A}}}{\otimes}}(\mathcal{A},\mathcal{M})$ symmetric monoidal.
\end{enumerate} 
\end{theorem}

In what follows, we write $\mathcal{D}(\mathcal{A})$ and $\mathcal{D}(\mathcal{A},\mathcal{M})$ as shorthand for $\mathcal{D}(\mathrm{Mod}_{\mathcal{A}}^{\mathrm{cond}})$ and $\mathcal{D}(\mathrm{Mod}_{(\mathcal{A},\mathcal{M})}^{\mathrm{cond}})$, respectively, and refer to $\mathcal{D}(\mathcal{A},\mathcal{M})$ as the \textit{derived category of $(\mathcal{A},\mathcal{M})$-modules}. We also call an object of $\mathrm{Mod}_{{\mathcal{A}}}^{\mathrm{cond}}$ (resp. $\mathcal{D}(\mathcal{A})$) \textit{$(\mathcal{A},\mathcal{M})$-complete} if and only if it lies in $\mathrm{Mod}_{(\mathcal{A},\mathcal{M})}^{\mathrm{cond}}$ (resp. $\mathcal{D}(\mathcal{A},\mathcal{M})$). Finally, both $-\underset{{\mathcal{A}}}{\otimes}(\mathcal{A},\mathcal{M})$ and $-\underset{\mathcal{A}}{\overset{L}{\otimes}}(\mathcal{A},\mathcal{M})$ are referred to as the \textit{completion functor of $(\mathcal{A},\mathcal{M})$}.

Let us now briefly discuss morphisms of (pre-)analytic rings, tensor products, and functors of measures. Pre-analytic and analytic rings can be organized into a category as follows:

\begin{definition}
\label{morphisms}
Let $(\mathcal{A},\mathcal{M})$ and $(\mathcal{B},\mathcal{N})$ be pre-analytic rings. A map $f:(\mathcal{A},\mathcal{M})\rightarrow (\mathcal{B},\mathcal{N})$ is a map ${\mathcal{A}}\rightarrow {\mathcal{B}}$ of condensed rings together with a natural transformation $\Phi^{(\mathcal{A},\mathcal{M})}_{(\mathcal{B},\mathcal{N})}:\mathcal{M}[S]\rightarrow \mathcal{N}[S]$ of condensed $\mathcal{A}$-modules.

Let $(\mathcal{A},\mathcal{M})$ and $(\mathcal{B},\mathcal{N})$ be analytic rings. A map $f:(\mathcal{A},\mathcal{M})\rightarrow (\mathcal{B},\mathcal{N})$ is a map ${\mathcal{A}}\rightarrow {\mathcal{B}}$ of condensed rings such that for every extremally disconnected set $S$, the condensed module $\mathcal{N}[S]$ is $(\mathcal{A},\mathcal{M})$-complete.
\end{definition}

It is readily verified that every map of analytic rings is also a map of pre-analytic rings in a natural way. The converse may not hold in general. However, under some mild assumptions on $(\mathcal{A},\mathcal{M})$ and $(\mathcal{B},\mathcal{N})$, it does. For example, it is true if $\mathcal{A}$ is discrete or the condensed module $\mathcal{N}[\ast]$ is $(\mathcal{A},\mathcal{M})$-complete (see \cite[Proposition 7.14]{Condensed}). All analytic rings considered in this paper are covered by these two cases.

Let us record the following basic result concerning base change of analytic rings:
\begin{proposition}[\cite{Condensed}, Proposition 7.7]
Let $(\mathcal{A},\mathcal{M})$ and $(\mathcal{B},\mathcal{N})$ be analytic rings and $f:(\mathcal{A},\mathcal{M})\rightarrow (\mathcal{B},\mathcal{N})$ a morphism of analytic rings.
\begin{enumerate}
    \item There exists a unique functor $\mathrm{Mod}^{\mathrm{cond}}_{(\mathcal{A},\mathcal{M})}\overset{}{\rightarrow}\mathrm{Mod}^{\mathrm{cond}}_{(\mathcal{B},\mathcal{N})}$, denoted $-\underset{{(\mathcal{A},\mathcal{M})}}{\otimes}(\mathcal{B},\mathcal{N})$, such that the composite functor
    \[\mathrm{Mod}^{\mathrm{cond}}_{\mathcal{A}}\xrightarrow{-\underset{\mathcal{A}}{\otimes}{\mathcal{B}}}\mathrm{Mod}^{\mathrm{cond}}_{{\mathcal{B}}}\xrightarrow{-\underset{{\mathcal{B}}}{\otimes}(\mathcal{B},\mathcal{N})}\mathrm{Mod}^{\mathrm{cond}}_{(\mathcal{B},\mathcal{N})}\]
    factors as 
    \[(-\underset{(\mathcal{A},\mathcal{M})}{\otimes}(\mathcal{B},\mathcal{N}))\circ (-\underset{\mathcal{A}}{\otimes}(\mathcal{A},\mathcal{M})).\]
    \item The composite functor
    \[\mathcal{D}(\mathcal{A})\xrightarrow{-\underset{\mathcal{A}}{\overset{L}{\otimes}}\mathcal{B}}\mathcal{D}(\mathcal{B})\xrightarrow{-\underset{\mathcal{B}}{\overset{L}{\otimes}}(\mathcal{B},\mathcal{N})}\mathcal{D}(\mathcal{B},\mathcal{N})\]
    factors as
    \[(-\underset{(\mathcal{A},\mathcal{M})}{\overset{L}{\otimes}}(\mathcal{B},\mathcal{N}))\circ (-\underset{\mathcal{A}}{\overset{L}{\otimes}}(\mathcal{A},\mathcal{M})),\]
    where $-\underset{(\mathcal{A},\mathcal{M})}{\overset{L}{\otimes}}(\mathcal{B},\mathcal{N})$ denotes the left derived functor of $-\underset{{(\mathcal{A},\mathcal{M})}}{\otimes}(\mathcal{B},\mathcal{N})$.
\end{enumerate}
\end{proposition}

Now let us explain how to extend the functor of measures of an analytic ring $(\mathcal{A},\mathcal{M})$ to a functor \[\{\mathrm{profinite\ sets}\}\rightarrow \mathcal{D}(\mathcal{A},\mathcal{M}).\]

\begin{definition}
Let $(\mathcal{A},\mathcal{M})$ be an analytic ring. Then for every profinite set $S$, $\mathcal{M}[S]\in \mathcal{D}(\mathcal{A},\mathcal{M})$ is defined as $\mathcal{A}[S]\underset{\mathcal{A}}{\overset{L}{\otimes}}(\mathcal{A},\mathcal{M})$.
\end{definition}

Explicitly, it may be computed as follows. Let $S$ be a profinite set and $S_{\bullet}\rightarrow S$ a hypercover such that for each $n\geq 0$, $S_n$ is extremally disconnected. Then $\mathcal{M}[S]$ is represented by the complex
\[\cdots\rightarrow \mathcal{M}[S_2]\rightarrow \mathcal{M}[S_1]\rightarrow \mathcal{M}[S_0]\rightarrow 0.\]
Recall that in Examples~\ref{examples}.\ref{1}, \ref{examples}.\ref{2}, \& \ref{examples}.\ref{3}, we already defined $\mathbb{Z}_{\blacksquare}[S],\ A_{\blacksquare}[S]$, and $\mathbb{Z}_{p,\blacksquare}$[S] for every profinite set $S$. However, as follows from \cite[Proposition 5.6]{Condensed}, the two definitions yield the same notion for $\mathbb{Z}$. Using the same argument, one can prove this for general finitely generated $\mathbb{Z}$-algebras. Finally, it follows from Theorem~\ref{freeobjects} that the same is true for $\mathbb{Z}_{p}$. Furthermore, for analytic rings over $\mathbb{Z}_{\blacksquare}$, the constructed functor enjoys the following remarkable property:

\begin{proposition}
\label{prprof}
Let $(\mathcal{A},\mathcal{M})$ be an analytic ring over $\mathbb{Z}_{\blacksquare}$ (e.g., $A_{\blacksquare}$ with $A$ discrete). Then for every profinite set $S$, the object $\mathcal{M}[S]:=\mathcal{A}[S]\underset{\mathcal{A}}{\overset{L}{\otimes}}(\mathcal{A},\mathcal{M})$ is concentrated in degree $0$, compact, and projective. Moreover, for any pair of profinite sets $S_1$ and $S_2$, the derived tensor product $\mathcal{M}[S_{1}]\underset{(\mathcal{A},\mathcal{M})}{\overset{L}{\otimes}} \mathcal{M}[S_2]$ is isomorphic to $\mathcal{M}[S_{1}\times S_{2}]$.
\end{proposition}
\begin{proof}
The case where $(\mathcal{A},\mathcal{M})$ is $\mathbb{Z}_{\blacksquare}$ follows from \cite[Corollary 5.5]{Condensed}, \cite[Proposition 5.6]{Condensed}, and \cite[Theorem 5.8(ii)]{Condensed}. We deduce the general case by a base change argument. To this end, we compute: 
\[\mathcal{A}[S]\underset{\mathcal{A}}{\overset{L}{\otimes}}(\mathcal{A},\mathcal{M})\cong {\mathbb{Z}}[S]\underset{{\mathbb{Z}}}{\overset{L}{\otimes}}(\mathcal{A},\mathcal{M})\cong {\mathbb{Z}}_{\blacksquare}[S]\underset{{\mathbb{Z}}_{\blacksquare}}{\overset{L}{\otimes}}(\mathcal{A},\mathcal{M})\cong {\mathbb{Z}}_{\blacksquare}[S]\underset{{\mathbb{Z}}_{\blacksquare}}{\overset{}{\otimes}}(\mathcal{A},\mathcal{M}),\]
where the last isomorphism follows from the fact that ${\mathbb{Z}}_{\blacksquare}[S]$ is projective. As for every pair of profinite sets $S_1$ and $S_2$, the tensor product $\mathbb{Z}_{\blacksquare}[S_{1}]\underset{\mathbb{Z}_{\blacksquare}}{\otimes} \mathbb{Z}_{\blacksquare}[S_2]$ is isomorphic to $\mathbb{Z}_{\blacksquare}[S_{1}\times S_{2}]$, this computation also implies that $\mathcal{M}[S_{1}]\underset{(\mathcal{A},\mathcal{M})}{\overset{L}{\otimes}} \mathcal{M}[S_2]$ is isomorphic to $\mathcal{M}[S_{1}\times S_{2}]$. Since the functor $-\underset{\mathbb{Z}_{\blacksquare}}{\otimes}(\mathcal{A},\mathcal{M}): \mathrm{Mod}^{\mathrm{cond}}_{\mathbb{Z}_{\blacksquare}}\rightarrow \mathrm{Mod}^{\mathrm{cond}}_{(\mathcal{A},\mathcal{M})}$ is left adjoint to the forgetful functor, which preserves epimorphisms, we conclude that it preserves projective objects. Finally, compactness follows from the fact that evaluation on $S$ commutes with filtered colimits (see {\cite[\href{https://stacks.math.columbia.edu/tag/0739}{Tag 0739}]{stacks-project}}).
\end{proof}

\begin{corollary}
\label{prpr}
Let $(\mathcal{A},\mathcal{M})$ be an analytic ring over $\mathbb{Z}_{\blacksquare}$. Then the tensor product of two projective condensed modules over $(\mathcal{A},\mathcal{M})$ is projective.
\end{corollary}

\begin{corollary}
\label{derivedtens}
Let $(\mathcal{A},\mathcal{M})$ be an analytic ring over $\mathbb{Z}_{\blacksquare}$. Then the symmetric monoidal tensor product $-\underset{(\mathcal{A},\mathcal{M})}{\overset{L}{\otimes}}-$ of $\mathcal{D}(\mathcal{A},\mathcal{M})$ is the left derived functor of $-\underset{(\mathcal{A},\mathcal{M})}{\overset{}{\otimes}}-$.
\end{corollary}

We now introduce a generalization of the notion of an analytic ring, namely the notion of an \textit{animated analytic ring}, and briefly explain how to extend the previous discussion to that context. The main advantage of this generalization is that it provides a better category-theoretic framework. For example, the usual concept of pushout does not produce anything meaningful in the context of analytic rings and one genuinely needs to apply the machinery of \textit{animation} in order to obtain a well-behaved notion. Therefore, we need to recall the following basic definition first:

\begin{definition}[\cite{Analytic}, Definition 11.4]
Let $\mathcal{C}$ be a category that admits all small colimits and is generated under small colimits by the full subcategory $\mathcal{C}^{\mathrm{cp}}$ of compact projective objects. The \textit{animation} of $\mathcal{C}$ is the $\infty$-category $\mathrm{Ani}(\mathcal{C})$ freely generated under sifted colimits by $\mathcal{C}^{\mathrm{cp}}$.
\end{definition}

For the category of sets, this definition yields the usual $\infty$-category of ``spaces'' (or ``anima''). The case of primary interest in this paper is that where $\mathcal{C}$ is either the category of condensed sets or condensed abelian groups.

\begin{definition}
Let $\mathcal{C}$ be the category of condensed sets or condensed abelian groups, respectively. Then the category $\mathrm{Ani}(\mathcal{C})$ is called the $\infty$-category of animated condensed sets or animated condensed abelian groups, respectively. 
\end{definition}

By the Dold-Kan correspondence, the category of animated condensed abelian groups is equivalent to $\infty$-derived category of condensed abelian groups $\mathcal{D}_{\geq 0}(\mathrm{Cond(Ab)})$. For every animated condensed ring $\mathcal{A}$, the $\infty$-category of animated condensed $\mathcal{A}$-modules $\mathcal{D}_{\geq 0}(\mathcal{A})$ is defined in the obvious way. All the definitions and statements above are also straightforwardly adapted to the animated context. Let us do that explicitly for the notion of analytic ring:

\begin{definition}[\cite{Analytic} Definitions 12.1]
An \textit{analytic animated ring} $(\mathcal{A},\mathcal{M})$ is an animated condensed ring $\mathcal{A}$ (called the \textit{underlying animated condensed ring of $(\mathcal{A},\mathcal{M})$}) equipped with a functor 
\[\{\mathrm{extremally\ disconnected\ sets}\}\rightarrow \mathcal{D}_{\geq 0}(\mathcal{A}),\ S\mapsto \mathcal{M}[S],\]
called the \textit{functor of measures of $(\mathcal{A},\mathcal{M})$}, that sends finite disjoint unions to products, and a natural transformation $\underline{S}\rightarrow\mathcal{M}[S]$ of animated condensed sets, such that for every object $C\in \mathcal{D}_{\geq 0}(\mathcal{A})$ that is a sifted colimit of objects of the form $\mathcal{M}[T]$ for various extremally disconnected sets $T$, the map $\mathrm{R}\underline{\mathrm{Hom}}_{\mathcal{A}}(\mathcal{M}[S], C) \rightarrow \mathrm{R}\underline{\mathrm{Hom}}_{\mathcal{A}}({\mathcal{A}}[S], C)$ is an isomorphism for all extremally disconnected sets $S$. 

An analytic ring $(\mathcal{A},\mathcal{M})$ is called \textit{complete}, if the map $\mathcal{A}\rightarrow\mathcal{M}[\ast]$ is an isomorphism.
\end{definition}

Observe that the only essential difference with Definition~\ref{analytic ring} is that $\mathcal{A}$ and $\mathcal{M}[S]$, with $S$ extremally disconnected, are allowed to be connective complexes, i.e., those whose positive cohomological degrees are all trivial. An analytic ring is said to be $0$-truncated, if both $\mathcal{A}$ and $\mathcal{M}[S]$ are concentrated in degree $0$ (note that in that case, only the language has changed).

We now record the following analog of Theorem~\ref{wichtig} in the present context:

\begin{theorem}[\cite{Analytic}, Proposition 12.4]
Let $(\mathcal{A},\mathcal{M})$ be an analytic animated ring. Let $\mathcal{D}_{\geq 0}(\mathcal{A},\mathcal{M})$ denote the full $\infty$-subcategory of $\mathcal{D}_{\geq 0}(\mathcal{A})$ of all objects $M\in \mathcal{D}_{\geq 0}(\mathcal{A})$, such that for every extremally disconnected set $S$, the map \[{\mathrm{Hom}}_{\mathcal{A}}(\mathcal{M}[S], M)\rightarrow{\mathrm{Hom}}_{\mathcal{A}}(\mathcal{A}[S], M)\] is an isomorphism. The category $\mathcal{D}_{\geq 0}(\mathcal{A},\mathcal{M})$ is generated under sifted colimits by objects of the form $\mathcal{M}[T]$ for various extremally disconnected sets $T$. Furthermore, for extremally disconnected set $S$, the object $\mathcal{M}[S]$ is compact and projective in $\mathcal{D}_{\geq 0}(\mathcal{A},\mathcal{M})$. The $\infty$-subcategory $\mathcal{D}_{\geq 0}(\mathcal{A},\mathcal{M})\subset \mathcal{D}_{\geq 0}(\mathcal{A})$ is stable under all limits and colimits and the inclusion functor $\mathcal{D}_{\geq 0}(\mathcal{A},\mathcal{M})\subset \mathcal{D}_{\geq 0}(\mathcal{A})$ admits a left adjoint
\[\mathcal{D}_{\geq 0}(\mathcal{A})\rightarrow \mathcal{D}_{\geq 0}(\mathcal{A},\mathcal{M}),\ C\mapsto C\overset{L}{\underset{\mathcal{A}}{\otimes}}(\mathcal{A},\mathcal{M}),\] which sends $\mathcal{A}[S]$ to $\mathcal{M}[S]$. The $\infty$-category $\mathcal{D}_{\geq 0}(\mathcal{A},\mathcal{M})$ is prestable. Its heart is the abelian category $\mathcal{D}^{\heartsuit}(\mathcal{A},\mathcal{M})$ that is the full subcategory of condensed $\pi_0\mathcal{A}$-modules\footnote{The functor $\pi_0: \mathrm{Ani}(\mathrm{Cond})\rightarrow \mathrm{Cond}$ is given by the universal property
of $\mathrm{Ani}(\mathrm{Cond})$.} generated under colimits by $\pi_0\mathcal{M}[T]$ for various extremally disconnected sets $T$. An object $C \in \mathcal{D}_{\geq 0}(\mathcal{A})$ lies in $\mathcal{D}_{\geq 0}(\mathcal{A},\mathcal{M})$ if and only if for each $i\in \mathbb{Z}_{\leq 0}$, the $i$-th cohomology group $H^i(C)$ lies in $\mathcal{D}^{\heartsuit}(\mathcal{A},\mathcal{M})$. Finally, there is a unique symmetric monoidal tensor product $-\underset{(\mathcal{A},\mathcal{M})}{\overset{L}{\otimes}}-$ on $\mathcal{D}_{\geq 0}(\mathcal{A},\mathcal{M})$ making the functor $-\overset{L}{\underset{{\mathcal{A}}}{\otimes}}(\mathcal{A},\mathcal{M})$ symmetric monoidal.
\end{theorem}

Passing to ``spectrum objects''\footnote{A formal procedure that recovers $\mathcal{D}(-)$ from $\mathcal{D}_{\geq 0}(-)$ in the current situation, see \cite[Chapter 1.4.2]{HA}.}, we obtain the stable $\infty$-category $\mathcal{D}(\mathcal{A})$ from $\mathcal{D}_{\geq 0}(\mathcal{A})$; analogous statements hold for a similarly defined full $\infty$-subcategory $\mathcal{D}(\mathcal{A},\mathcal{M})\subset \mathcal{D}(\mathcal{A})$. The statements above about morphisms and base change of analytic rings are easily adapted to this framework. As we will always work with animated analytic rings throughout the remainder of the paper, or rather the distinction, as a rule (although not always), will not be important, we will also slightly abuse terminology from now on by generally omitting the word ``animated''.

One advantage of the generalized notion of an analytic ring is that it allows us to consider the induced analytic ring structure for any map of condensed rings, as made precise by the following:

\begin{proposition}[\cite{Analytic}, Proposition 12.8\footnote{The second half of the theorem is not stated explicitly in \cite{Analytic}; however, it readily follows from the proof.}]
\label{minimal0}
Let $(\mathcal{A},\mathcal{M})$ be an analytic ring and $\phi:{\mathcal{A}}\rightarrow\mathcal{B}$ a map of condensed rings. Then the functor
\[\{\mathrm{extremally\ disconnected\ sets}\}\rightarrow \mathcal{D}(\mathcal{B}),\ S\mapsto \mathcal{N}[S]:=\mathcal{B}[S] \underset{\mathcal{A}}{\overset{L}{\otimes}}(\mathcal{A},\mathcal{M})\]
defines an analytic ring structure on $\mathcal{B}$, which is said to be \textit{induced from $(\mathcal{A},\mathcal{M})$}, such that there is a (necessarily unique) map of analytic rings $\varphi:(\mathcal{A},\mathcal{M})\rightarrow (\mathcal{B},\mathcal{N})$ for which $\phi$ serves as the map of the underlying rings. Furthermore, an object $C$ in $\mathcal{D}(\mathcal{B})$ lies in $\mathcal{D}(\mathcal{B},\mathcal{N})$ if and only if its restriction to $\mathcal{D}(\mathcal{A})$ lies in $\mathcal{D}(\mathcal{A},\mathcal{M})$. In particular, the map $\varphi$ is initial among all maps of analytic rings of the form $(\mathcal{A},\mathcal{M})\rightarrow (\mathcal{B},\mathcal{N}')$.
\end{proposition}

Before we proceed, we need to mention that there is actually a subtlety in the definition of a commutative analytic ring. Recall that all rings, including condensed ones, are assumed commutative. However, this alone is not sufficient for a good notion of commutativity. As explained in \cite[Lectute XII]{Analytic}, a few other properties are desirable. They are all implied by the following definition: 

\begin{definition}[\cite{Analytic}, Definition 12.10]
A commutative analytic ring $(\mathcal{A},\mathcal{M})$ is an analytic ring $(\mathcal{A},\mathcal{M})$ such that, for every prime $p\in \mathbb{Z}$, the Frobenius map $\varphi_p : \mathcal{A} \rightarrow \mathcal{A}/^{\mathbb{L}}p$ induces a map of analytic rings $\varphi_p : (\mathcal{A},\mathcal{M})\rightarrow (\mathcal{A}/^{\mathbb{L}}p,\mathcal{M}/^{\mathbb{L}}p)$.
\end{definition}

This extra condition, however, is automatic in a number of situations. For example, it is satisfied by all 0-truncated analytic rings or analytic rings over $\mathbb{Z}_{\blacksquare}$ (see \cite[Appendix to Lecture XII: Frobenius]{Analytic}), which covers all examples encountered in this paper. Applying Proposition~\ref{minimal0} to the map $\mathcal{A}\rightarrow \mathcal{M}[\ast]$\footnote{Actually, it is one instance where the extra property in the definition of a commutative analytic ring is essential; see \cite[Appendix to Lecture XII]{Analytic} for details.}, where $(\mathcal{A},\mathcal{M})$ is an analytic ring, we obtain \textit{the completion of $(\mathcal{A},\mathcal{M})$}, which is uniquely characterized by the obvious universal property. It readily follows from the second half of that proposition that completion does not change the category of complete objects.

Let $\mathrm{AnRing}$ denote the category of complete commutative animated analytic rings with morphisms given by maps of animated analytic rings. We conclude the section by discussing pushouts of complete commutative analytic rings and steady maps, which will be essential to some of our proofs and constructions in the following sections.

\begin{proposition}[\cite{Analytic}, Proposition 12.12]
\label{pushout}
Let $(\mathcal{B},\mathcal{N}) \leftarrow (\mathcal{A},\mathcal{M}) \rightarrow (\mathcal{C},\mathcal{L})$ be a diagram of complete analytic commutative rings. The pushout $(\mathcal{E},\mathcal{K})\overset{\mathrm{def}}{=}(\mathcal{B},\mathcal{N})\underset{(\mathcal{A},\mathcal{M})}{\overset{L}{\otimes}} (\mathcal{C},\mathcal{L})$ of this diagram in $\mathrm{AnRing}$ exists. It can be defined as the completion of an analytic ring structure on $\mathcal{B}\underset{\mathcal{A}}{\overset{L}{\otimes}}\mathcal{C}$.
The corresponding functor
\[\mathcal{D}(\mathcal{E},\mathcal{K})\rightarrow \mathcal{D}(\mathcal{B}\underset{\mathcal{A}}{\overset{L}{\otimes}}\mathcal{C})\]
is fully faithful, with essential image given by all $C\in \mathcal{D}(\mathcal{B}\underset{\mathcal{A}}{\overset{L}{\otimes}}\mathcal{C})$ whose images in $\mathcal{D}(\mathcal{B})$ and $\mathcal{D}(\mathcal{C})$ lie in $\mathcal{D}(\mathcal{B},\mathcal{N})$ and $\mathcal{D}(\mathcal{C},\mathcal{L})$, respectively.
The left adjoint \[-\underset{\mathcal{B}\underset{\mathcal{A}}{\overset{L}{\otimes}}\mathcal{C}}{\overset{L}{\otimes}}(\mathcal{E},\mathcal{K}):\mathcal{D}(\mathcal{B}\underset{\mathcal{A}}{\overset{L}{\otimes}}\mathcal{C}) \rightarrow \mathcal{D}(\mathcal{E},\mathcal{K})\]
is given by the sequential colimit
\[-\rightarrow -\underset{\mathcal{B}}{\overset{L}{\otimes}}(\mathcal{B},\mathcal{N})\rightarrow(-\underset{\mathcal{B}}{\overset{L}{\otimes}}(\mathcal{B},\mathcal{N}))\underset{\mathcal{C}}{\overset{L}{\otimes}}(\mathcal{C},\mathcal{L})\rightarrow((-\underset{\mathcal{B}}{\overset{L}{\otimes}}(\mathcal{B},\mathcal{N}))\underset{\mathcal{C}}{\overset{L}{\otimes}}(\mathcal{C},\mathcal{L}))\underset{\mathcal{B}}{\overset{L}{\otimes}}(\mathcal{B},\mathcal{N})\rightarrow\cdots\]
where each functor is regarded as an endofunctor of $\mathcal{D}(\mathcal{B}\underset{\mathcal{A}}{\overset{L}{\otimes}}\mathcal{C})$.
\end{proposition}

As is obvious from the proposition, pushouts are somewhat subtle in general. However, they can be computed explicitly if one of the maps is \textit{steady}: 

\begin{definition}\cite[Definition 12.13]{Analytic}
A map $f:(\mathcal{A},\mathcal{M})\rightarrow (\mathcal{B},\mathcal{N})$ in $\mathrm{AnRing}$ is steady if for all maps $g:(\mathcal{A},\mathcal{M})\rightarrow (\mathcal{C},\mathcal{L})$ in $\mathrm{AnRing}$, the functor $M\mapsto M\underset{\mathcal{B}}{\overset{L}{\otimes}}(\mathcal{B},\mathcal{M})$ preserves the full $\infty$-subcategory of $\mathcal{D}(\mathcal{B}\underset{\mathcal{A}}{\overset{L}{\otimes}}\mathcal{C})$ spanned by objects whose image in $\mathcal{D}(\mathcal{C})$ lies in $\mathcal{D}(\mathcal{C},\mathcal{L})$. Equivalently, for all $M\in  \mathcal{D}(\mathcal{C},\mathcal{L})$, the object $M\underset{\mathcal{A}}{\overset{L}{\otimes}}(\mathcal{B},\mathcal{N})=M\underset{(\mathcal{A},\mathcal{M})}{\overset{L}{\otimes}}(\mathcal{B},\mathcal{N})$, regarded as an object in $\mathcal{D}(\mathcal{B}\underset{\mathcal{A}}{\overset{L}{\otimes}}\mathcal{C})$, lies in $\mathcal{D}(\mathcal{C},\mathcal{L})$ when restricted to $\mathcal{C}$.
\end{definition}

In fact, the condition that $f$ is steady means precisely that for all $g$ the colimit
\[-\rightarrow -\underset{\mathcal{C}}{\overset{L}{\otimes}}(\mathcal{C},\mathcal{L})\rightarrow(-\underset{\mathcal{C}}{\overset{L}{\otimes}}(\mathcal{C},\mathcal{L}))\underset{\mathcal{B}}{\overset{L}{\otimes}}(\mathcal{B},\mathcal{N})\rightarrow((-\underset{\mathcal{C}}{\overset{L}{\otimes}}(\mathcal{C},\mathcal{L}))\underset{\mathcal{B}}{\overset{L}{\otimes}}(\mathcal{B},\mathcal{N}))\underset{\mathcal{C}}{\overset{L}{\otimes}}(\mathcal{C},\mathcal{L})\rightarrow \cdots\]
stabilizes at $(-\underset{\mathcal{C}}{\overset{L}{\otimes}}(\mathcal{C},\mathcal{L}))\underset{\mathcal{B}}{\overset{L}{\otimes}}(\mathcal{B},\mathcal{N})$. In this case, the pushout $(\mathcal{E},\mathcal{K})$ of $(\mathcal{B},\mathcal{N}) \leftarrow (\mathcal{A},\mathcal{M}) \rightarrow (\mathcal{C},\mathcal{L})$ is given by the functor $\mathcal{K}: S \mapsto \mathcal{L}[S] \underset{(\mathcal{A},\mathcal{M})}{\overset{L}{\otimes}}(\mathcal{B},\mathcal{N})$. Yet another equivalent characterization is given by the following proposition.

\begin{proposition}[\cite{Analytic}, Proposition 12.14]
\label{useful completion}
A map $f: (\mathcal{A},\mathcal{M})\rightarrow (\mathcal{B},\mathcal{N})$ in $\mathrm{AnRing}$ is steady if and only if for all
pushout diagrams

\begin{center}
\begin{tikzcd}
(\mathcal{A},\mathcal{M})\arrow{d}\arrow{r} &   (\mathcal{B},\mathcal{N}) \arrow{d}\\
(\mathcal{C},\mathcal{L}) \arrow{r} & (\mathcal{E},\mathcal{K})
\end{tikzcd}
\end{center}
of complete commutative analytic rings and all $M \in \mathcal{D}(\mathcal{C},\mathcal{L})$, the base change map \[(M|_{\mathcal{A}})\underset{(\mathcal{A},\mathcal{M})}{\overset{L}{\otimes}}(\mathcal{B},\mathcal{N})\rightarrow M\underset{(\mathcal{C},\mathcal{L})}{\overset{L}{\otimes}}(\mathcal{E},\mathcal{K})\] is an isomorphism.
\end{proposition}
 
The class of steady maps enjoys the following properties:
 
\begin{proposition}[\cite{Analytic}, Proposition 12.15]\label{steady meow}
The class of steady maps is stable under base change and composition. If a steady map $f:(\mathcal{A},\mathcal{M})\rightarrow (\mathcal{B},\mathcal{N})$ in $\mathrm{AnRing}$ factors over a map $f':(\mathcal{A}',\mathcal{M}')\rightarrow (\mathcal{B},\mathcal{N})$, then $f'$ is also steady. Moreover, the class of steady maps is closed under all colimits.
\end{proposition}
\section{Analytic rings associated to complete Huber pairs}\label{general}
In this section, we establish some technical lemmas and develop the tools necessary for the main proofs. First, we introduce certain technical notions that are used later to describe explicitly the spaces of measures of the analytic ring associated to a complete Huber pair. Next, we collect some technical facts about analytic rings associated to certain discrete Huber pairs, which we will use later to prove the promised descent results. In the last subsection, we first define analytic rings associated to discrete Huber pairs and establish their basic properties, after that we construct an analytic ring for every complete Huber pair and provide a description of its functor of measures, and then we prove that the obtained functor from the category of complete Huber pairs to the category of complete commutative analytic rings is fully faithful.  

There is, however, one general statement about prodiscrete abelian groups that does not fit into any of the described subsections, but, nonetheless, is important for our analysis of analytic rings associated to complete Huber pairs. Therefore, we record it here for later use:
\begin{lemma}
\label{surj}
Let $g: M\rightarrow N$ be an open surjective continuous map of complete first countable linearly topologized abelian groups (in particular, both $M$ and $N$ are prodiscrete). Then the induced map $g_S:\underline{M}(S)\rightarrow\underline{N}(S)$ is surjective for any compact Hausdorff topological space $S$. 
\end{lemma}

\begin{proof}
Suppose that $f:S\rightarrow N$ is a continuous map and let and let $(M_n)_{n\in\mathbb{N}}$ and $(N_n)_{n\in\mathbb{N}}$ be fundamental systems of open submodules of $M$ and $N$ respectively. Observe that any locally constant map from $S$ to $N$ factors through a finite quotient of $S$. Therefore, every locally constant map trivially admits a lift to $M$. Using the same observation again, it is easy to check that $f$ can be written as an infinite sum of locally constant maps $\underset{i\geq 0}{\sum}f_i$ such that ${\forall n\geq 0 \ f_n(S)\subset N_n}$. Therefore, it evidently suffices to lift the individual maps $f_i$'s in such a way that the resulting sum still converges. For this, we need the following property: let $m$ be the maximal natural number such that $g^{-1}(N_n)\subset M_m$; then $m$ goes to infinity as $n$ does. This, however, easily follows from the fact that $M\rightarrow N$ is open and $(N_n)_{n\in\mathbb{N}}$ forms a basis of neighborhoods of $0$.
\end{proof}

We also record here the following classical result for later reference:

\begin{theorem}[Open mapping theorem]
\label{open mapping}
Let $A$ be an analytic complete Huber ring and $M,N$ complete first countable topological $A$-modules. Suppose that $f: M \rightarrow N$ is a continuous surjective $A$-module map. Then $f$ is open.
\end{theorem}

\begin{proof}
See \cite[Theorem 1.1.9]{Kdl}
\end{proof}
\subsection{Quasi-finitely generated submodules}\label{general1}
The object of this subsection is to study the additive structure of complete Huber rings from the point of view of topological algebra. We introduce the notion of \textit{quasi-finitely generated submodule}, which we will use later to give an explicit description of the spaces of measures $(A,A^+)_{\blacksquare}[S]$ of the analytic ring associated to a complete Huber pair. 
\begin{lemma}
\label{finitelygeneratedproto}
Let $A$ be a complete Huber ring, $(A_0,I)$ a pair of definition, and $B\subset A_0$ a subring (considered with the discrete topology). Consider a tower $\{M_n\}_{n\in \mathbb{N}}$ of finitely generated {$B$-modules} such that $M_n\subset A/I^n$, the image of $M_n$ under the canonical projection ${A/I^{n}\rightarrow A/I^{n-1}}$ lies inside $M_{n-1}$, and the induced map $M_n\rightarrow M_{n-1}$ is surjective. Then $\varprojlim {M_n}$ is a bounded closed subset of $A$ and, for any profinite set $S$, $\underline{A}(S)\cong \underset{(M_n)_{n\in\mathbb{N}}}{\colim}(\varprojlim\underline{M_n}(S))$, where the colimit is taken over all such towers.
\end{lemma}

\begin{proof}
We first observe that $M:=\varprojlim{M_n}$ is bounded. Indeed, by definition, we have to show that, for any $k\geq 0$, there exists $l\geq 0$ such that $I^l M\subset I^k$. This holds because $M_n$ is finitely generated for every $n\geq 0$ and $B$ is bounded (because $B\subset A_0$). Moreover, $M$ is closed, as it is the projective limit of a tower of closed subspaces.

Now let us show that $\underline{A}(S)\cong \underset{(M_n)_{n\in\mathbb{N}}}{\colim}(\varprojlim\underline{M_n}(S))$. This statement is equivalent to the fact that any continuous map $S\rightarrow A\cong \varprojlim A/I^n$ factors through $\varprojlim M_n$ for some tower $(M_n)_{n\in \mathbb{N}}$ as above. For every $n\geq 0$, consider the $B$-submodule generated by the image of $S$ in $A/I^n$ (which is finite as $S$ is compact and $A/I^n$ is discrete). It is easy to see that the constructed tower satisfies the desired properties and the image of $S$ in $A$ lies in its inverse limit.
\end{proof}

Note that this construction does depend on the choice of $B$. It also seems to depend on the choice of a pair of definition containing $B$. However, it does not: if $M\subset A$ comes from a tower as in Lemma \ref{finitelygeneratedproto} for some pair of definition $(A_0,I)$ containing $B$, then it actually comes from such a tower for any other pair of definition containing $B$. Indeed, if $(A'_0,I')$ is another such pair of definition, then it is straightforward to check that $M$ comes from the tower $\{M'_n\}_{n\in \mathbb{N}},\ M'_n=M/(M\cap (I')^n)$, which satisfies all the aforementioned properties. Therefore, we have the following well-behaved notion:

\begin{definition}
\label{quasifin}
Let $A$ be a complete Huber ring and $B$ a finitely generated subring of $A^{\circ}$. We say that a $B$-submodule $M=\varprojlim {M_n}$ of $A$ is \textit{quasi-finitely generated over $B$} if it comes from a tower as in Lemma \ref{finitelygeneratedproto} for some (and, hence, for any) pair of definition containing $B$.
\end{definition}

\begin{remark}
As $B$ is a finitely generated subring of $A^{\circ}$, there always exists a pair of definition containing $B$, so the previous definition makes sense.
\end{remark}

\begin{theorem}
\label{finitelygenerated}
Let $A$ be a complete Huber ring and $B\subset A^{\circ}$ a finitely generated subring. Then $\underline{A}\cong \underset{M}{\colim}\underline{M}$, where the colimit is taken over all quasi-finitely generated $B$-submodules $M\subset A$.
\end{theorem}

\begin{proof}
The statement immediately follows from Lemma \ref{finitelygeneratedproto}.
\end{proof}

\begin{remark}
It is easy to see that the collection $\{M_B  | \ {B\subset A^{\circ},\ M_B\subset A}\}$ forms a directed poset, where $B$ is a finitely generated $\mathbb{Z}$-subalgebra and $M_B$ is a quasi-finitely generated $B$-submodules of $A$.
\end{remark}

\begin{lemma}
\label{subquasi}
Let $A$ be a complete Huber pair and $B$ a finitely generated subring of $A^{\circ}$. Then every closed submodule of a quasi-finitely generated $B$-submodule $M\subset A$ is quasi-finitely generated over $B$.
\end{lemma}

\begin{proof}
Fix a pair of definition $(A_0,I)$ containing $B$ and let $M'\subset M$ be a closed submodule. Consider the image of $M'$ under the canonical map $M\rightarrow M_n$ and denote it $M_n'$. As the ring $B$ is Noetherian and $M_n$ is finitely generated over $B$, we see that $M_n'$ is finitely generated over $B$. It is clear from the definition that the canonical maps $M_{n+1}'\rightarrow M_n'$ are surjective. Thus, we obtain a quasi-finitely generated $B$-submodule $\Tilde{M}:=\varprojlim M_n'$. We claim that $M'=\Tilde{M}$. Indeed, it is clear from the definition that $M'$ is a dense subset of $\Tilde{M}$. As $M'$ is closed, it is equal to $\Tilde{M}$, hence quasi-finitely generated.
\end{proof}

In fact, if we allow infinite sums, a quasi-finitely generated module admits a countable system of generators, which is made precise by the following:

\begin{lemma}
\label{countbas}
Let $A$ be a complete Huber ring, $B\subset A^{\circ}$ a finitely generated subring, and $M\subset A$ a quasi-finitely generated submodule over $B$. Fix a pair of definition $(A_0,I)$ containing $B$ and, for every $n\geq 0$, denote the quotient $M/(M\cap I^n)$ by $M_n$.
\begin{enumerate}[(i)]
    \item There exists a family of finite sets $\{K_n\}_{n\in\mathbb{N}}$, where each $K_n\subset A/I^n$, such that

\begin{itemize}
    \item every element of $K_n$ is non-zero;
    \item ${K}_n$ generates $M_n$ over $B$;
    \item for every $n> 0$, ${K}_n$ is a disjoint union of ${K}_n'$ and ${K}_n''$ such that ${K}_n'$ is mapped one to one to ${K}_{n-1}$ and ${K}_n''$ generates $(M_n\cap I^{n-1})/I^n$ over $B$.
\end{itemize}

\item Every such family $(K_n)_{n\in \mathbb{N}}$ yields a countable ``set of generators'' of $M$. For an element $x\in K_n''$, consider the unique $x_m\in K_{n+m}$ that is mapped to $x$. The resulting sequence converges to an element of $M$. Applying this procedure to every $x\in K_n''$ for all $n>0$, we obtain an at most countable set $K\subset M$ with the following properties:

\begin{itemize}
    \item the image of $K$ in $M_n$ is exactly $K_n$; in addition, $K$ can be written as the disjoint union of $(K^n)'$ and $(K^n)''$ such that $(K^n)'$ is mapped bijectively onto $K_n$ and $(K^n)''$ is mapped to zero;
    \item every element of $M$ is an (infinite) sum of elements of $K$ with coefficients from $B$; in addition, every such sum converges and defines an element of $M$. 
 \end{itemize}
\end{enumerate}
\end{lemma}

\begin{proof}
First, observe that the second part of the lemma trivially follows from the first one. Thus, we need only construct, for each $n$, a finite set of $B$-generators of $M_n$ with the aforementioned special properties. Let $K_0$ be an arbitrary finite set of non-zero elements generating $M_0$ over $B$. Take an arbitrary lift of $K_0$ to $M_1$, which we denote by $K_1'$, and consider the $B$-submodule generated by the lifts, which we denote by $\tilde{M}_1$. It is immediate that every element of $M_1$ lies inside $I^0/I=A_0/I$ modulo $\tilde{M}_1$. As $B$ is Noetherian, the submodule $M_1\cap A_0/I$ is also finitely generated. Let $K_1''$ be an arbitrary finite set of its non-zero generators and define $K_1$ as the union $K_1'\cup K_1''$. Repeating this procedure inductively, we obtain a system of finite sets ${K}_n\subset M_n$ with the desired properties.
\end{proof}
\subsection{Technical lemmas}\label{general2}
In this subsection, we collect some technical lemmas, which will form a basis for our computations in what follows. The reader is advised to skip it on first reading, coming back as needed later.

\begin{lemma}
\label{steady maps}
The maps of analytic rings

\begin{equation}
    \mathbb{Z}_{\blacksquare}\rightarrow (\mathbb{Z}[T],\mathbb{Z})_{\blacksquare}\tag{1}
\end{equation}
\begin{equation}
    (\mathbb{Z}[T],\mathbb{Z})_{\blacksquare}\rightarrow \mathbb{Z}[T]_{\blacksquare} \tag{2}
\end{equation}
\begin{equation}
    \mathbb{Z}_{\blacksquare}\rightarrow \mathbb{Z}[T]_{\blacksquare} \tag{3}
\end{equation}
\begin{equation}
    (\mathbb{Z}[T],\mathbb{Z})_{\blacksquare}\rightarrow (\mathbb{Z}[T,T^{-1}],\mathbb{Z})_{\blacksquare} \tag{4}
\end{equation}
\begin{equation}
    (\mathbb{Z}[T,T^{-1}],\mathbb{Z})_{\blacksquare}\rightarrow (\mathbb{Z}[T,T^{-1}],\mathbb{Z}[T^{-1}])_{\blacksquare} \tag{5}
\end{equation}
\begin{equation}
    (\mathbb{Z}[T],\mathbb{Z})_{\blacksquare}\rightarrow (\mathbb{Z}[T,T^{-1}],\mathbb{Z}[T^{-1}])_{\blacksquare} \tag{6}
\end{equation}
are steady. 
\end{lemma}

\begin{proof}
The map (1) is steady by \cite[Proposition 13.14]{Analytic} (see Definition~\ref{definition nuclear} below and the discussion afterwards for a reminder on nuclearity) and (2) is steady by \cite[Example 13.15(2)]{Analytic}. The map (3) is steady by virtue of Proposition~\ref{steady meow} as it is the composition of (1) and (2). Invoking \cite[Proposition 13.14]{Analytic} once again, we deduce that the map $(4)$ is steady. Using it one more time together with (2), we deduce (5). Finally, (6) follows from (4) and (5).
\end{proof}

\begin{lemma}
\label{polynomialization}
For every $M\in \mathcal{D}(\mathbb{Z}_{\blacksquare})$, the object $M\underset{\mathbb{Z}_{\blacksquare}}{\overset{L}{\otimes}}(\mathbb{Z}[T],\mathbb{Z})_{\blacksquare}$ is isomorphic to $M\underset{\mathbb{Z}}{\overset{L}{\otimes}}\mathbb{Z}[T]$. Moreover, it can be explicitly given by $N^{\bullet}\underset{\mathbb{Z}}{\overset{}{\otimes}}\mathbb{Z}[T]$, where $N^{\bullet}$ is any complex representing $M$.
\end{lemma}

\begin{proof}
Note that the second part of the statement follows from the first as $-\underset{\mathbb{Z}}{\overset{L}{\otimes}}\mathbb{Z}[T]$ is exact (since arbitrary direct sums are exact); the first part itself is a special case of Lemma~\ref{minimal0}.
\end{proof}

\begin{lemma}
\label{localizationT}
For every $M$ in $\mathcal{D}((\mathbb{Z}[T],\mathbb{Z})_{\blacksquare})$ (resp. $\mathcal{D}(\mathbb{Z}[T]_{\blacksquare})$), the object $M\underset{(\mathbb{Z}[T],\mathbb{Z})_{\blacksquare}}{\overset{L}{\otimes}}(\mathbb{Z}[T,T^{-1}],\mathbb{Z})_{\blacksquare}$ (resp. $M\underset{\mathbb{Z}[T]_{\blacksquare}}{\overset{L}{\otimes}}(\mathbb{Z}[T,T^{-1}],\mathbb{Z}[T])_{\blacksquare}$) is isomorphic to $M\underset{\mathbb{Z}[T]}{\overset{L}{\otimes}} \mathbb{Z}[T,T^{-1}]$. Moreover, it can be explicitly given by $N^{\bullet}\underset{\mathbb{Z}[T]}{\overset{}{\otimes}}\mathbb{Z}[T,T^{-1}]$, where $N^{\bullet}$ is any complex representing $M$.
\end{lemma}

\begin{proof}
Arguing as in Lemma~\ref{polynomialization} and using Lemma~\ref{minimal0} again, we see that it suffices to show that the functor $-\underset{\mathbb{Z}[T]}{\overset{L}{\otimes}}\mathbb{Z}[T,T^{-1}]$ is exact and preserves $\mathrm{Mod}^{\mathrm{cond}}_{\mathbb{Z}[T]_{\blacksquare}}$. Both claims readily follow from the fact that, for any $L\in \mathrm{Mod}^{\mathrm{cond}}_{\mathbb{Z}[T]}$, $L\underset{\mathbb{Z}[T]}{\overset{}{\otimes}}\mathbb{Z}[T,T^{-1}]$ is isomorphic to the cokernel of $L[U]\xrightarrow[]{TU-1}L[U]$.
\end{proof}

\begin{lemma}
\label{powerseries}
Let $S$ be the profinite set $\mathbb{N}\cup \{\infty \}$ and $T$ a formal variable. Consider the condensed module $\mathbb{Z}_{\blacksquare}[S]/\mathbb{Z}$, where the inclusion $\mathbb{Z}\hookrightarrow \mathbb{Z}[\mathbb{N}\cup \{\infty\}]$ is induced by the map $\{*\}\rightarrow \mathbb{N}\cup \{\infty\},\ *\mapsto \infty$. Then it is compact and projective in $\mathrm{Mod}^{\mathrm{cond}}_{\mathbb{Z}_{\blacksquare}}$, and isomorphic to $\mathbb{Z}[[T]]$ via the unique morphism $\mathbb{Z}[S]\mapsto \mathbb{Z}[[T]]$ corresponding to the continuous mapping $S\rightarrow \mathbb{Z}[[U]]$ given by the formula $n\mapsto T^n,\infty\mapsto 0$.
\end{lemma}
\begin{proof}
For each $n\geq 0$, consider the following map of condensed abelian groups: 
\[\underset{i=0}{\overset{n}{\bigoplus}}\mathbb{Z}\xrightarrow{f_n} \underset{i=0}{\overset{n+1}{\bigoplus}}\mathbb{Z}, (a_0,\dots,a_n)\mapsto (a_0,\dots,a_n,a_n).\] It is readily verified that the colimit of $(\underset{i=0}{\overset{n}{\bigoplus}}\mathbb{Z},f_n)_{n\in\mathbb{N}}$ is isomorphic to the space of continuous functions $\underline{C(S,\mathbb{Z})}$. For every $n\geq 0$, let $\phi_n:\underset{i=0}{\overset{n}{\bigoplus}}\mathbb{Z}\rightarrow \mathbb{Z}$ denote the map that sends $(a_0,\dots,a_n)$ to $a_n$. The collection $(\phi_n)_{n\in \mathbb{N}}$ clearly defines a map $\phi$ from the colimit of $\{\underset{i=0}{\overset{n}{\bigoplus}}\mathbb{Z},f_n\}_{n\in\mathbb{N}}$ to $\mathbb{Z}$. Under the identifications of the colimit with $\underline{C(S,\mathbb{Z})}$ and $\mathbb{Z}$ with $\underline{C(*,\mathbb{Z})}$, the map $\phi$ is induced by the map $\{*\}\rightarrow S$ given by $*\mapsto {\infty}$. Observe that for each $n\geq 0$, the kernel of $\phi_n$ is $\underset{i=0}{\overset{n-1}{\bigoplus}}\mathbb{Z}$, and the induced map $\ker \phi_n \rightarrow \ker \phi_{n+1}$ is the canonical inclusion $\underset{i=0}{\overset{n-1}{\bigoplus}}\mathbb{Z}\rightarrow\underset{i=0}{\overset{n}{\bigoplus}}\mathbb{Z}$. Therefore, we clearly have \[\underset{\mathbb{N}}{\bigoplus}\mathbb{Z}\cong \ker(\underline{C(S,\mathbb{Z})}\rightarrow \underline{C(*,\mathbb{Z})}).\] As $\mathbb{Z}$ is projective, $\underline{C(S,\mathbb{Z})}$ is isomorphic to $\underset{\mathbb{N}}{\bigoplus}\mathbb{Z} \oplus \mathbb{Z}$, and hence we have
\[\mathbb{Z}[[T]]\cong \underset{\mathbb{N}}{\prod}\mathbb{Z}\cong\underline{\mathrm{Hom}}(\underset{\mathbb{N}}{\bigoplus}\mathbb{Z},\mathbb{Z})\cong \underline{\mathrm{Hom}}(\underline{C(S,{\mathbb{Z}})},\mathbb{Z})/\mathbb{Z}\cong \mathbb{Z}_{\blacksquare}[S]/\mathbb{Z}.\]
\end{proof}

\begin{proposition}
\label{completion derived}
Let $M$ be an object of $\mathcal{D}((\mathbb{Z}[T],\mathbb{Z})_{\blacksquare})$. Then the completion $M\underset{(\mathbb{Z}[T],\mathbb{Z})_{\blacksquare}}{\overset{L}{\otimes}}\mathbb{Z}[T]_{\blacksquare}$ can be given explicitly by $\mathrm{R}\underline{\mathrm{Hom}}_{R}(R_{\infty}/R,M)[1]$, where $R=\mathbb{Z}[T],\ R_{\infty}=\mathbb{Z}((T^{-1}))$.
\end{proposition}

\begin{proof}
The statement is proved in \cite[Lecture VIII]{Condensed}.
\end{proof}

\begin{proposition}
\label{completion derived 2}
Let $M$ be an object of $\mathcal{D}(\mathbb{Z}_{\blacksquare})$. Then the completion $M\underset{\mathbb{Z}_{\blacksquare}}{\overset{L}{\otimes}}\mathbb{Z}[T]_{\blacksquare}$ can be given explicitly by $\mathrm{R}\underline{\mathrm{Hom}}_{\mathbb{Z}}(R_{\infty}/R,M)$, where $R=\mathbb{Z}[T],\ R_{\infty}=\mathbb{Z}((T^{-1}))$.
\end{proposition}

\begin{proof}
We note that the object $R_{\infty}/R$ is in fact isomorphic to $\mathbb{Z}[[T^{-1}]]$ as a $\mathbb{Z}_{\blacksquare}$-module. As the latter is compact (by Lemma~\ref{powerseries}), it suffices to prove the statement for objects of the form $\mathbb{Z}_{\blacksquare}[S]$ with extremally disconnected $S$. As such a module is isomorphic to $\prod_J \mathbb{Z}$ for some set $J$, we need only show that we have 
\[\mathrm{R}\underline{\mathrm{Hom}}_{\mathbb{Z}}(R_{\infty}/R,\mathbb{Z})\cong \mathbb{Z}[T],\]
which follows from Lemma~\ref{powerseries}.
\end{proof}

We conclude this subsection by the following result:

\begin{proposition}
\label{completion derived topological}
Let $M$ be a prodiscrete topological abelian group. Then $M\underset{\mathbb{Z}_{\blacksquare}}{\overset{L}{\otimes}}\mathbb{Z}[T]_{\blacksquare}$ is isomorphic to $\underline{M\langle T\rangle}$, where $M\langle T\rangle$ denotes the abelian group of formal power series $\underset{i\geq 0}{\sum} m_iT^i$ with coefficients in $M$ such that $(m_i)_{i\in\mathbb{N}}$ converges to zero (endowed with the compact-open topology).
\end{proposition}

\begin{proof}
This is an immediate consequence of Proposition~\ref{completion derived 2}. Indeed, we write $M$ as the projective limit $\varprojlim M_i$ of an inverse system of discrete modules. It can be easily checked that the module $M\langle T \rangle$ is isomorphic to $\varprojlim M_i[T]$. Therefore, as the functor $-\underset{\mathbb{Z}_{\blacksquare}}{\overset{L}{\otimes}}\mathbb{Z}[T]_{\blacksquare}$ commutes with projective limits, it suffices to check that the completion $M_i\underset{\mathbb{Z}_{\blacksquare}}{\overset{L}{\otimes}}\mathbb{Z}[T]_{\blacksquare}$ is isomorphic to $M_i[T]$, which is clear.
\end{proof}
\subsection{Analytic rings associated to complete Huber pairs and their functors of measures}\label{general3}
The main goal of this subsection is to construct a $0$-truncated analytic ring for every complete Huber pair and give an explicit description of its functor of measures. In fact, we will construct an analytic ring for every pair $(A,A^+)$, where $A$ is a complete Huber ring and $A^+$ is an arbitrary subring of $A^{\circ}$. However, this level of generality is rather redundant for our purposes: as explained at the end of this section, for every such fake ``Huber pair'' $(A,A^+)$, our construction will produce the same analytic ring as for $(A,\widetilde{A^+})$, where $\widetilde{A^+}$ is the minimal open integrally closed subring containing $A^+$. Before discussing the general case, we first examine in detail the construction of the analytic ring associated to a discrete Huber pair, see \cite[Definition 9.1.]{Condensed}. In particular, we give a complete proof that such a condensed analytic ring is indeed analytic presupposing this only for finitely generated $\mathbb{Z}$-algebras (see \cite[Theorem 8.13]{Condensed} for a proof of the latter). Recall that for a finitely generated $\mathbb{Z}$-algebra $A$, its associated analytic ring is denoted by $A_{\blacksquare}$ (see Example~\ref{examples}.\ref{2}). Its construction can be generalized to arbitrary discrete $\mathbb{Z}$-algebras as follows.

\begin{definition}
Let $A$ be a discrete $\mathbb{Z}$-algebra. The pre-analytic ring $A_{\blacksquare}$ is defined as the pre-analytic ring with underlying condensed ring $A$ and the functor of measures given by
\[\{\mathrm{extremally\ disconnected\ sets}\}\rightarrow \mathrm{Mod}^{\mathrm{cond}}_{A},\ S\mapsto A_{\blacksquare}[S]:=\colim_{A'\rightarrow A}A'_{\blacksquare}[S],\]
where the colimit is taken over all finitely generated $\mathbb{Z}$-subalgebras $A'$ of $A$.
\end{definition}

\begin{proposition}\label{general_discrete}
Let $A$ be a discrete $\mathbb{Z}$-algebra. Then the pre-analytic ring $A_{\blacksquare}$ is analytic. Furthermore, an object $M$ in $\mathrm{Mod}^{\mathrm{cond}}_{A}$ (resp. $\mathcal{D}(A)$) is $A_{\blacksquare}$-complete if and only if it is $A'_{\blacksquare}$-complete for every finitely generated subring $A'$ of $A$.
\end{proposition}

\begin{proof}
We first observe that the second part of the statement follows straightforwardly from the construction of $A_{\blacksquare}$. Therefore, we need only show that $A_{\blacksquare}$ is analytic. Unwinding the definition, we see that it amounts to proving that for every pair of index sets $I,J$, families $\{S_i\}_{i\in I},\{S_j\}_{j\in J}$ of extremally disconnected sets and map
\[f: \underset{i\in I}{\bigoplus}{A_{\blacksquare}}[S_i]\rightarrow \underset{j\in J}{\bigoplus}{A_{\blacksquare}}[S_j]\]
with kernel $K$, the map \[\mathrm{R}\underline{\mathrm{Hom}}_{A}(A_{\blacksquare}[S], K) \rightarrow \mathrm{R}\underline{\mathrm{Hom}}_{A}(A[S], K)\] is an isomorphism for all extremally disconnected sets $S$. As the free module $A[S]$ can be written as $\colim_{A'\rightarrow A}A[S]$, it suffices to prove that for every finitely generated $\mathbb{Z}$-subalgebra $A'$ of $A$, the canonical map \[\mathrm{R}\underline{\mathrm{Hom}}_{A'}(A'_{\blacksquare}[S], K) \rightarrow \mathrm{R}\underline{\mathrm{Hom}}_{A'}(A'[S], K)\] is an isomorphism. Since for every finitely generated $\mathbb{Z}$-subalgebra $A''$ of $A$ containing $A'$ and every extremally disconnected set $T$, the module $A_{\blacksquare}''[T]$ is $A_{\blacksquare}'$-complete and $\mathcal{D}(A_{\blacksquare}')$ is stable under both colimits and limits, the desired claim follows.
\end{proof}

In order to generalize the definition above to discrete Huber pairs, we need to establish first some auxiliary technical results, which are also used later to construct analytic rings for general complete Huber pairs.

\begin{lemma}
\label{discrete}
Let $A$ be a discrete ring and $M$ a discrete module over $A$. Then the associated condensed module $\underline{M}$ is $A_{\blacksquare}$-complete.
\end{lemma}

\begin{proof}
Consider the étale site $\ast_{\mathrm{\acute{e}t}}$ of the point $\ast$, i.e., the category $\mathrm{Fin}$ of finite sets with covers given by finite families of jointly surjective maps. The category of sheaves on this site is canonically isomorphic to the category of sets. Observe that there is an obvious functor $\ast_{\mathrm{\acute{e}t}}\rightarrow \ast_{\mathrm{pro\acute{e}t}}$, which sends a finite set to itself. This functor is evidently continuous and preserves finite products. Therefore, we have a pullback functor from the category of $A$-modules to the category of condensed modules over $\underline{A}$, which sends a (discrete) module $M$ to $\underline{M}$. As this functor commutes with colimits and the category of $A_{\blacksquare}$-complete modules is closed under colimits and contains $\underline{A}$, the desired result follows.
\end{proof}

\begin{lemma}
\label{fingenmodtens}
Let $A$ be a finitely generated $\mathbb{Z}$-algebra and $M$ a discrete module over $A$. Then, for every set $I$, we have
\[\underline{M} \underset{\underline{A}}{\overset{}{\otimes}}\underset{I}{\prod}\underline{A}\cong \underline{M} \underset{\underline{A}}{\overset{L}{\otimes}}\underset{I}{\prod}\underline{A}\cong \underline{M} \underset{A_{\blacksquare}}{\overset{L}{\otimes}}\underset{I}{\prod}\underline{A}\cong \underline{M} \underset{A_{\blacksquare}}{\overset{}{\otimes}}\underset{I}{\prod}\underline{A}.\]
If $M$ is finitely generated, then, in fact, we have $\underline{M} \underset{A_{\blacksquare}}{\overset{}{\otimes}}\underset{I}{\prod}\underline{A}\cong \underset{I}{\prod}\underline{M}$.
\end{lemma}

\begin{proof}
For the first isomorphism, see the very beginning of \cite[Appendix to Lecture VIII]{Condensed}. It is also shown there that for every finitely generated module $M$, we have $\underline{M} \underset{\underline{A}}{\overset{}{\otimes}}\underset{I}{\prod}\underline{A}\cong\underset{I}{\prod}M$. Therefore, to prove the remaining two isomorphisms, we first reduce to that case by using that $M\underset{A_{\blacksquare}}{\overset{L}{\otimes}}\underset{I}{\prod}\underline{A}$ commutes with filtered colimits in $M$, and then simply note that $\underset{I}{\prod}M$ is complete.
\end{proof}

\begin{lemma}
\label{alternative}
Let $A^+\subset A$ be a pair of discrete rings. Then for every profinite set $S$, we have \[{A}[S]\underset{{A}^+}{\overset{L}{\otimes}}A^+_{\blacksquare}\cong{A}[S]\underset{{A}^+}{\overset{}{\otimes}}A^+_{\blacksquare}\cong {A}\underset{A^+_{\blacksquare}}{\otimes}A^+_{\blacksquare}[S]\cong {A}\underset{{A}^+}{\overset{}{\otimes}}A^+_{\blacksquare}[S].\]
\end{lemma}

\begin{proof}
Note that the condensed module $A$ is $A^+_{\blacksquare}$-complete by virtue of Lemma~\ref{discrete}. Therefore, we have
\[A[S]\underset{A^+}{\overset{L}{\otimes}}A^+_{\blacksquare}\cong (A\underset{A^+}{\overset{L}{\otimes}} A^+[S])\underset{A^+}{\overset{L}{\otimes}}A^+_{\blacksquare}\cong A\underset{A^+_{\blacksquare}}{\overset{L}{\otimes}}A^+_{\blacksquare}[S].\]
To finish the proof, we must show that $A\underset{A^+}{\overset{L}{\otimes}}A^+_{\blacksquare}[S]$ is concentrated in degree $0$ and $A^+_{\blacksquare}$-complete. Writing $A^+_{\blacksquare}[S]$ as $\colim_{A'\rightarrow A^+}A'_{\blacksquare}[S]$, where the colimit is taken over finitely generated $\mathbb{Z}$-subalgebras $A'$ of $A^+$, we see that $A\underset{A^+}{\overset{L}{\otimes}}A^+_{\blacksquare}[S]$ is isomorphic to $\colim_{A'\rightarrow A^+}A\underset{A'}{\overset{L}{\otimes}}A'_{\blacksquare}[S]$. Therefore, the desired claim follows readily from Lemma~\ref{fingenmodtens}.
\end{proof}

Combining Proposition~\ref{minimal0} and Lemma~\ref{alternative}, we obtain the following $0$-truncated analytic ring:

\begin{definition}
Let $A^+\subset A$ be a pair of discrete rings. The analytic ring $(A,A^+)_{\blacksquare}$ is defined as the underlying ring $A$ with analytic structure induced from $A^{+}_{\blacksquare}$ along the inclusion $A^+\rightarrow A$. In particular, its functor of measures is given by the functor
\[\{\mathrm{extremally\ disconnected\ sets}\}\rightarrow \mathrm{Mod}^{\mathrm{cond}}_{\underline{A}},\ S\mapsto (A,A^+)_{\blacksquare}[S]:=A\underset{A^+}{\overset{L}{\otimes}}A^+_{\blacksquare}[S].\]
\end{definition}

Before we proceed to the general case, let us record a few useful properties of the construction in the discrete case:

\begin{proposition}
\label{relative_completion}
Let $A^+\subset A$ be a pair of discrete rings. Then an object $M$ in $\mathrm{Mod}^{\mathrm{cond}}_{A}$ (resp. $\mathcal{D}(A)$) lies in $\mathrm{Mod}^{\mathrm{cond}}_{(A,A^+)_{\blacksquare}}$ (resp. $\mathcal{D}((A,A^+)_{\blacksquare})$) if and only if its restriction to $A^+$ lies in $\mathrm{Mod}^{\mathrm{cond}}_{A^+_{\blacksquare}}$ (resp. $\mathcal{D}(A^+_{\blacksquare})$).
\end{proposition}

\begin{proof}
This is a special case of Proposition~\ref{minimal0}.
\end{proof}

\begin{proposition}
\label{integral induced}
Let $f:A'\rightarrow A$ be an integral map of (discrete) rings. Then the analytic structure $A_{\blacksquare}$ on $A$ is induced from the structure $A'_{\blacksquare}$ on $A'$. In particular, an object $M$ in $\mathrm{Mod}^{\mathrm{cond}}_{A}$ (resp. $\mathcal{D}(A)$) lies in $\mathrm{Mod}^{\mathrm{cond}}_{A_{\blacksquare}}$ (resp. $\mathcal{D}(A_{\blacksquare})$) if and only if its restriction to $A'$ lies in $\mathrm{Mod}^{\mathrm{cond}}_{A'_{\blacksquare}}$ (resp. $\mathcal{D}(A'_{\blacksquare})$).
\end{proposition}

\begin{proof} 
Writing $A'$ as the union $\colim_{R'\subset A'} R'$ of its finitely generated $\mathbb{Z}$-subalgebras and considering the integral closures $\overline{f(R')}$ of their images in $A$, we reduce proof to the case where $A'$ is a finitely generated $\mathbb{Z}$-algebra. Indeed, it readily follows from the fact that any finitely generated subring of $A$ is contained in $\overline{f(R')}$ for some finitely generated subalgebra $R'\subset A'$. Applying the second part of Proposition~\ref{general_discrete}, we further reduce to the case where $f$ is finite. But in that case, the desired statement is a straightforward consequence of Lemma~\ref{fingenmodtens}.
\end{proof}

We note that Proposition~\ref{integral induced} implies that for every pair of discrete rings $A^+\subset A$, the analytic ring structure $(A,A^+)_{\blacksquare}$ on $A$ is isomorphic to $(A,\bar{A}^+)_{\blacksquare}$, where $\bar{A}^+$ denotes the integral closure of $A^+$ in $A$. Therefore, as in the context of adic spaces, it suffices to consider only discrete Huber pairs.

Let us now construct an analytic ring for every complete Huber pair. In what follows, let $A$ be a complete Huber ring and $A^+$ an arbitrary subring of $A^{\circ}$ and let $(\mathcal{A},\mathcal{M})$ denote the analytic ring $(A_{\mathrm{disc}},A^+_{\mathrm{disc}})_{\blacksquare}$ and $(\mathcal{A}^+,\mathcal{M}^+)$ the analytic ring $(A^+_{\mathrm{disc}})_{\blacksquare}$. 
\begin{proposition}
\label{projdisc}
For every profinite set $S$, $\mathcal{M}[S]$ (resp. $\mathcal{M}^+[S]$) is concentrated in degree $0$, compact, and projective. For every pair of profinite sets $S_1$ and $S_2$, the tensor product $\mathcal{M}[S_1]\underset{(\mathcal{A},\mathcal{M})}{\overset{L}{\otimes}}\mathcal{M}[S_2]$ (resp. $\mathcal{M}^+[S_1]\underset{(\mathcal{A}^+,\mathcal{M}^+)}{\overset{L}{\otimes}}\mathcal{M}^+[S_2]$) is isomorphic to $\mathcal{M}[S_1\times S_2]$ (resp. $\mathcal{M}^+[S_1\times S_2]$); therefore, the tensor product of projective objects is projective. Finally, the symmetric monoidal tensor product $-\underset{(\mathcal{A},\mathcal{M})}{\overset{L}{\otimes}}-$ (resp. $-\underset{(\mathcal{A}^+,\mathcal{M}^+)}{\overset{L}{\otimes}}-$) of $\mathcal{D}(\mathcal{A},\mathcal{M})$ (resp. $\mathcal{D}(\mathcal{A}^+,\mathcal{M}^+)$) is the left derived functor of $-\underset{(\mathcal{A},\mathcal{M})}{\overset{}{\otimes}}-$ (resp. $-\underset{(\mathcal{A}^+,\mathcal{M}^+)}{\overset{}{\otimes}}-$).
\end{proposition}

\begin{proof}
The statement follows immediately from Proposition~\ref{prprof}, Corollary~\ref{prpr} and Corollary~\ref{derivedtens}.
\end{proof}

\begin{lemma}
\label{going_back}
The condensed ring $\underline{A}$ is $(\mathcal{A},\mathcal{M})$- and $(\mathcal{A}^+,\mathcal{M}^+)$-complete.
\end{lemma}

\begin{proof}
By Proposition~\ref{relative_completion}, it suffices to show that $\underline{A}$ is $(\mathcal{A}^+,\mathcal{M}^+)$-complete. Observe that for every extremally disconnected set $S$, we have
\[\mathrm{Hom}_{\mathcal{A}^+}(\mathcal{M}[S], \underline{A})\cong\varprojlim\mathrm{Hom}_{B}(B_{\blacksquare}[S], \underline{A})\overset{}{\cong} \varprojlim\mathrm{Hom}_{B}(B[S], \underline{A})\cong \underline{A}(S),\]
where the limits are taken over all finitely generated subrings $B$ of $A^+$ (endowed with the discrete topology). Indeed, as $B$ is finitely generated, there exists a ring of definition of $A$ containing $B$. Therefore, we can write $A$ as an inverse limit of discrete $B$-modules, and thus $\underline{A}$ is $B_{\blacksquare}$-complete by Lemma~\ref{discrete}, which implies the desired isomorphism.
\end{proof}

Let $S$ be a profinite set. In what follows, we will denote the object $\underline{A}[S] \underset{\mathcal{A}^+}{\overset{L}{\otimes}} (\mathcal{A}^+,\mathcal{M}^+)\in \mathcal{D}(\underline{A})$ by $(A,A^+)_{\blacksquare}[S]$.

\begin{lemma}
\label{freeguyshuber}
For every profinite set $S$, we have 
\[(A,A^+)_{\blacksquare}[S]\cong  \underline{A}\underset{(\mathcal{A}^+,\mathcal{M}^+)}{\overset{}{\otimes}}\mathcal{M}^+[S]\cong \underline{A}\underset{(\mathcal{A},\mathcal{M})}{\overset{}{\otimes}}\mathcal{M}[S]\cong \underline{A}[S]\underset{A_{\mathrm{disc}}}{\overset{L}{\otimes}}(\mathcal{A},\mathcal{M}).\]

In particular, $(A,A^+)_{\blacksquare}[S]$ is concentrated in degree $0$.
\end{lemma}

\begin{proof}
Applying Lemma~\ref{going_back}, we obtain 
\[\underline{A}[S] \underset{\mathcal{A}^+}{\overset{L}{\otimes}} (\mathcal{A}^+,\mathcal{M}^+)\cong (\underline{A}\underset{\mathcal{A}^+}{\overset{L}{\otimes}}\mathcal{A}^+[S]) \underset{\mathcal{A}^+}{\overset{L}{\otimes}} (\mathcal{A}^+,\mathcal{M}^+)\cong \underline{A}\underset{(\mathcal{A}^+,\mathcal{M}^+)}{\overset{L}{\otimes}}\mathcal{M}^+[S].\]
A similar computation shows that $\underline{A}[S]\underset{\mathcal{A}}{\overset{L}{\otimes}}(\mathcal{A},\mathcal{M})\cong  \underline{A}\underset{(\mathcal{A},\mathcal{M})}{\overset{L}{\otimes}}\mathcal{M}[S]$. Using Lemma~\ref{alternative}, we now compute
\[\underline{A}\underset{(\mathcal{A}^+,\mathcal{M}^+)}{\overset{L}{\otimes}}\mathcal{M}^+[S]\cong (\underline{A}\underset{(\mathcal{A},\mathcal{M})}{\overset{L}{\otimes}}\mathcal{A})\underset{(\mathcal{A}^+,\mathcal{M}^+)}{\overset{L}{\otimes}}\mathcal{M}[S]\cong \underline{A}\underset{(\mathcal{A},\mathcal{M})}{\overset{L}{\otimes}}(\mathcal{A}\underset{(\mathcal{A}^+,\mathcal{M}^+)}{\overset{L}{\otimes}}\mathcal{M}^+[S]) \cong  \underline{A}\underset{(\mathcal{A},\mathcal{M})}{\overset{L}{\otimes}}\mathcal{M}[S].\]

Therefore, it remains to verify that $\underline{A}\underset{(\mathcal{A}^+,\mathcal{M}^+)}{\overset{L}{\otimes}}\mathcal{M}^+[S]$ is concentrated in degree $0$. By definition, we have $\mathcal{M}^+[S]\cong\colim_{B\rightarrow A^+} B_{\blacksquare}[S]$. A direct computation shows that ${\colim_{B\rightarrow A^+} (\underline{A}\underset{B_{\blacksquare}}{\overset{L}{\otimes}} B_{\blacksquare}[S])}$ is $(\mathcal{A}^+,\mathcal{M}^+)$-complete; hence, it is isomorphic to $\underline{A}\underset{(\mathcal{A}^+,\mathcal{M}^+)}{\overset{L}{\otimes}}\mathcal{M}^+[S]$. Therefore, it suffices to show that $\underline{A}\underset{B_{\blacksquare}}{\overset{L}{\otimes}} B_{\blacksquare}[S]$ is concentrated in degree $0$ for every (discrete) finitely generated subring $B\subset A^+$. In fact, we claim that it is isomorphic to $\colim_M \underset{J}{\prod} \underline{M}$, where $J$ is a set such that $\mathbb{Z}_{\blacksquare}[S]\cong \underset{J}{\prod} \mathbb{Z}$ and the colimit is taken over all quasi-finitely generated $B$-submodules of $A$. Applying Theorem~\ref{finitelygenerated}, we see that it is enough to show that $\underline{M}\underset{B_{\blacksquare}}{\overset{L}{\otimes}} \underset{J}{\prod} B \cong \underset{J}{\prod} \underline{M}$, where $M=\varprojlim M_n$ is a quasi-finitely generated $B$-submodule of $A$, which follows from Lemma~\ref{computation}.
\end{proof}

\begin{lemma}
\label{computation}
Let $B$ be a finitely generated $\mathbb{Z}$-algebra and $M=\varprojlim M_i$ the projective limit of an inverse system of finitely generated (discrete) $B$-modules. Then for every set $J$, we have
\[ \underline{M}\underset{B_{\blacksquare}}{\overset{L}{\otimes}} \underset{J}{\prod} B \cong \underset{J}{\prod} \underline{M}.\]
\end{lemma}

\begin{proof}
By Lemma~\ref{fingenmodtens}, we have
\[\varprojlim( \underline{M_n}\underset{B_{\blacksquare}}{\overset{L}{\otimes}} \underset{J}{\prod} B)\overset{}{\cong}\varprojlim( \underline{M_n}\underset{B_{\blacksquare}}{\overset{}{\otimes}} \underset{J}{\prod} B)\overset{}{\cong} \varprojlim \underset{J}{\prod} \underline{M_n}\cong \underset{J}{\prod} \underline{M}. \] Moreover, we have $\varprojlim M_n\cong \mathrm{R}\varprojlim M_n$ and $\varprojlim (\underline{M_n}\underset{B_{\blacksquare}}{\overset{L}{\otimes}} \underset{J}{\prod} B)\cong \mathrm{R}\varprojlim (\underline{M_n}\underset{B_{\blacksquare}}{\overset{L}{\otimes}} \underset{J}{\prod} B)$, since the maps in the towers are surjective. We claim that the canonical map \[(\varprojlim \underline{M_n})\underset{B_{\blacksquare}}{\overset{L}{\otimes}} \underset{J}{\prod} B\rightarrow\mathrm{R}\varprojlim( \underline{M_n}\underset{B_{\blacksquare}}{\overset{L}{\otimes}} \underset{J}{\prod} B)\] is in fact an isomorphism. As we have the following canonical morphism of triangles
\begin{center}
\begin{tikzcd}[]
(\mathrm{R}\varprojlim \underline{M_n}) \underset{B_{\blacksquare}}{\overset{L}{\otimes}} B_{\blacksquare}[S] \arrow{r}\arrow{d} & (\underset{\mathbb{N}}{\prod} \underline{M_n}) \underset{B_{\blacksquare}}{\overset{L}{\otimes}} B_{\blacksquare}[S] \arrow{r}\arrow{d} & (\underset{\mathbb{N}}{\prod} \underline{M_n})\underset{B_{\blacksquare}}{\overset{L}{\otimes}} B_{\blacksquare}[S] \arrow{r}\arrow{d} & (\mathrm{R}\varprojlim \underline{M_n}[1])\underset{B_{\blacksquare}}{\overset{L}{\otimes}} B_{\blacksquare}[S] \arrow{d} \\
\mathrm{R}\varprojlim (\underline{M_n} \underset{B_{\blacksquare}}{\overset{L}{\otimes}} B_{\blacksquare}[S])\arrow{r} & \underset{\mathbb{N}}{\prod} (\underline{M_n}\underset{B_{\blacksquare}}{\overset{L}{\otimes}} B_{\blacksquare}[S]) \arrow{r} & \underset{\mathbb{N}}{\prod} (\underline{M_n}\underset{B_{\blacksquare}}{\overset{L}{\otimes}} B_{\blacksquare}[S]) \arrow{r} & \mathrm{R}\varprojlim (\underline{M_n}\underset{B_{\blacksquare}}{\overset{L}{\otimes}} B_{\blacksquare}[S])[1]
\end{tikzcd}
\end{center}
it suffices to show that the two middle vertical maps are isomorphisms. As already noted, \[\underline{M_n}\underset{B_{\blacksquare}}{\overset{L}{\otimes}} B_{\blacksquare}[S]\cong \underset{J}{\prod} \underline{M_n}.\] Thus, we need to prove that 
\[(\underset{\mathbb{N}}{\prod} \underline{M_n}) \underset{B_{\blacksquare}}{\overset{L}{\otimes}} B_{\blacksquare}[S]\cong \underset{\mathbb{N},J}{\prod} \underline{M_n}.\] 
Suppose first that $M_n$ is finite free for every $n\geq 0$. In this special case, the claim is a consequence of a more general statement: for every set $I$, we have
\[\underset{I}{\prod} B \underset{B_{\blacksquare}}{\overset{L}{\otimes}} \underset{J}{\prod} B\cong \underset{I,J}{\prod} B.\]
This can be deduced easily by examining the proof of \cite[Proposition 6.3]{Condensed}. Now, for each $n\geq 0$, pick a projective resolution of $M_n$ by finite free modules and consider their product. Since arbitrary direct products are exact and products of free modules are projective in $\mathrm{Mod}_{B_{\blacksquare}}^{\mathrm{cond}}$, we get a projective resolution of $\underset{\mathbb{N}}{\prod} M_n$ which can be used to compute $(\underset{\mathbb{N}}{\prod} M_n) \underset{B_{\blacksquare}}{\overset{L}{\otimes}} B_{\blacksquare}[S]$. Therefore, we are readily reduced to the special case above, which implies the desired isomorphism, which implies the desired isomorphism.
\end{proof}

As a consequence of the proof, we have the following explicit description of $(A,A^+)_{\blacksquare}[S]$:

\begin{theorem}
\label{freeobjects}

For every profinite set $S$, we have \[(A,A^+)_{\blacksquare}[S]\cong \underset{B\rightarrow A^+,\ M}{\colim} \underset{J}{\prod} \underline{M},\] where $J$ is a set such that $\mathbb{Z}_{\blacksquare}[S]\cong \underset{J}{\prod} {\mathbb{Z}}$ and the colimit is taken over all finitely generated subrings $B\subset A^+$ and all quasi-finitely generated $B$-submodules $M$ of $A$.
\end{theorem}

\begin{remark}
General abstract nonsense yields a canonical injective map \[\underset{B\rightarrow A^+,\ M}{\colim} \underset{J}{\prod} \underline{M}\rightarrow \underset{J}{\prod} \underline{A}.\]
\end{remark}

Combining Proposition~\ref{minimal0} and Lemma~\ref{freeguyshuber}, we finally obtain the following result:

\begin{theorem}
The functor
\[\{\mathrm{extremally\ disconnected\ sets}\}\rightarrow \mathrm{Mod}^{\mathrm{cond}}_{\underline{A}},\ S\mapsto (A,A^+)_{\blacksquare}[S]\]
defines a complete analytic ring $(A,A^+)_{\blacksquare}$ with underlying condensed ring $\underline{A}$, such that there is a (necessarily unique) map of analytic rings $\varphi:(A_{\mathrm{disc}},A^+_{\mathrm{disc}})_{\blacksquare}\rightarrow (A,A^+)_{\blacksquare}$ for which the canonical map $\underline{A}_{\mathrm{disc}}\rightarrow\underline{A}$ serves as the map of the underlying condensed rings.
\end{theorem}

Using Lemma~\ref{freeguyshuber}, we see that the analytic structure of $(A,A^+)_{\blacksquare}$ can be described as the analytic structure induced from $(A_{\mathrm{disc}},A^+_{\mathrm{disc}})_{\blacksquare}$ or from $(A^+_{\mathrm{disc}})_{\blacksquare}$. Therefore, we the following useful consequence of Proposition~\ref{minimal0}:

\begin{proposition}
\label{completeness+str}
An object $M$ in $\mathrm{Mod}_{\underline{A}}^{\mathrm{cond}}$ (resp. $\mathcal{D}(\underline{A})$) lies in $\mathrm{Mod}_{(A,A^+)_{\blacksquare}}^{\mathrm{cond}}$ (resp. $\mathcal{D}((A,A^+)_{\blacksquare})$) if and only if it lies in $\mathrm{Mod}_{(A^+)_{\blacksquare}}^{\mathrm{cond}}$ (resp. $\mathcal{D}((A^+)_{\blacksquare})$) if and only if it lies in  $\mathrm{Mod}_{(A_{\mathrm{disc}},A^+_{\mathrm{disc}})_{\blacksquare}}^{\mathrm{cond}}$ (resp. $\mathcal{D}((A_{\mathrm{disc}},A^+_{\mathrm{disc}})_{\blacksquare})$) if and only if it lies in $\mathrm{Mod}_{(A^+_{\mathrm{disc}})_{\blacksquare}}^{\mathrm{cond}}$ (resp. $\mathcal{D}((A^+_{\mathrm{disc}})_{\blacksquare})$).
\end{proposition}

Our next goal is to prove a useful criterion of $(A,A^+)_{\blacksquare}$-completeness for condensed modules over $\underline{A}$ and then explain why it is sufficient to consider only the case where $(A,A^+)$ is a complete Huber pair. First, we introduce a bit of terminology.

\begin{definition}
Let $I$ be an index set and $f_I=\{f_i|i\in I\}$ a collection of elements of $A^{\circ}$. \textit{The ring of integral elements topologically generated by $f_I$} is defined as the minimal open integrally closed subring of $A$ containing $f_I$ .
\end{definition}

\begin{remark}
This notion is well-defined. Indeed, every open integrally closed subring contains the ring of topologically nilpotent elements $A^{\circ \circ}$, which is open as it contains every ideal of definition for every ring of definition.
\end{remark}

Before we formulate and prove the criterion, we prove the following technical lemma:

\begin{lemma}
\label{nilpotent completeness}
Let $M$ be an object in $\mathrm{Mod}^{\mathrm{cond}}_{\underline{A}}$ or $\mathcal{D}(\underline{A})$. If $M$ is $\mathbb{Z}_{\blacksquare}$-complete, then it is also $\mathbb{Z}[f]_{\blacksquare}$-complete for every topologically nilpotent element $f\in A^{\circ \circ}$.
\end{lemma}

\begin{proof}
We note that it is enough to prove the statement for $\mathcal{D}(\underline{A})$. Applying Proposition~\ref{integral induced}, we deduce that it suffices to show that $M$ is $\mathbb{Z}[T]_{\blacksquare}$-complete, where $M$ is endowed with the $\mathbb{Z}[T]$-module structure given by $T\mapsto f$. As $f$ is topologically nilpotent, we get a map \[\underline{\mathbb{Z}[[T]]}\rightarrow \underline{A},\ T\mapsto f\] (where ${\mathbb{Z}[[T]]}$ is endowed with the $T$-adic topology). Therefore, $M$ is a $\mathbb{Z}_\blacksquare$-complete $\underline{\mathbb{Z}[[T]]}$-module, hence a $(\mathbb{Z}[[T]],\mathbb{Z})_{\blacksquare}$-module. Thus, we need to prove that $(\mathbb{Z}[[T]],\mathbb{Z})_{\blacksquare}[S]$ is $\mathbb{Z}[T]_{\blacksquare}$-complete for every profinite set $S$. Using Theorem~\ref{freeobjects}, we obtain an isomorphism 
\[(\mathbb{Z}[[T]],\mathbb{Z})_{\blacksquare}[S]\cong \underset{J}{\prod}\mathbb{Z}[[T]].\]
Therefore, it is enough to show that $\mathbb{Z}[[T]]$ is $\mathbb{Z}[T]_{\blacksquare}$-complete. Writing it as the inverse limit $\varprojlim{\mathbb{Z}[T]/(T^n)}$, we see that it is sufficient to prove that $\mathbb{Z}[T]/(T^n)$ is $\mathbb{Z}[T]_{\blacksquare}$-complete, which clearly holds as $\mathbb{Z}[T]/(T^n)$ is discrete.
\end{proof}

\begin{proposition}
\label{generated and complete}
Let $I$ be an index set and $f_I=\{f_i|i\in I\}$ a collection of elements of $A^{\circ}$. Assume that $A^+$ is topologically generated by $f_i$ (in particular, we assume that $(A,A^+)$ is a complete Huber pair). Let $M$ be an object of $\mathcal{D}(\underline{A})$. It is $(A,A^+)_{\blacksquare}$-complete if and only if it is $\mathbb{Z}[f_i]_{\blacksquare}$-complete for every $i\in I$. 
\end{proposition}

\begin{proof}
Assume that $M$ is $(A,A^+)_{\blacksquare}$-complete; then the fact that it is also $\mathbb{Z}[f_i]_{\blacksquare}$-complete follows trivially from the discussion under Definition~\ref{morphisms}.

Conversely, assume that $M$ is $\mathbb{Z}[f_i]_{\blacksquare}$-complete for every $i\in I$. By Lemma~\ref{nilpotent completeness}, we may assume that $A^{\circ \circ}\subset f_I$. Thus, we may assume that $A^+$ is integral over its subring generated by $f_I$. Applying Proposition~\ref{completeness+str} and then Proposition~\ref{integral induced}, we see that it suffices to prove that $M$ is $A'_{\blacksquare}$-complete, where $A'\subset A^+$ is the subring generated by $f_I$. Using the second part of Proposition~\ref{general_discrete} and Proposition~\ref{integral induced} again, we reduce to the case where $A'=\mathbb{Z}[x_1,\dots,x_n]$, $I=\{1,\dots,n\}$, and $f_1=x_1,\dots, f_n=x_n$. In this case, we argue by induction on $n$. We observe that the desired statement is equivalent to the following isomorphism of analytic rings: 
\[\mathbb{Z}[x_1,\dots,x_{n-1}]_{\blacksquare}\underset{\mathbb{Z}_{\blacksquare}}{\overset{L}{\otimes}}\mathbb{Z}[x_n]_{\blacksquare}\cong\mathbb{Z}[x_1,\dots,x_{n}]_{\blacksquare}.\] 
Recall that the map $\mathbb{Z}_{\blacksquare}\rightarrow \mathbb{Z}[T]_{\blacksquare}$ is steady (by Lemma~\ref{steady maps}). Therefore, it suffices to show that the following holds for every set $J$:
\[(\underset{J}{\prod}\mathbb{Z}[x_1,\dots,x_{n-1}])\underset{\mathbb{Z}_{\blacksquare}}{\overset{L}{\otimes}}\mathbb{Z}[x_n]_{\blacksquare}\cong\underset{J}{\prod}\mathbb{Z}[x_1,\dots,x_{n}].\]
But this follows directly from Proposition~\ref{completion derived topological}.
\end{proof}

Proposition~\ref{generated and complete} implies that the analytic ring $(A,A^+)_{\blacksquare}$ is isomorphic to the one associated to $(A,\widetilde{A^+})$, where $\widetilde{A^+}$ is the minimal open integrally closed subring containing $A^+$. Thus, as in the discrete case, we may restrict our attention to complete Huber pairs. 

As a formal consequence of the proof above, we also obtain the following isomorphism of analytic rings:

\[\mathbb{Z}[x_1,\dots,x_n,y_1,\dots,y_m]_{\blacksquare}\underset{\mathbb{Z}[x_1,\dots,x_n]_{\blacksquare}}{\overset{L}{\otimes}}\mathbb{Z}[x_1,\dots,x_n,z_1,\dots,z_l]_{\blacksquare}\cong \mathbb{Z}[x_1,\dots,x_n,y_1,\dots,y_m,z_1,\dots,z_l]_{\blacksquare}.\]

In fact, even more is true. Let $S$ be a (discrete) ring and $R$ and $R'$ (discrete) rings over $S$ such that $R\underset{S}{\overset{L}{\otimes}} R'\cong R\underset{S}{\overset{}{\otimes}} R'$. Then, applying Proposition~\ref{generated and complete}, we obtain the following:

\begin{proposition}
We have the following isomorphism of analytic rings:
\[R_{\blacksquare}\underset{S_{\blacksquare}}{\overset{L}{\otimes}} R'_{\blacksquare}\cong (R\underset{S}{\overset{L}{\otimes}} R')_{\blacksquare}.\]
\end{proposition}

Applying Theorem~\ref{freeobjects}, we see that the association $(A,A^{+})\mapsto (A,A^{+})_{\blacksquare}$ in fact defines a functor from the category of complete Huber pairs to the category of complete commutative analytic rings. We close this section by proving that it is fully faithful.

\begin{proposition}
Let $\mathrm{cAff}$ denote the category of complete Huber pairs. The functor
\[\mathrm{cAff}\rightarrow \mathrm{AnRing},\ (A,A^+)\rightarrow (A,A^+)_{\blacksquare}\]
is fully faithful.
\end{proposition}
\begin{proof}
As complete Huber rings are compactly generated Hausdorff (in fact, metrizable) topological spaces, the functor is clearly faithful. Now let $(A,A^+)$ and $(B,B^+)$ be complete Huber pairs and let $f:A\rightarrow B$ be a continuous ring morphism. Suppose that for every extremally disconnected set $S$, the condensed module $(B,B^+)_{\blacksquare}[S]$ is $(A,A^+)_{\blacksquare}$-complete. We claim that $f(A^+)\subset B^+$. Let us assume the opposite: namely, that there exist an element $a\in A^+$ such that $b:=f(a)\not \in B^+$. Precomposing $f$ with the map \[(\mathbb{Z}[T],\mathbb{Z}[T])\rightarrow (A,A^+),\ T\mapsto a,\] we may assume that $(A,A^+)$ is $(\mathbb{Z}[T],\mathbb{Z}[T])$. Let $S$ denote the profinite set $\mathbb{N}\cup \{\infty \}$ and consider the map $\mathbb{Z}[T]_{\blacksquare}[S]\rightarrow (A,A^+)_{\blacksquare}[S]$. Using Theorem~\ref{freeobjects}, we see that it is explicitly given by
\[\underset{\mathbb{N}}{\prod}\mathbb{Z}[T]\rightarrow\underset{B\rightarrow A^+,\ M}{\colim} \underset{\mathbb{N}}{\prod} \underline{M},\  T\mapsto b.\]
Considering the image of the element $(1,T,T^2,\dots)$, we deduce that there exists a finitely generated subring $R\subset B^+$ and a quasi-finitely generated $R$-submodule $M\subset B$ such that the subring $\mathbb{Z}[b]$ lies inside $M$.  As we may clearly assume that $M$ contains $R$, we see that for each $n\geq 0$, the element $b$ is integral over $R$ modulo $I^n$, where $(B_0,I)$ is a pair of definition such that $R\subset B_0$. As $I$ is a subset of $B^+$, it implies that $b$ is integral over $B^+$, in contradiction to the fact that $B^+$ is integrally closed.
\end{proof}

\section{Analytic descent for animated solid modules on analytic adic spaces}\label{descentchapter}
The main goal of this section is to prove the following result:

\begin{theorem}
\label{maindescent}
Let $X$ be an analytic adic space and let $U$ denote an arbitrary affinoid subspace of $X$. Then the association $U\mapsto \mathcal{D}((\mathcal{O}_X(U),\mathcal{O}_X^+(U))_{\blacksquare})$ defines a sheaf of $\infty$-categories on $X$.
\end{theorem}

The proof of Theorem~\ref{maindescent} will follow the strategy of \cite[Lecture X]{Condensed} very closely. There is, however, one important difference. The main ingredients in the discrete case are certain statements about localization and base change (see \cite[Proposition 10.1]{Condensed} and \cite[Proposition 10.2]{Condensed}). To prove them, one argues by induction and reduces to the case where the rational open subspace under consideration is of a special form. As is, this approach does not work in the analytic case. The main problem is that it is required in the definition of a rational subspace $X(\frac{f_1,\dots, f_n}{g})$ that the elements $f_1,\dots,f_n,g$ generate an open ideal. This is, of course, always true in the discrete case but obviously wrong in general. Therefore, we need to modify the argument. We proceed as follows: we first introduce open covers of a certain special form that have a number of technical advantages and then verify the aforementioned statements in that case. To this end, we introduce the following notion:

\begin{definition}[\cite{Kdl}, Definition 1.6.6]
Let $X$ be an affinoid adic space and $f$ an element of $\mathcal{O}_X(X)$. The open covering $\{X(\frac{f}{1}),X(\frac{1}{f})\}$ of $X$ is called \textit{the simple Laurent covering defined by $f$} and the open covering $\{X(\frac{1}{f}),X(\frac{1}{1-f})\}$ is called \textit{the simple balanced covering defined by $f$}. A covering is called \textit{nice}, if it is a composition of coverings, each of which is either a simple Laurent covering or a simple balanced covering. A rational open subspace is called \textit{nice} if it is an element of some nice covering.
\end{definition}
\begin{remark}
\leavevmode
\begin{enumerate}[label=(\roman*)]
    \item We probably need to explain the word ``composition'' here. We start with a simple Laurent/balanced covering $U_1,U_2$; both of these subsets are affinoid. Therefore, we may form simple Laurent coverings $U_1=U_{11}\cup U_{12}$ and $U_2=U_{21}\cup U_{22}$, which correspond to some $f_1\in \mathcal{O}_X(U_1)$ and $f_2\in \mathcal{O}_X(U_2)$ respectively. In the same fashion, we construct Laurent coverings of $U_{11}$, $U_{12}$, $U_{21}$, and $U_{22}$, and so on. Any such covering obtained in a finite number of steps is called a \textit{composition} of Laurent coverings.
    \item Note that a covering all of whose elements are nice, is not necessarily nice.
\end{enumerate}
\end{remark}

Besides their simple form, the main technical advantage of nice open subsets is that they are sufficiently ``fine'' in the following sense: 

\begin{proposition}[\cite{Kdl}, Lemma 1.6.13 and Lemma 1.9.14] \label{kedlnice} Let $X$ be an analytic affinoid adic space. Then the family of its nice rational subsets forms a basis of quasi-compact open subsets, stable under intersections. Furthermore, every open covering of $X$ can be refined by a nice covering.
\end{proposition}

We will also need the following result of Kedlaya-Liu:

\begin{proposition}[\cite{Kdl}, Lemma 1.8.2]
\label{open kedl}
Let $(A,A^+)$ be a sheafy analytic complete Huber pair and $f,g$ elements of $A$ that generate an open ideal. Then multiplication by $(gT-f)$ defines a closed embedding of $A\langle T\rangle$ into itself.
\end{proposition}

\begin{proof}
Suppose first that either $f=1$ and $g\in A^{\circ}$ or $f\in A^{\circ}$ and $g=1$. In both of these cases, the statement follows from a straightforward computation. 

In the general case, we note that the elements $f$ and $g$ generate the unit ideal as the ring $A$ is sheafy. Therefore, $X:=\mathrm{Spa}(A,A^+)$ is covered by $X(\frac{f}{g})$ and $X(\frac{g}{f})$. Consider the commutative diagram
 \begin{center}
    \begin{tikzcd}
    & 0 \arrow[d, ] & 0 \arrow[d, ] & 0 \arrow[d, ] & \\
    0 \arrow[r,] & A\langle T\rangle  \arrow[r,"\cdot (gT-f)"]\arrow[d, ] & A\langle T\rangle  \arrow[r,]\arrow[d,] & A(\frac{f}{g}) \arrow[d, ]\arrow[r,] & 0\\
    0 \arrow[r,] & A(\frac{f}{g})\langle T\rangle \oplus A(\frac{g}{f})\langle T\rangle \arrow[r,"\cdot (gT-f)"]\arrow[d, ] & A(\frac{f}{g})\langle T\rangle \oplus A(\frac{g}{f})\langle T\rangle  \arrow[r,]\arrow[d,] & A(\frac{f}{g})\oplus A(\frac{f}{g},\frac{g}{f}) \arrow[d, ]\arrow[r,] & 0\\
    0 \arrow[r,] & A(\frac{f}{g},\frac{g}{f})\langle T\rangle  \arrow[r,"\cdot (gT-f)"]\arrow[d, ] & A(\frac{f}{g},\frac{g}{f})\langle T\rangle \arrow[r,] \arrow[d,]& A(\frac{f}{g},\frac{g}{f}) \arrow[d, ]\arrow[r,] & 0\\
    & 0 & 0 & 0 &
    \end{tikzcd}
\end{center}

Applying Theorem~\ref{open mapping}, we see that its columns are exact, since the ring $A$ is sheafy. Therefore, we are reduced to proving that its middle and bottom rows are exact, which follows from the cases handled above.
\end{proof}

Before we begin with the proof of Theorem~\ref{maindescent}, we need to establish some auxiliary results. We begin by recording the following consequence of Propositions~\ref{going_back}:

\begin{lemma}
Let $A$ be a complete Huber ring. Then $\underline{A}$ is $\mathbb{Z}_{\blacksquare}$-complete and $\underline{A\langle T \rangle}$ is \mbox{$\mathbb{Z}[T]_{\blacksquare}$-complete}. 
\end{lemma}

We also have the following consequence of Lemma~\ref{polynomialization}:

\begin{lemma}
The condensed ring $\underline{A}[T]$ is $(\mathbb{Z}[T],\mathbb{Z})_{\blacksquare}$-complete.
\end{lemma}

We now establish the following useful technical result:

\begin{lemma}
\label{trivial_bundle}
Let $A$ be a complete Huber ring. Then we have the following isomorphism of analytic rings: 
\[(A\langle T\rangle,A^+\langle T\rangle)_{\blacksquare} \cong (A,A^+)_{\blacksquare}\underset{\mathbb{Z}_{\blacksquare}}{\overset{L}{\otimes}}\mathbb{Z}[T]_{\blacksquare}.\]
In particular, for any profinite set $Q$, we have the following isomorphism in $\mathcal{D}(\underline{A}[T])$:
\[(A\langle T\rangle,A^+\langle T\rangle)_{\blacksquare}[Q] \cong (A,A^+)_{\blacksquare}[Q]\underset{\mathbb{Z}_{\blacksquare}}{\overset{L}{\otimes}}\mathbb{Z}[T]_{\blacksquare}.\]
Furthermore, the map \[(A,A^+)_{\blacksquare}\rightarrow (A\langle T\rangle,A^+\langle T\rangle)_{\blacksquare}\] is steady.
\end{lemma}

\begin{proof}
As the map $\mathbb{Z}_{\blacksquare}\rightarrow \mathbb{Z}[T]_{\blacksquare}$ is steady, both the second and the last parts of the statement follow from the first. Applying Proposition~\ref{completion derived topological}, we obtain an isomorphism
\[\underline{A\langle T \rangle}= \underline{A}\underset{\mathbb{Z}_{\blacksquare}}{\overset{L}{\otimes}}\mathbb{Z}[T]_{\blacksquare}.\]

To prove the lemma, it suffices to show that the categories $\mathcal{D}((A,A^+)_{\blacksquare}\underset{\mathbb{Z}_{\blacksquare}}{\overset{L}{\otimes}}\mathbb{Z}[T]_{\blacksquare})$ and $\mathcal{D}((A\langle T\rangle,A^+\langle T \rangle)_{\blacksquare})$, both of which are full subcategories of $\mathcal{D}(\underline{A\langle T\rangle})$ by the above isomorphism, are equal. Indeed, the desired claim will follow as the completion functors in both cases are given by the left adjoints of the inclusion functors. By Propositions~\ref{pushout} and \ref{completeness+str}, an object $M\in \mathcal{D}(\underline{A}\langle T\rangle)$ lies in $\mathcal{D}((A,A^+)_{\blacksquare}\underset{\mathbb{Z}_{\blacksquare}}{\overset{L}{\otimes}}\mathbb{Z}[T]_{\blacksquare})$ if and only if it is $(A^+)_{\blacksquare}$- and $\mathbb{Z}[T]$-complete. By Proposition~\ref{generated and complete} this is equivalent to being $A^+\langle T\rangle_{\blacksquare}$-complete; thus, we obtain the desired statement applying Proposition~\ref{completeness+str} one more time.
\end{proof}

We now recall the following definition:

\begin{definition}
Let $f:(\mathcal{A},\mathcal{M})\rightarrow (\mathcal{B},\mathcal{N})$ be a map of analytic rings. It is said to be a \textit{localization} if the forgetful functor $\mathcal{D}(\mathcal{B},\mathcal{N})\rightarrow \mathcal{\mathcal{D}(\mathcal{A},\mathcal{M})}$ is fully faithful.
\end{definition}

The class of localizations enjoys the following property:

\begin{proposition}[\cite{Analytic}, Exercise 12.17]
The class of localizations is stable under base change and composition and satisfies the two out of three property.
\end{proposition}

The following two examples of steady localizations will serve as a basis for constructing new ones:

\begin{proposition}
\label{important localizations}
The following two maps are steady localizations:
\[ (\mathbb{Z}[T],\mathbb{Z})_{\blacksquare}\rightarrow (\mathbb{Z}[T,T^{-1}],\mathbb{Z})_{\blacksquare},\]
\[ (\mathbb{Z}[T],\mathbb{Z})_{\blacksquare}\rightarrow \mathbb{Z}[T]_{\blacksquare}.\]
\end{proposition}

\begin{proof}
The first part follows from \cite[Example 13.15(1)]{Analytic} and the second from \cite[Example 13.15(2)]{Analytic}.
\end{proof}

One of the main ingredients of our proof is the following auxiliary result:

\begin{proposition}
\label{cool rational subsets}
Let $(A,A^+)$ be a complete Huber pair and $f,g$ elements of $A$ that generate an open ideal. Assume that the map
\[A\langle U\rangle\xrightarrow{gU-f}A\langle U\rangle\] is a closed embedding (this is true, for example, if $A$ is discrete or analytic and sheafy). Endow the condensed ring $A[g^{-1}]\overset{\mathrm{def}}{=} A\underset{\mathbb{Z}[T]}{\overset{L}{\otimes}}\mathbb{Z}[T, T^{-1}]$, where the $\mathbb{Z}[T]$-module structure on $A$ is given by $T\mapsto g$, with the $\mathbb{Z}[U]$-module structure given by $U\mapsto \frac{f}{g}$. Then we have the following isomorphism of analytic rings:
\[(A(\tfrac{f}{g}),A^+(\tfrac{f}{g}))_{\blacksquare}\cong ((A,A^+)_{\blacksquare}\underset{(\mathbb{Z}[T],\mathbb{Z})_{\blacksquare}}{\overset{L}{\otimes}}(\mathbb{Z}[T,T^{-1}],\mathbb{Z})_{\blacksquare})\underset{(\mathbb{Z}[U],\mathbb{Z})_{\blacksquare}}{\overset{L}{\otimes}}\mathbb{Z}[U]_{\blacksquare}.\]
Furthermore, the map \[(A,A^+)_{\blacksquare}\rightarrow(A(\tfrac{f}{g}),A^+(\tfrac{f}{g}))_{\blacksquare}\] is a steady localization.
\end{proposition}

\begin{proof}
Note that the second part of the statement follows from the first and Proposition~\ref{important localizations}. Arguing as in the proof of Lemma~\ref{trivial_bundle}, we see that it suffices to prove that we have

\[\underline{A(\tfrac{f}{g})}\cong (\underline{A}\underset{(\mathbb{Z}[T],\mathbb{Z})_{\blacksquare}}{\overset{L}{\otimes}}(\mathbb{Z}[T,T^{-1}],\mathbb{Z})_{\blacksquare})\underset{(\mathbb{Z}[U],\mathbb{Z})_{\blacksquare}}{\overset{L}{\otimes}}\mathbb{Z}[U]_{\blacksquare}.\]

By assumption and Lemma~\ref{surj}, the complex
\[\underline{A\langle U\rangle}\xrightarrow{gU-f}\underline{A\langle U\rangle}\] is quasi-isomorphic to $\underline {A(\frac{f}{g})}$. As $g$ is invertible in $A(\frac{f}{g})$ and the functor $-\underset{(\mathbb{Z}[T],\mathbb{Z})_{\blacksquare}}{\overset{L}{\otimes}}(\mathbb{Z}[T,T^{-1}],\mathbb{Z})_{\blacksquare}$ is exact, this complex is also quasi-isomorphic to
\[\underline{A\langle U\rangle}[g^{-1}]\xrightarrow{gU-f}\underline{A\langle U\rangle}[g^{-1}],\]
which is in turn isomorphic to
\[\underline{A\langle U\rangle}[g^{-1}]\xrightarrow{U-\frac{f}{g}}\underline{A\langle U\rangle}[g^{-1}].\]
Using the isomorphism $\underline{A\langle U\rangle}\cong \underline{A}[U]\underset{(\mathbb{Z}[U],\mathbb{Z})_{\blacksquare}}{\overset{L}{\otimes}}\mathbb{Z}[U]_{\blacksquare}$ and the fact that the functor $-\underset{(\mathbb{Z}[U],\mathbb{Z})_{\blacksquare}}{\overset{L}{\otimes}}\mathbb{Z}[U]_{\blacksquare}$ commutes with colimits, we see that the last complex is isomorphic to
\[\underline{A}[g^{-1}][U]\underset{(\mathbb{Z}[U],\mathbb{Z})_{\blacksquare}}{\overset{L}{\otimes}}\mathbb{Z}[U]_{\blacksquare}\xrightarrow{U-\frac{f}{g}}\underline{A}[g^{-1}][T]\underset{(\mathbb{Z}[U],\mathbb{Z})_{\blacksquare}}{\overset{L}{\otimes}}\mathbb{Z}[U]_{\blacksquare},\]
which is evidently isomorphic to $A[g^{-1}]\underset{(\mathbb{Z}[U],\mathbb{Z})_{\blacksquare}}{\overset{L}{\otimes}}\mathbb{Z}[U]_{\blacksquare}$, as desired.
\end{proof}

\begin{remark}
Let $A$ be a general (i.e., not necessarily sheafy) analytic complete Huber ring and let $f,g$ be elements in $A$ that generate an open ideal (equivalently, the unit ideal). Then the proof shows that
\[(A\underset{(\mathbb{Z}[T],\mathbb{Z})_{\blacksquare}}{\overset{L}{\otimes}}(\mathbb{Z}[T,T^{-1}],\mathbb{Z})_{\blacksquare})\underset{(\mathbb{Z}[U],\mathbb{Z})_{\blacksquare}}{\overset{L}{\otimes}}\mathbb{Z}[U]_{\blacksquare}\]
is quasi-isomorphic to the complex 
\[\underline{A\langle T \rangle}/^{\mathbb{L}}(gT-f):=[\underline{A\langle T\rangle}\xrightarrow{gT-f}\underline{A\langle T\rangle}].\]
\end{remark}

We now state and prove the main lemma, which lies at the heart of our proof. In what follows, we assume that $(A,A^+)$ is a sheafy analytic complete Huber pair unless otherwise specified, and denote the adic space $\mathrm{Spa}(A,A^+)$ by $X$ and the analytic ring $(A,A^+)_{\blacksquare}$ by $(\mathcal{A},\mathcal{M})$; additionally, for every affinoid open subspace $U\subset X$, we denote the ring $\mathcal{O}_{X}(U)$ by $A_U$, $\mathcal{O}_{X}^+(U)$ by $A_U^+$ and the analytic ring $(A_U,A_U^+)_{\blacksquare}$ by $(\mathcal{A}_{U},\mathcal{M}_U)$. We will also write $(\mathcal{A}(\frac{f}{g}),\mathcal{M}(\frac{f}{g}))$ as shorthand for $(\mathcal{A}_{X(\frac{f}{g})},\mathcal{M}_{X(\frac{f}{g})})$ and $L_U$ for the localization functor $-\overset{L}{\underset{(\mathcal{A},\mathcal{M})}{\otimes }}(\mathcal{A}_U,\mathcal{M}_U):\mathcal{C}_X\rightarrow C_U$.

\begin{proposition}
\leavevmode
\label{basechange}
\begin{enumerate}[label=(\roman*)]
    \item Let $U$ be an affinoid open subspace of $X$. Then the forgetful functor \[\mathcal{D}(\mathcal{A}_U,\mathcal{M}_U)\hookrightarrow \mathcal{D}(\mathcal{A},\mathcal{M})\] is a steady localization.
    \item Let $\phi:(A,A^+)\rightarrow (B,B^+)$ be a map of complete Huber pairs such that $(B,B^+)$ is analytic and sheafy, and $(A,A^+)$ is \textbf{either} analytic and sheafy \textbf{or} discrete. With the same notation as above, let $U$ be an nice rational open subspace of $X$ and denote its preimage in $Y:=\mathrm{Spa}(B,B^+)$ (a nice subspace as well) by $V$; also denote the analytic rings $(B,B^+)_{\blacksquare}$ and $(\mathcal{O}_{Y}(V),\mathcal{O}_{Y}^+(V))_{\blacksquare}$ by $(\mathcal{B},\mathcal{N})$ and $(\mathcal{B}_V,\mathcal{N}_V)$ respectively. Then have the following isomorphism of analytic rings
    \begin{equation}
        (\mathcal{B}_V,\mathcal{N}_V)\cong (\mathcal{B},\mathcal{N})\underset{(\mathcal{A},\mathcal{M})}{\overset{L}{\otimes}} (\mathcal{A}_U,\mathcal{M}_U).\tag{$\ast$}
    \end{equation}
    In particular, the first of the following diagrams is Cartesian and the second naturally commutes:
    
    \begin{center}
    \begin{tikzcd}
    \mathcal{D}(\mathcal{B}_V,\mathcal{N}_V) \arrow[r,] \arrow[d, ]
    & \mathcal{D}(\mathcal{B},\mathcal{N})\arrow[d, ] \\
    \mathcal{D}(\mathcal{A}_U,\mathcal{M}_U) \arrow[r,]
    & \mathcal{D}(\mathcal{A},\mathcal{M})
    \end{tikzcd}\hspace{40pt}%
\begin{tikzcd}[]
\mathcal{D}(\mathcal{B}) \arrow{r}{-\underset{\mathcal{B}}{\overset{L}{\otimes}}\mathcal{B}_V} \arrow[d,]
& \mathcal{D}(\mathcal{B}_V) \arrow[d, ] \\
\mathcal{D}(\mathcal{A}) \arrow{r}{-\underset{\mathcal{A}}{\overset{L}{\otimes}}\mathcal{A}_U}
& \mathcal{D}(\mathcal{A}_U)
\end{tikzcd}
\end{center}
    \item In the situation of $(ii)$, if $(A,A^+)$ is analytic and sheafy, then the isomorphism $(\ast)$ holds for an arbitrary affinoid subspace $U$.
    \item Let $U,V\subset X$ be affinoid open subspaces. Then we have \[\mathcal{D}(\mathcal{A}_{U\cap V},\mathcal{M}_{U\cap V})=\mathcal{D}(\mathcal{A}_U,\mathcal{M}_U)\cap \mathcal{D}(\mathcal{A}_V,\mathcal{M}_V)\] and \[L_{U\cap V}=L_{U}\circ L_V=L_V\circ L_U.\]
    \item Let $\{U_i\}_{i=1}^{n}$ be an open affinoid covering of $X$. Then the functor 
    \[\mathcal{D}(\mathcal{A},\mathcal{M}) \rightarrow \underset{i=1}{\overset{n}{\prod}}\mathcal{D}(\mathcal{A}_{U_i},\mathcal{M}_{U_i})\]
    is conservative. 
\end{enumerate}
\end{proposition}

\begin{proof}
We first prove (i)-(iii) in the case where $U$ is a nice rational open subspace. Arguing by induction, we easily reduce the proof of both parts to the case where $U$ is either $X(\frac{f}{1})$ or $X(\frac{1}{f})$. In this case, both $(i)$ and the isomorphism $(\ast)$ follow immediately from Proposition~\ref{cool rational subsets}. The rest of $(ii)$ follows from Propositions~\ref{pushout} and \ref{useful completion}.

We now deduce $(v)$ from the analogous statement for discrete Huber pairs. We need to show that for any $M\in \mathcal{D}(\mathcal{A},\mathcal{M})$ such that $L_{U_i}(M)=0$ for every $i=1,\ldots,n$, we have $M=0$. Applying Proposition~\ref{kedlnice}, we may assume that the covering is nice. Arguing by induction, we easily reduce the proof to the case where the covering $\{U\}_{i=1}^{n}$ is either simple Laurent or simple balanced. As the two cases are completely analogous, let us assume that $U$ is the simple Laurent covering $\mathcal{U}:=\{X(\frac{f}{1}),X(\frac{1}{f})\}$ defined by $f\in A$. Consider the canonical morphism $\phi: X\rightarrow X^{\mathrm{disc}}$ and the open covering ${\mathcal{V}:=\{X^{\mathrm{disc}}(\frac{f}{1}),X^{\mathrm{disc}}(\frac{1}{f})\}}$ of $X^{\mathrm{disc}}$. Observe that we have $X(\frac{f}{1})=\phi^{-1}(X^{\mathrm{disc}}(\frac{f}{1})),\ X(\frac{1}{f})=\phi^{-1}(X^{\mathrm{disc}}(\frac{1}{f}))$. Applying $(ii)$ and viewing $M$ as an element of $\mathcal{D}((A^{\mathrm{disc}},(A^{\mathrm{disc}})^{+})_{\blacksquare})$, we see that its restrictions to the elements of $\mathcal{V}$ (i.e., the images under the localization functors) are all zero. Therefore, it follows from the discrete case (see \cite[Lemma 10.3(ii)]{Condensed}) that $M$ is zero as an element of $\mathcal{D}((A^{\mathrm{disc}},(A^{\mathrm{disc}})^{+})_{\blacksquare})$. Therefore, $M$ is zero in $\mathcal{D}(\mathcal{A},\mathcal{M})$, since the forgetful functor $\mathcal{D}(\mathcal{A})\rightarrow \mathcal{D}(\mathcal{A}^{disc})$ has trivial kernel.

Applying Propositions~\ref{useful completion} and \ref{kedlnice}, we easily deduce $(i)$ for general $U$ from $(v)$.

To prove $(iii)$ for general affinoid $U$, we note that any morphism of complete Huber pairs with analytic domain is automatically adic. Therefore, the adic space $V$ is affinoid. We need to show that the natural map  \[\phi_S:(\mathcal{B},\mathcal{N})\underset{(\mathcal{A},\mathcal{M})}{\overset{L}{\otimes}} (\mathcal{A}_U,\mathcal{M}_U)[S] \rightarrow (\mathcal{B}_V,\mathcal{N}_V)[S]\]
is an isomorphism for every profinite set $S$. To do so, we pick an open cover of $U$ whose elements are nice open subspaces of $X$ (which is always possible by Proposition~\ref{kedlnice}). It follows from $(ii)$ that $\phi_S$ becomes an isomorphism after localizing to the elements of the cover, therefore we conclude by applying $(v)$.

Finally, $(iv)$ follows immediately from $(iii)$.
\end{proof}

We now can prove Theorem~\ref{maindescent}. As in the discrete case, it is a consequence of a more general categorical statement. Using Proposition~\ref{kedlnice} and Lemma~\ref{basechange}, we formally deduce Theorem~\ref{maindescent} from the following:

\begin{proposition}[\cite{Condensed}, Proposition 10.5]
Let $X$ be a spectral topological space equipped with a basis $B$ of quasi-compact open subsets stable under intersections. Let $\mathcal{C}$ be a stable $\infty$-category and $U\mapsto \mathcal{C}_U\subset C$ a covariant functor from $B$ to full subcategories of $\mathcal{C}$ such that every inclusion map $\mathcal{C}_U\hookrightarrow \mathcal{C}$ admits a left adjoint $L_U$. Furthermore, assume that
\begin{enumerate}[label=(\roman*)]
    \item for every pair $U,V\in B$, we have $\mathcal{C}_{U\cap V}=\mathcal{C}_U\cap \mathcal{C}_V$ and $L_{U\cap V}=L_{U}\circ L_V=L_V\circ L_U$;
    \item let $U$ be an element of $B$ and suppose that $\{U_i\}_{i=1}^{n}, U_i\in B$ form a covering of $U$; for every $M\in \mathcal{C}_X$ such that $L_{U_i}(M)=0$ for every $i=1,\ldots,n$, we have $M=0$.
\end{enumerate}
Then the covariant functor $U\mapsto \mathcal{C}_U$, with the left adjoints $L_U$ as the restriction functors, defines a sheaf of $\infty$-categories.
\end{proposition}
\section{Descent for pseudocoherent and perfect complexes and finite projective modules}\label{vb}
Let us begin by briefly recalling the relevant definitions. In what follows, we assume that $(A,A^+)$ is a sheafy analytic complete Huber pair unless otherwise specified, and denote the adic space $\mathrm{Spa}(A,A^+)$ by $X$ and the analytic ring $(A,A^+)_{\blacksquare}$ by $(\mathcal{A},\mathcal{M})$; additionally, for every affinoid open subspace $U\subset X$, we denote the rings $\mathcal{O}_{X}(U)$ and $\mathcal{O}_{X}^+(U)$ by $A_U$ and $A_U^+$, respectively, and the analytic ring $(A_U,A_U^+)_{\blacksquare}$ by $(\mathcal{A}_{U},\mathcal{M}_U)$. 

\begin{definition}
An object $M$ of $\mathcal{D}(A)$ is said to be a \textit{pseudocoherent complex} if it is quasi-isomorphic to a complex of the form
\begin{equation}
\label{form1}
\ldots\rightarrow P_{n+m}\rightarrow \ldots\rightarrow P_{n+1}\rightarrow P_{n}\rightarrow 0,~\tag{$\ddagger$}
\end{equation}
where $P_{i}$ are finite projective modules over $A$. 

An object $M$ is said to be a \textit{perfect complex} if it is quasi-isomorphic to a finite complex of the form~(\ref{form}). A perfect complex $M$ is said to have \textit{tor-amplitude in $[a,b]$}, where $a\leq b$ is a pair of integers, if it is quasi-isomorphic to a complex of the form~(\ref{form}) all of whose terms in cohomological degrees outside the interval $[a,b]$ are trivial.
\end{definition}

\begin{remark}
It can be easily shown that a complex $M$ is pseudocoherent if and only if it is quasi-isomorphic to a bounded above complex with finite free terms. Note, however, the analogous statement for perfect complexes is false.
\end{remark}

For the sake of completeness, we also recall the notion of vector bundle:

\begin{definition}
Let $(Y,\mathcal{O}_Y)$ be a ringed space. A sheaf $\mathcal{F}$ of $\mathcal{O}_Y$-modules on $Y$ is called a \textit{vector bundle}, if every point $y\in Y$ has an open neighborhood $V$ such that $\mathcal{F}|_{V}$ is isomorphic to $\mathcal{O}_Y^n|_{V}$ for some non-negative integer $n$. The category of vector bundles on $Y$ is denoted $\mathrm{Vect}_Y$.
\end{definition}
Whenever we talk about vector bundles in the context of adic spaces, by a vector bundle on an adic space $(Y,\mathcal{O}_Y,|\cdot(y)|_{y\in Y})$, we mean a vector bundle on the ringed space $(Y,\mathcal{O}_Y)$.

We denote the full $\infty$-subcategory of pseudocoherent complexes in $\mathcal{D}(A)$ by $\mathrm{PCoh}_A$, the $\infty$-\hspace*{0pt}subcategory of perfect complexes by $\mathrm{Perf}_A$ and the subcategory of finite projective modules by $\mathrm{FP}_A$. In addition, for every integer $n\in \mathbb{Z}$, we denote the full $\infty$-subcategory of pseudocoherent complexes in $\mathcal{D}(A)^{\leq n}$ by $\mathrm{PCoh}_A^{\leq n}$, and for every pair of integers $a\leq b$, we denote the $\infty$-subcategory of perfect complexes with tor-amplitude in $[a,b]$ by $\mathrm{Perf}_A^{[a,b]}$ (note that $\mathrm{FP}_A$ is clearly equal to $\mathrm{Perf}_A^{[0,0]}$). Our ultimate objective is to prove the following theorem:

\begin{theorem}
\label{descent of everything}
The functors sending every affinoid subspace $U\subset X$ to the $\infty$-categories $\mathrm{PCoh}_{A_U}$, $\mathrm{PCoh}_{A_U}^{\leq n}$, $\mathrm{Perf}_{A_U}$, $\mathrm{Perf}_{A_U}^{[a,b]}$, and $\mathrm{FP}_{A_U}$, respectively, define sheaves of $\infty$-categories on $X$.
\end{theorem}
 
As a formal consequence of Theorem~\ref{descent of everything}, we deduce the following result at the end of Section~\ref{vb2}:
 
\begin{theorem}
\label{vector bundles descent}
There is a natural equivalence of categories $\Upsilon:\mathrm{FP}_A\xrightarrow[]{} \mathrm{Vect}_X$, whose quasi-inverse is given by the functor of global sections.
\end{theorem}

The section is organized as follows. In the first subsection, we explain the construction of the functor $\Upsilon:\mathrm{Mod}_A\rightarrow \mathcal{O}_X$-$\mathrm{Mod}$ and put the discussion into the context of condensed mathematics by introducing the \textit{condensification functor}; in particular, we introduce the notion of \textit{(relatively) discrete condensed modules}. Then we introduce the notion of \textit{pseudocoherence} and discuss its different manifestations and their interrelationships in the second subsection. In the third subsection, we develop the technical tools necessary for our proof of Theorem~\ref{descent of everything}. Specifically, we study nuclear, compact, and dualizable objects in detail and prove that they satisfy descent with respect to the analytic topology on $X$. Finally, we prove Theorem~\ref{descent of everything} and discuss its corollaries in the concluding subsection.
\subsection{The functor \texorpdfstring{$\Upsilon$}{} and condensification}\label{vb1}

Let $U$ be an affinoid subspace of $X$ and $M$ a module over $A$. In what follows, we often refer to $M\underset{A}{\otimes}A_U$ as \textit{the restriction of $M$ to $U$}.
\begin{definition}\label{sheaves}
Let $M$ be a module over $A$. The sheaf of modules associated to the presheaf sending every affinoid open subspace $U$ to the restriction of $M$ to $U$ is called \textit{the sheaf of $\mathcal{O}_X$-modules associated to $M$} and denoted $\Tilde{M}$. The functor $\mathrm{Mod}_A\rightarrow \mathcal{O}_X$-$\mathrm{Mod}$ given by $M\mapsto \Tilde{M}$ is denoted by $\Upsilon$.
\end{definition}

It is not a priory clear if the space of global section of $\Tilde{M}$ is equal to $M$; in fact, it may be false in general. However, we will later prove that this is actually true for finite projective modules, or more generally for stably pseudocoherent modules. The following lemma guarantees that the functor $\Upsilon:\mathrm{FP}_A\xrightarrow[]{} \mathrm{Vect}_X$ is well-defined:

\begin{lemma}
Let $P$ be a finite projective module over $A$. Then there exists a rational covering of $X$ such that the restriction of $P$ to any of its elements is finite and free.
\end{lemma}

\begin{proof}
Consider the module $P$ as a vector bundle on the scheme $\spec A$. Pick elements $f_1,\dots,f_n$ generating the unit ideal such that $P_{f_i}$ is free for each $i=1,\dots,n$, where $P_{f_i}$ denotes the localization of $P$ at the multiplicative set generated by $f$ in the usual sense of commutative algebra. Let $s_1,\dots, s_n$ be elements of $A$ such that $s_1f_1+\dots +s_nf_n=1$. Since $P_{f_i}$ is free, so is $P_{s_if_i}$. Observe that the rational open subspaces $X(\frac{1}{s_1f_1}),\dots, X(\frac{1}{s_nf_n})$ form a covering of $X$. As the formation of the restriction of $P$ to $X(\frac{1}{s_if_i})$ factors through the formation of $P_{s_if_i}$, we see that the restriction of $P$ to each $X(\frac{1}{s_if_i})$ is finite and free.
\end{proof}
 
As our proof of the main theorem uses techniques of condensed mathematics, we need to transfer the problem to that context first. For this, we prove the following auxiliary result: 
 
\begin{lemma}
\label{condensification}
Let $R$ be a discrete ring and endow every module over $R$ with the discrete topology. Then the functor $C:\mathrm{Mod}_R\rightarrow \mathrm{Mod}_{\underline{R}}^{\mathrm{cond}}$ given by $M\mapsto \underline{M}$ extends to a fully faithful exact functor $\mathrm{dCond}_R:D(R)\hookrightarrow D(\underline{R})$, which we call \textit{the discrete condensification functor}. Moreover, it can be explicitly described as follows: for every object $P\in D(R)$, $\mathrm{dCond}_R(P)$ is given by $\underline{N^{\bullet}}$, where $N^{\bullet}$ is any complex representing $P$. In addition, $\mathrm{dCond}_R$ commutes with filtered colimits and its image lies in $\mathcal{D}(R_{\blacksquare})$.
\end{lemma}

\begin{proof}
Most of the proposition is a formal consequence of the fact that the functor $C$ is exact (it is left exact trivially and right exact as the ground ring and modules are discrete) and fully faithful (by \cite[Proposition 1.7]{Condensed}) and preserves projective objects (as it preserves direct sums and $\underline{R}$ is projective). To show that $\mathrm{dCond}_R$ commutes with filtered colimits, it suffices to prove that $C$ does, which follows from the fact all modules are endowed with the discrete topology. The very last part of the theorem follows immediately from the above and the fact that the category $\mathcal{D}(R_{\blacksquare})$ is closed under colimits.
\end{proof}

This is, however, not sufficient for our purposes. In the context of Huber rings, we also need to take into account their topological structures. The statements we need, though, can be formulated for arbitrary condensed rings. 

\begin{definition}
Let $R$ be a condensed ring and denote its underlying ring endowed with the discrete topology by $R(\ast)$. \textit{The condensification functor} is the functor $\mathrm{Cond}_{R}:D(R(\ast))\rightarrow D(R)$ given by $M\mapsto \mathrm{dCond}_{R(\ast)}(M)\underset{\underline{R(\ast)}}{\overset{L}{\otimes}}R$. An object of $D(R)$ is said to be \textit{(relatively) discrete} if it lies in the essential image of $\mathrm{Cond}_R$.
\end{definition}

In the case where the condensed ring $R$ is given by $\underline{S}$ for some topological ring $S$, we will generally abuse notation by denoting the condensification functor $\mathrm{Cond}_{S}$.

\begin{theorem}
\label{embed}
Let $R$ be a condensed ring. Then the condensification functor is fully faithful and exact and commutes with filtered colimits. In addition, $\mathrm{Cond}_R(R(\ast))$ is quasi-isomorphic to $R$, and for any $M\in D(R(\ast))$ and $i\in \mathbb{Z}$, $H^i(\mathrm{Cond}_R (M))= 0$ implies that $H^i(M)=0$. 

Furthermore, if $R$ is given by $\underline{A}$ for an analytic complete Huber ring $A$, the image of $\mathrm{Cond}_{A}$ lies in $\mathcal{D}((A,A^+)_{\blacksquare})$.\footnote{For any Huber pair $(A,A^+)$.} In addition, if $M \in D(A)$ is quasi-isomorphic to a complex of finite free modules, then for each $i\in \mathbb{Z}$, $H^i(\mathrm{Cond}_A (M))= 0$ if and only if $H^i(M)=0$. 
\end{theorem}

\begin{proof}
As both $\mathrm{dCond}_{R(\ast)}$ and $-\underset{\underline{R(\ast)}}{\overset{L}{\otimes}}R$ are exact and commute with filtered colimits, the same is true for $\mathrm{Cond}_R$. Since $\underline{R(\ast)}$ is projective in $\mathrm{Mod}^{\mathrm{cond}}_{\underline{R(\ast)}}$, we have $\mathrm{Cond}_R(R(\ast))\cong R$. Consider a pair of chain complexes of $R(\ast)$-modules $M^{\bullet},N^{\bullet}$ with free terms. As we have $\mathrm{Hom}_{R(\ast)}(R(\ast)^{\oplus I},R(\ast)^{\oplus J})\cong \mathrm{Hom}_{R}(R^{\oplus I},R^{\oplus J})$ for every pair of sets $I,J$, we easily deduce that $\mathrm{Hom}_{D(R(\ast))}(M^{\bullet},N^{\bullet})$ is isomorphic to $\mathrm{Hom}_{D(R)}(\underline{M}^{\bullet},\underline{N}^{\bullet})\cong \mathrm{Hom}_{D(R)}(\mathrm{Cond}_R{(M}^{\bullet}),\mathrm{Cond}_R{(N}^{\bullet}))$. Therefore, the condensification functor is fully faithful. Now let $M$ be an object of $D(R(\ast))$ and $M^{\bullet}$ a chain complex with free terms that represents it. Note that for every $i\in \mathbb{Z}$, we evidently have $H^i(M^{\bullet})\cong H^i(\underline{M^{\bullet}})(*)$. Therefore, if $H^i(\mathrm{Cond}_R(M))$ is trivial, so is $H^i(M)$. 

We now assume that $R$ is given by $\underline{A}$ for an analytic complete Huber ring $A$. The second to last claim follows immediately from the above and the fact that the category $\mathcal{D}((A,A^+)_{\blacksquare})$ is closed under colimits. Finally, suppose that all the terms of $M^{\bullet}$ are finite. It remains to prove that if $H^i(M)$ is trivial, then so is $H^i(\mathrm{Cond}_A(M))$. As the kernel of $d^i:M^i\rightarrow M^{i+1}$ is clearly closed, we see that it is complete and first countable. Therefore, the desired claim follows from Theorem~\ref{open mapping} and Lemma~\ref{surj}.
\end{proof}

\begin{remark}
Whenever there is no ambiguity what base ring $A$ is meant, we suppress it from the notation and simply write $\mathrm{Cond}$ for the condensification functor.
\end{remark}

We record the following property of the condensification functor for later use:

\begin{lemma}
\label{discrete_localization}
Let $R\rightarrow S$ be a map of condensed rings. Then the following diagram naturally commutes:

\begin{center}
\begin{tikzcd}[column sep=4em]
D(R(\ast))\arrow[r,hook,"\mathrm{Cond}_R"]\arrow{d}{-\underset{R(\ast)}{\overset{L}{\otimes}} S(\ast)} & D(R)\arrow{d}{-\underset{R}{\overset{L}{\otimes}} S}\\
D(S(\ast))\arrow[r,hook,"\mathrm{Cond}_{S}"] & D(S)
\end{tikzcd}
\end{center}

In particular, the following diagram naturally commutes for any affinoid subspace $U$ of $X$:

\begin{center}
\begin{tikzcd}[column sep=4em]
D(A)\arrow[r,hook,"\mathrm{Cond}_A"]\arrow{d}{-\underset{A}{\overset{L}{\otimes}} A_U} & D(\mathcal{A},\mathcal{M})\arrow{d}{-\underset{(\mathcal{A},\mathcal{M})}{\overset{L}{\otimes}} (\mathcal{A}_U,\mathcal{M}_U)}\\
D(A_U)\arrow[r,hook,"\mathrm{Cond}_{A_U}"] & D(\mathcal{A}_U,\mathcal{M}_U)
\end{tikzcd}
\end{center}
\end{lemma}

\begin{proof}
The claim readily follows from the fact that every complex of $R(\ast)$-modules is quasi-isomorphic to a complex with free terms.
\end{proof}

This lemma will allow us to prove Theorem~\ref{descent of everything} by means of condensed mathematics. However, we first record the following property of the condensification functor:

\begin{lemma}
\label{retrdisc}
Let $X$ be a discrete object of $D(\mathcal{A},\mathcal{M})$. Then any retract of $X$ is again discrete.  
\end{lemma}

\begin{proof}
As retracts of $X$ are in one to one correspondence with idempotent endomorphisms of $X$, the claim follows from the fact that $\mathrm{Cond}$ is fully-faithful and exact. 
\end{proof}

We may now describe the essential image of $\mathrm{FP}_A$ under $\mathrm{Cond}$:

\begin{proposition}
\label{different projectives}
Let $P$ be a condensed module over $(\mathcal{A},\mathcal{M})$. It is a retract of $\underline{A}^n$ for some $n\geq 0$ if and only of it is discrete and comes from a finite projective module over $A$. 
\end{proposition}

\begin{proof}
The implication ``$\implies$'' was proved in Lemma~\ref{retrdisc}; the reverse implication is trivial. 
\end{proof}

\subsection{Pseudocoherence}\label{vb2}

The goal of this section is to discuss the notion of pseudocoherence in the context of analytic adic spaces. We begin with the following general categorical definition:

\begin{definition}
\label{cat ps}
Let $\mathcal{A}$ be an abelian category that admits all small colimits and has a family of compact projective generators. An object $M\in\mathcal{A}$ is called \textit{pseudocoherent} if for each $i\geq 0$, the functor $\mathrm{Ext}^i_{\mathcal{A}}(M,-):\mathcal{A}\rightarrow\mathrm{Ab}$ commutes with filtered colimits. An object $X\in \mathcal{D}^{-}(\mathcal{A})$ is called \textit{pseudocoherent} if for each $i\in\mathbb{Z}$, the functor $\mathrm{Hom}_{\mathcal{D}(\mathcal{A})}(X,-):\mathcal{D}^{\geq i}(\mathcal{A})\rightarrow\mathrm{Ab}$ commutes with filtered colimits in $\mathcal{D}^{\geq i}(\mathcal{A})$.
\end{definition}

It can be reformulated in more explicit terms as follows:

\begin{proposition}\label{pseudo}
\leavevmode
\begin{enumerate}[label=(\roman*)]
    \item An object $M\in\mathcal{A}$ is pseudocoherent if and only if there exists a projective resolution of the form
    \[\ldots\rightarrow\underset{i=1}{\overset{k_n}{\bigoplus}}P_{n,i}\rightarrow \ldots\rightarrow \underset{i=1}{\overset{k_1}{\bigoplus}}P_{1,i}\rightarrow\underset{i=1}{\overset{k_0}{\bigoplus}}P_{0,i}\rightarrow M\rightarrow 0,\]
    where all $P_{ji}$ are from a fixed family of compact projective generators and all $k_i$ are non-negative integers. In addition, every partial resolution of this form may be prolonged. 
    \item An object $X\in \mathcal{D}^-(\mathcal{A})$ is pseudocoherent if and only if it is quasi-isomorphic to a bounded above complex with terms of the form $\underset{i=1}{\overset{k}{\bigoplus}}P_{i}$, where all $P_{i}$'s are from a fixed family of compact projective generators and $k$ is a non-negative integer. In addition, if $H^i(M)=0$ for all $i>n$, then this complex may be chosen so that its terms of degree $i>n$ are trivial.
\end{enumerate}
In particular, an object of $\mathcal{A}$ is pseudocoherent if and only if it is pseudocoherent as an object of $\mathcal{D}(\mathcal{A})$.
\end{proposition}

\begin{proof}
For part $(i)$, see \cite[Proposition 4.12]{Condensed}. Part $(ii)$ is analogous to \cite[\href{https://stacks.math.columbia.edu/tag/094B}{Tag 094B}]{stacks-project}.
\end{proof}

Proposition~\ref{pseudo} motivates the following definitions: 
\begin{definition}
\label{ps mod}
A module $M$ over $A$ is called \textit{pseudocoherent} if it admits a projective resolution by finite free $A$-modules. We say that $M$ is \textit{stably pseudocoherent} if, in addition, for every affinoid open subspace $U\subset X$, we have $M\overset{L}{\underset{A}{\otimes}} A_U\cong M\overset{}{\underset{A}{\otimes}} A_U$. The full subcategory of stably pseudocoherent modules is denoted $\mathrm{PCoh}^0_A$. 
\end{definition}

\begin{definition}
\label{top ps mod}
A module $M$ over $A$ is called \textit{topological pseudocoherent} if it is pseudocoherent and complete with respect to the natural topology. We say that $M$ is \textit{topological stably pseudocoherent} if, in addition, it is stably pseudocoherent and for every affinoid open subspace $U\subset X$, the restriction $M\overset{}{\underset{A}{\otimes}} A_U$ is complete with respect to the natural topology. The full subcategory of topological stably pseudocoherent modules of $A$ is denoted $\mathrm{tPCoh}^0_A$.
\end{definition}

These notions have the following obvious counterparts in the framework of sheaves on $X$, which can be regarded as analogues of coherent sheaves in the context of algebraic geometry:

\begin{definition}
\label{ps sh}
An $\mathcal{O}_X$-module $\mathcal{F}$ is called \textit{pseudocoherent} if every point $x\in X$ has a affinoid open neighborhood $U\subset X$ such that $\mathcal{F}|_{U}\cong \Tilde{M}$ for some stably pseudocoherent module $M$ over $A_U$. The full subcategory of pseudocoherent sheaves on $X$ is denoted $\mathrm{PCoh}^0_X$.
\end{definition}

\begin{definition}
\label{top ps sh}
An $\mathcal{O}_X$-module $\mathcal{F}$ is called \textit{topological pseudocoherent} if every point $x\in X$ has a affinoid open neighborhood $U\subset X$ such that $\mathcal{F}|_{U}\cong \Tilde{M}$ for some topological stably pseudocoherent module $M$ over $A_U$. The full subcategory of topological pseudocoherent sheaves on $X$ is denoted $\mathrm{tPCoh}^0_X$.
\end{definition}

Finally, let us specialize the notion of pseudocoherence to the setting of modules over $(\mathcal{A},\mathcal{M})$:

\begin{definition}
\label{condensed ps}
A condensed module $M\in\mathrm{Mod}_{(\mathcal{A},\mathcal{M})}$ is called \textit{pseudocoherent} if it admits a projective resolution with terms of the form $\mathcal{M}[S]$ for various profinite sets $S$. A condensed module $M$ is called \textit{stably pseudocoherent}, if for every affinoid open subspace $U\subset X$, we have \[M\overset{L}{\underset{(\mathcal{A},\mathcal{M})}{\otimes}} (\mathcal{A}_U,\mathcal{M}_U)\cong M\overset{}{\underset{(\mathcal{A},\mathcal{M})}{\otimes}} (\mathcal{A}_U,\mathcal{M}_U).\] The full subcategory of stably pseudocoherent condensed modules over $(\mathcal{A},\mathcal{M})$ is denoted $\mathrm{PCoh}^0_{(\mathcal{A},\mathcal{M})}$. 

An object of $\mathcal{D}(\mathcal{A},\mathcal{M})$ is called \textit{pseudocoherent} if it is quasi-isomorphic to a bounded above complex with terms of the form $\mathcal{M}[S]$ for various profinite sets $S$. The full $\infty$-subcategory of pseudocoherent objects in $\mathcal{D}(\mathcal{A},\mathcal{M})$ is denoted $\mathrm{PCoh}_{(\mathcal{A},\mathcal{M})}$.
\end{definition}

\begin{remark}
The notions of a pseudocoherent module and pseudocoherent sheaf are strongly inspired by the ideas of Kedlaya and Liu, see \cite[Definition 1.4.5,Definition 1.4.17.]{Kdl}. Nevertheless, our definitions and those of Kedlaya and Liu are not quite the same. Namely, the modules that they call \textit{pseudocoherent} we also call pseudocoherent, but those that they term \textit{stably pseudocoherent} essentially correspond to \textit{topological stably pseudocoherent modules} in the terminology of this paper. However, there is a minor but important difference: we also require that any affinoid localization of a topological stably pseudocoherent module is still concentrated in degree $0$.
\end{remark}

Observe that pseudocoherent condensed modules admit the following equivalent description: 

\begin{proposition}
A condensed module $M$ is pseudocoherent if and only if, for each $i\in\mathbb{Z}$, the functor $\underline{\mathrm{Ext}}^i_{(\mathcal{A},\mathcal{M})}(M,-):\mathrm{Mod}_{(\mathcal{A},\mathcal{M})}\rightarrow \mathrm{Mod}_{(\mathcal{A},\mathcal{M})}$ commutes with all filtered colimits.
\end{proposition}

\begin{proof}
The implication ``$\implies$'' follows directly from Definition~\ref{condensed ps}. For the converse, we first observe that for every pair of objects $X,Y\in \mathrm{Mod}_{(\mathcal{A},\mathcal{M})}$, we have \[\mathrm{Ext}_{(\mathcal{A},\mathcal{M})}^{i}(X,Y)\cong\underline{\mathrm{Ext}}^i_{(\mathcal{A},\mathcal{M})}(X,Y)(*),\] and then use Definition~\ref{cat ps} and Proposition~\ref{pseudo} to conclude that $M$ is pseudocoherent.
\end{proof}

The following lemma ensures that the topological pseudocoherent sheaves are indeed a generalization of vector bundles:

\begin{lemma}
A finite projective module $P$ over $A$ is topological stably pseudocoherent.
\end{lemma}

\begin{proof}
It clearly suffices to show that $P$ is complete with respect to its natural topology, since for every affinoid open subspace $U\subset X$, the localization functor $-\overset{L}{\underset{A}{\otimes}}A_U$ sends finite projective modules to finite projective ones. Pick a surjection $A^n\twoheadrightarrow P$ and consider the following short exact sequence:
\[0\xrightarrow{}\ker\phi\xrightarrow{} A^m\xrightarrow{\phi} P\xrightarrow{} 0.\]
As the sequence is split, $\ker\phi$ is also finitely generated and we endow it with the natural topology. Using that $P$ is projective again, we obtain an algebraic section $s: A^n\xrightarrow{} \ker\phi$. Since $\ker\phi$ is a topological $A$-module, this map is also continuous. As every retract of a Hausdorff space is closed in that space, we see that $\ker\phi$ is closed in $A^n$; by general properties of linearly topologized modules, it implies that $P$ is complete and first countable with respect to the natural topology and $\phi$ is open (see \cite[\href{https://stacks.math.columbia.edu/tag/0AMT}{Tag 0AMT}]{stacks-project}).
\end{proof}

We now relate the notion of pseudocoherence over $A$ to that over $(\mathcal{A},\mathcal{M})$ by means of the condensification functor.

\begin{proposition}
\label{pc to cond pc}
The functor $\mathrm{Cond}$ sends pseudocoherent and stably pseudocoherent modules over $A$ to pseudocoherent and stably pseudocoherent modules over $(\mathcal{A},\mathcal{M})$, respectively. Furthermore, it sends pseudocoherent objects in $\mathcal{D}(A)$ to pseudocoherent objects in $\mathcal{D}(\mathcal{A},\mathcal{M})$. Conversely, every discrete pseudocoherent (resp. stably pseudocoherent) module over $(\mathcal{A},\mathcal{M})$ comes from a pseudocoherent (resp. stably pseudocoherent) module over $A$. Moreover, every discrete pseudocoherent object in $\mathcal{D}(\mathcal{A},\mathcal{M})$ comes from a pseudocoherent object in $\mathcal{D}(A)$.
\end{proposition}

\begin{proof}
The first two claims readily follow from Proposition~\ref{embed} and Proposition~\ref{discrete_localization}. For the last two, it clearly suffices to prove that a pseudocoherent object $X\in \mathcal{D}(\mathcal{A},\mathcal{M})$ that is discrete as an object of $\mathcal{D}(\mathcal{A},\mathcal{M})$ actually comes from a pseudocoherent object of $\mathcal{D}(A)$. It is easily verified that there exists a bounded above complex $N^{\bullet}\in \mathrm{Ch}^-(A)$ with free terms such that $X$ is quasi-isomorphic to $\underline{N^{\bullet}}$. As the condensification functor is fully faithful and commutes with colimits, it immediately follows from the description of pseudocoherent objects provided by Definition~\ref{cat ps} that $N^{\bullet}$ is pseudocoherent as an object of $\mathcal{D}(A)$.
\end{proof}

\begin{lemma}
\label{condensing_fpr}
Let $M$ be a pseudocoherent module over $A$ that is complete with respect to the natural topology. Pick a finite presentation of $M$, say 
\[A^n\xrightarrow{} A^m\xrightarrow{\phi} M\xrightarrow{} 0.\]
Then $\underline{M}$ is isomorphic to the cokernel of the induced map $\underline{\phi}:\underline{A}^n\xrightarrow{} \underline{A}^m$. Therefore, $\underline{M}$ is $(\mathcal{A},\mathcal{M})$-complete. Furthermore, $\underline{M}$ and $\mathrm{Cond}(M)$ are isomorphic as objects of $\mathcal{D}(\mathcal{A},\mathcal{M})$. Finally, if $M$ is finite projective, then $\underline{M}$ is compact and projective as well.
\end{lemma}
\begin{proof}
By definition, $\phi$ is a continuous surjective map of first-countable and complete modules over $A$. Therefore, it is open by virtue of the Theorem~\ref{open mapping}. Applying Lemma~\ref{surj}, we see that $\underline{M}$ is isomorphic to the cokernel of $\underline{\phi}$. As $H^i(\mathrm{Cond}(M))=0$ if and only if $H^i(M)=0$ (Theorem~\ref{embed}), we also have $\underline{M}\cong \mathrm{Cond}(M)$. To prove the last statement of the proposition, observe that $\underline{M}$ is a retract of a compact and projective object (by Lemma~\ref{different projectives}), so it is compact and projective itself.
\end{proof}

For (topological) stably pseudocoherent modules, we have the following analogue of Theorem~\ref{descent of everything}:

\begin{theorem}
\label{pseudo stack}
The associations sending every affinoid subspace $U\subset X$ to the categories $\mathrm{PCoh}^0_{A}$ and $\mathrm{tPCoh}^0_{A}$, respectively, define sheaves of categories on $X$.
\end{theorem}

We prove Theorem~\ref{pseudo stack} in the last subsection. As for vector bundles, we also have the following result:

\begin{theorem}
\label{pseudocoherent descent}
\leavevmode
\begin{enumerate}[label=(\roman*)]
    \item The functor $\Upsilon$ restricts to an equivalence $\Upsilon:\mathrm{PCoh}^0_A\xrightarrow[]{} \mathrm{PCoh}^0_X$, whose quasi-inverse is given by the functor of global sections.
    \item The functor $\Upsilon$ restricts to an equivalence $\Upsilon:\mathrm{tPCoh}^0_A\xrightarrow[]{} \mathrm{tPCoh}^0_X$, whose quasi-inverse is given by the functor of global sections.
\end{enumerate}
\end{theorem}

Let us explain how Theorems~\ref{descent of everything} and \ref{pseudo stack} imply Theorems~\ref{vector bundles descent} and \ref{pseudocoherent descent}. We need to show that for any stably pseudocoherent module $M$ over $A$, the presheaf sending an affinoid open subspace $U$ to $M\underset{A}{\otimes}A_U$ is in fact a sheaf. Note that this is obvious if $M$ is finite projective, as in this case, $M$ is a retract of a finite free module. For general stably pseudocoherent $M$, we identify the elements of $M$ with maps $A\rightarrow M$. Invoking Theorem~\ref{pseudo stack}, we see that any local datum of such maps uniquely defines a global map, which implies the sheaf axioms, as desired.
\subsection{Nuclear, compact, dualizable, and pseudocoherent objects in \texorpdfstring{$\mathcal{D}(\mathcal{A},\mathcal{M})$}{}}\label{vb3}
In this subsection, we develop the main technical tools that we need for our proofs of Theorems~\ref{descent of everything} and \ref{pseudo stack}. Following \cite{Analytic}, we introduce the notion of a \textit{nuclear object} and that of a \textit{dualizable object}, establish their basic properties and relate them to the notions of compact and pseudocoherent objects. After that, we prove that these four classes of objects in $D(\mathcal{A},\mathcal{M})$ satisfy descent (see Theorems~\ref{compact_descent}, \ref{nuclear_descent}, \ref{dual_descent} and \ref{pseudocoherent objects descent}). First, we need to introduce the notion of duality in the present context.

\begin{definition}
\label{definition duality}
Let $P$ be an object of $\mathcal{D}(\mathcal{A})$. Then the object $\mathrm{R}\underline{\mathrm{Hom}}_{\mathcal{A}}(P,\mathcal{A})$ is called \textit{the dual of $P$} and denoted $P^{\vee}$.
\end{definition}

We record the following easy result for later use:

\begin{lemma}
\label{diffdescofdual}
\leavevmode
\begin{enumerate}[label=(\roman*)]
    \item For every $X\in \mathcal{D}(\mathcal{A})$ and $Y\in \mathcal{D}(\mathcal{A},\mathcal{M})$, we have 
    \[\mathrm{R}\underline{\mathrm{Hom}}_{\mathcal{A}}(X,Y) \cong \mathrm{R}\underline{\mathrm{Hom}}_{\mathcal{A}}(X\underset{\mathcal{A}}{\overset{L}{\otimes}}(\mathcal{A},\mathcal{M}),Y).\]
     \item For every pair of objects $X,Y\in \mathcal{D}(\mathcal{A},\mathcal{M})$, we have 
    \[\mathrm{R}\underline{\mathrm{Hom}}_{(\mathcal{A},\mathcal{M})}(X,Y)\cong \mathrm{R}\underline{\mathrm{Hom}}_{\mathcal{A}}(X,Y).\]
    \item For every profinite set $S$, we have $\mathcal{M}[S]^{\vee}\cong \underline{C(S,A)}$.
\end{enumerate}
In particular, for every $P\in \mathcal{D}(\mathcal{A})$, we have $P^{\vee}\cong (P\underset{\mathcal{A}}{\overset{L}{\otimes}}(\mathcal{A},\mathcal{M}))^{\vee}\in \mathcal{D}(\mathcal{A},\mathcal{M})$.
\end{lemma}

\begin{proof}
Parts $(i)$ and $(ii)$ readily follow from the fact that the forgetful functor $\mathcal{D}(\mathcal{A},\mathcal{M})\xrightarrow[]{} \mathcal{D}(\mathcal{A})$ is fully faithful and the functor  $-\underset{\mathcal{A}}{\overset{L}{\otimes}}(\mathcal{A},\mathcal{M})$ is left adjoint to it. Part $(iii)$ is an easy consequence of $(i)$ and $(ii)$. Indeed, as $\mathcal{M}[S]$ is projective for every profinite set $S$, the object $\mathcal{M}[S]^{\vee}$ is concentrated in degree $0$. We now compute
\[\mathcal{M}[S]^{\vee}\cong \underline{\mathrm{Hom}}_{\mathcal{A}}(\mathcal{A}[S],\mathcal{A})\cong \underline{C(S,A)}.\]
\end{proof}

The notion of a dualizable object is in fact very general and applicable to every symmetric monoidal category.  

\begin{definition}
Let $(\mathcal{C},\otimes,1)$ be a symmetric monoidal category. An object $X$ is called \textit{dualizable} if there exists an object $X^{\vee}$ together with maps $\mathrm{ev}_X:X^{\vee}\otimes X\rightarrow 1$, $i_X: 1\rightarrow X\otimes X^{\vee}$ such that the following diagrams commute:
\begin{center}
\begin{tikzcd}[column sep=4em]
X^{\vee} \otimes (X\otimes X^{\vee})\arrow{d}{\alpha^{-1}_{X^{\vee},X,X^{\vee}}} &  \arrow{l}[above]{\mathrm{id}_{X^{\vee}}\otimes i_X} X^{\vee}\otimes 1 \arrow{d}{l^{-1}_{X^{\vee}}\circ r_{X^{\vee}}}\\
(X^{\vee} \otimes X)\otimes X^{\vee} \arrow{r}{\mathrm{ev}_X\otimes\mathrm{id}_{X^{\vee}}} & 1 \otimes X^{\vee} 
\end{tikzcd}
\hspace{40pt}%
\begin{tikzcd}[column sep=4em]
(X \otimes X^{\vee})\otimes X \arrow{d}{\alpha_{X,X^{\vee},X}} &\arrow{l}[above]{i_X\otimes \mathrm{id}_{X}} X\otimes 1 \arrow{d}{r^{-1}_{X}\circ l_{X}}\\
X \otimes (X^{\vee}\otimes X) \arrow{r}{\mathrm{id}_{X} \otimes \mathrm{ev}_X} & 1 \otimes X 
\end{tikzcd}
\end{center}
where $\alpha$ denotes the associator of $\mathcal{C}$ and $l,r$ denote the left and right unitors respectively.
\end{definition}

\begin{remark}
Whenever we write $M^{\vee}$ later in the text, we follow the notation of Definition~\ref{definition duality} and \textit{do not} imply that $M$ is dualizable. 
\end{remark}

Dualizable objects enjoy the following useful property:

\begin{lemma}
\label{dualizable}
Let $(\mathcal{C},\otimes,1)$ be a closed symmetric monoidal category, $X$ a dualizable object and $X^{\vee}$ its dual. Then for any object $Z$, there is a canonical isomorphism $X^{\vee}\otimes Z\cong \underline{\mathrm{Hom}}(X,Z)$.
\end{lemma}

\begin{proof}
We prove the statement by constructing mappings \[\phi: X^{\vee}\otimes Z\rightarrow \underline{\mathrm{Hom}}(X,Z)\ \mathrm{and}\ {\psi: \underline{\mathrm{Hom}}(X,Z) \rightarrow X^{\vee}\otimes Z}\] that are inverse to each other. By tensor-hom adjunction, we may specify a map $X\otimes X^{\vee}\otimes Z\rightarrow Z$ in order to define $\phi$. There is a canonical choice for such a map: $\mathrm{ev}_X\otimes \mathrm{id}_Z$. The mapping $\psi$ is defined as the composition of the canonical mapping $\underline{\mathrm{Hom}}(X,Z)\rightarrow \underline{\mathrm{Hom}}(X^{\vee}\otimes X,X^{\vee}\otimes Z)$ with $i_X^{*}$. An easy diagram chase shows that these maps are indeed inverse to each other and their construction is functorial in $Z$.
\end{proof}

We now define nuclear objects:

\begin{definition}[\cite{Analytic}, Definition 13.10]
\label{definition nuclear}
An object $M\in \mathcal{D}(\mathcal{A},\mathcal{M})$ is called \textit{nuclear} if for all extremally disconnected sets $S$, the natural map \[(\mathcal{M}[S]^{\vee}\overset{L}{\underset{(\mathcal{A},\mathcal{M})}{\otimes}} M)(*)\rightarrow M(S)\] is an isomorphism in $\mathcal{D}(\mathrm{Ab})$.
\end{definition}

To state some of their basic properties, we need the following definitions:

\begin{definition}[\cite{Analytic}, Definition 13.11]
A map $f:P\rightarrow Q$ between compact objects of $\mathcal{D}(\mathcal{A},\mathcal{M})$ is called \textit{trace class} if there is some map $g: \mathcal{A} \rightarrow P^{\vee}\overset{L}{\underset{\mathcal{A}}{\otimes}} Q$ such that $f$ is the composite $P\xrightarrow{1\otimes g} P\overset{L}{\underset{(\mathcal{A},\mathcal{M})}{\otimes}}P^{\vee}\overset{L}{\underset{(\mathcal{A},\mathcal{M})}{\otimes}} Q\xrightarrow{} Q$, where the second map contracts the first two factors.
\end{definition}
\begin{remark}
This notion is stable under localization. Indeed, let $U$ be a affinoid open subspace of $X$. As localization preserves compact objects (see Theorem~\ref{compact_descent}), we may use the canonical map \[\mathrm{R}\underline{\mathrm{Hom}}_{(\mathcal{A},\mathcal{M})}(P,\mathcal{A}) \underset{(\mathcal{A},\mathcal{M})}{\overset{L}{\otimes}} (\mathcal{A}_U,\mathcal{M}_U)\rightarrow \mathrm{R}\underline{\mathrm{Hom}}_{(\mathcal{A}_U,\mathcal{M}_U)}(P \underset{(\mathcal{A},\mathcal{M})}{\overset{L}{\otimes}} (\mathcal{A}_U,\mathcal{M}_U),\mathcal{A}_U),\]
to construct the required mapping $\mathcal{A}_U\rightarrow (P \underset{(\mathcal{A},\mathcal{M})}{\overset{L}{\otimes}} (\mathcal{A}_U,\mathcal{M}_U))^{\vee}\overset{L}{\underset{(\mathcal{A}_U,\mathcal{M}_U)}{\otimes}}(Q \underset{(\mathcal{A},\mathcal{M})}{\overset{L}{\otimes}} (\mathcal{A}_U,\mathcal{M}_U))$.
\end{remark}

\begin{definition}[\cite{Analytic}, Definition 13.12]
An object $M\in \mathcal{D}(\mathcal{A},\mathcal{M})$ is called \textit{basic nuclear}, if it can be written as the colimit of a diagram of the form
\[ P_0\xrightarrow{f_0} P_1 \xrightarrow[]{f_1} P_2\xrightarrow[]{f_2}\cdots, \]
where for each $i\geq 0$, $P_i$ is compact and $f_i$ is trace class. 
\end{definition}

Let us record the following general properties of nuclear objects for later use:

\begin{proposition}[\cite{Analytic}, Proposition 13.13]
\label{nuclear colimits}
An object $M\in \mathcal{D}(\mathcal{A},\mathcal{M})$ is nuclear if and only if it can be written as a filtered colimit of basic nuclear objects. In particular, every basic nuclear object is nuclear. Furthermore, the class of basic nuclear objects is stable under all countable colimits, and the class of nuclear objects is stable under all colimits. 
\end{proposition}

\begin{proposition}[\cite{Analytic}, Proposition 13.14]
\label{extra nuclear}
Let $M\in \mathcal{D}(\mathcal{A},\mathcal{M})$ be a nuclear object. Then for all $N\in \mathcal{D}(\mathcal{A},\mathcal{M})$ and extremally disconnected sets $S$, the canonical map 
\[(\mathrm{R}\underline{\mathrm{Hom}}_{(\mathcal{A},\mathcal{M})}(\mathcal{M}[S],N)\overset{L}{\underset{(\mathcal{A},\mathcal{M})}{\otimes}} M)(*)\rightarrow (N\overset{L}{\underset{(\mathcal{A},\mathcal{M})}{\otimes}} M)(S)\]
is an isomorphism in $\mathcal{D}(\mathrm{Ab})$.
\end{proposition}

There is also an obvious ``internal'' version of the definition of nuclearity. However, the following proposition ensures that in the present context, it is actually equivalent to the introduced one:

\begin{proposition}
\label{nuclinternal}
Let $M\in \mathcal{D}(\mathcal{A},\mathcal{M})$ be a nuclear object. Then for every profinite set $S$ and every $N\in \mathcal{D}(\mathcal{A},\mathcal{M})$, the natural mapping \[\mathrm{R}\underline{\mathrm{Hom}}_{(\mathcal{A},\mathcal{M})}(\mathcal{M}[S],N)\overset{L}{\underset{(\mathcal{A},\mathcal{M})}{\otimes}} M\rightarrow \mathrm{R}\underline{\mathrm{Hom}}_{(\mathcal{A},\mathcal{M})}(\mathcal{M}[S],N\overset{L}{\underset{(\mathcal{A},\mathcal{M})}{\otimes}} M)\] is an isomorphism in $\mathcal{D}(\mathcal{A},\mathcal{M})$. In particular, for every profinite set $S$, the natural mapping \[{\mathcal{M}[S]^{\vee}\overset{L}{\underset{(\mathcal{A},\mathcal{M})}{\otimes}} M\rightarrow \mathrm{R}\underline{\mathrm{Hom}}_{(\mathcal{A},\mathcal{M})}(\mathcal{M}[S],M)}\] is an isomorphism.
\end{proposition}

\begin{proof}
We first show that for all profinite sets $S$, the natural map \[\mathrm{R}\mathrm{Hom}_{(\mathcal{A},\mathcal{M})}(\mathcal{A},\mathrm{R}\underline{\mathrm{Hom}}_{(\mathcal{A},\mathcal{M})}(\mathcal{M}[S],N)\overset{L}{\underset{(\mathcal{A},\mathcal{M})}{\otimes}} M)\rightarrow \mathrm{R}\mathrm{Hom}_{(\mathcal{A},\mathcal{M})}(\mathcal{M}[S], N\overset{L}{\underset{(\mathcal{A},\mathcal{M})}{\otimes}} M)\] is an isomorphism in $\mathcal{D}(\mathrm{Ab})$. As $\mathcal{M}[S]$ is projective, it is a direct summand of $\mathcal{M}[Q]$ for some extremally disconnected set $Q$. Applying Proposition~\ref{extra nuclear}, we obtain an isomorphism \[\mathrm{R}\mathrm{Hom}_{(\mathcal{A},\mathcal{M})}(\mathcal{A},\mathrm{R}\underline{\mathrm{Hom}}_{(\mathcal{A},\mathcal{M})}(\mathcal{M}[Q],N)\overset{L}{\underset{(\mathcal{A},\mathcal{M})}{\otimes}} M)\xrightarrow{\sim} \mathrm{R}\mathrm{Hom}_{(\mathcal{A},\mathcal{M})}(\mathcal{M}[Q],  N\overset{L}{\underset{(\mathcal{A},\mathcal{M})}{\otimes}} M),\]
which evidently implies the desired one. 

To prove the proposition, it suffices to show that for every extremally disconnected set $Q$, the natural map \[\mathrm{R}{\mathrm{Hom}}_{(\mathcal{A},\mathcal{M})}(\mathcal{M}[Q],\mathrm{R}\underline{\mathrm{Hom}}_{(\mathcal{A},\mathcal{M})}(\mathcal{M}[S],N)\overset{L}{\underset{(\mathcal{A},\mathcal{M})}{\otimes}} M)\rightarrow\mathrm{R}{\mathrm{Hom}}_{(\mathcal{A},\mathcal{M})}(\mathcal{M}[Q], \mathrm{R}\underline{\mathrm{Hom}}_{(\mathcal{A},\mathcal{M})}(\mathcal{M}[S], N\overset{L}{\underset{(\mathcal{A},\mathcal{M})}{\otimes}} M))\] is an isomorphism in $\mathcal{D}(\mathrm{Ab})$. Applying Proposition~\ref{extra nuclear} again and using Lemma~\ref{diffdescofdual}, we see that the left hand side is isomorphic to $(\mathrm{R}\underline{\mathrm{Hom}}_{(\mathcal{A},\mathcal{M})}(\mathcal{M}[S\times Q],N)\overset{L}{\underset{(\mathcal{A},\mathcal{M})}{\otimes}} M)(*)$. 
By Lemma~\ref{diffdescofdual}, we can rewrite the right hand side as ${\mathrm{R}{\mathrm{Hom}}_{(\mathcal{A},\mathcal{M})}(\mathcal{M}[S\times Q],N\overset{L}{\underset{(\mathcal{A},\mathcal{M})}{\otimes}} M)}$. The desired result now follows by the claim above.
\end{proof}

We now want to describe the relationship between nuclear and dualizable objects. For this, we need the following result: 

\begin{lemma}
\label{compact}
Let $M$ be an object of $\mathcal{D}(\mathcal{A},\mathcal{M})$. Then $M$ is compact if an only if it is a retract of a finite complex with terms of the form $\underset{i=1}{\overset{n}{\bigoplus}}\mathcal{M}[S_i]$ for some non-negative integer $n$ and profinite sets $S_i$.
\end{lemma}

\begin{proof}
See \cite[\href{https://stacks.math.columbia.edu/tag/094B}{Tag 094B}]{stacks-project}.
\end{proof}

We now record the following general observation regarding the relationship between nuclear and dualizable objects:

\begin{proposition}
\label{dual}
An object $M\in \mathcal{D}(\mathcal{A},\mathcal{M})$ is dualizable if and only if it is compact and nuclear.
\end{proposition}

\begin{proof}
It is a general fact about closed symmetric monoidal categories that if the tensor unit is compact, then dualizable objects are compact as well. Indeed, we compute:

\[\underline{\mathrm{Hom}}(M,\underset{i\in I}{\bigoplus} N_i)\cong \underline{\mathrm{Hom}}(1,M^{\vee}\otimes\underset{i\in I}{\bigoplus} N_i)\cong \underline{\mathrm{Hom}}(1,\underset{i\in I}{\bigoplus}(M^{\vee}\otimes N_i))\]\[\cong \underset{i\in I}{\bigoplus}\underline{\mathrm{Hom}}(1,M^{\vee}\otimes N_i)\cong \underset{i\in I}{\bigoplus}\underline{\mathrm{Hom}}(M,N_i).\]
Therefore, every dualizable object is in fact basic nuclear. Indeed, it readily follows from the axioms that the identity map $\mathrm{id}_M:M\rightarrow M$ is trace class.

For the reverse implication, assume that $M$ is nuclear and compact. Nuclearity implies that for every compact object $P$, we have $(P^{\vee}\underset{(\mathcal{A},\mathcal{M})}{\otimes}M)(*)\cong \mathrm{RHom}_{(\mathcal{A},\mathcal{M})}(P,M)$ (easily verified by induction using Lemma~\ref{compact}). Invoking the isomorphism \[(P^{\vee}\underset{(\mathcal{A},\mathcal{M})}{\otimes}M)(*)\cong \mathrm{RHom}_{(\mathcal{A},\mathcal{M})}(\mathcal{A},P^{\vee}\underset{\mathcal{A}}{\otimes}M)\] and substituting $P=M$, we obtain a mapping $i_M:\mathcal{A}\rightarrow M^{\vee}\underset{(\mathcal{A},\mathcal{M})}{\overset{L}{\otimes}}M$, which corresponds to the identity mapping $\mathrm{id}_M:M\rightarrow M$. It is straightforward to check that this map, together with the obvious evaluation map $\mathrm{ev}_M:M^{\vee}\underset{(\mathcal{A},\mathcal{M})}{\overset{L}{\otimes}}M\rightarrow \mathcal{A}$, satisfies the axioms.
\end{proof}

We now prove the following property of the functor $\mathrm{R}\underline{\mathrm{Hom}}_{(\mathcal{A},\mathcal{M})}(-, -)$, which will be essential for our proofs later on; it can be regarded as ``nuclearity'' of the structure sheaf of $X$:

\begin{proposition}
\label{basicdualbasechange}
Let $U$ be an affinoid open subspace of $X$. Then for every profinite set $S$ and object $M\in \mathcal{D}(\mathcal{A},\mathcal{M})$, we have
\[\mathrm{R}\underline{\mathrm{Hom}}_{(\mathcal{A},\mathcal{M})}(\mathcal{M}[S], M)\overset{L}{\underset{(\mathcal{A},\mathcal{M})}{\otimes}}(\mathcal{A}_U,\mathcal{M}_U)\xrightarrow{\sim}\mathrm{R}\underline{\mathrm{Hom}}_{(\mathcal{A}_U,\mathcal{M}_U)}(\mathcal{M}_{U}[S], M\overset{L}{\underset{(\mathcal{A},\mathcal{M})}{\otimes}}(\mathcal{A}_U,\mathcal{M}_U))\]
\end{proposition}

\begin{proof}
Suppose first that $U$ is a nice rational open subspace. Arguing by induction, we easily reduce the proof to the case where $U$ is either $X(\frac{1}{f})$ or $X(\frac{f}{1})$ for some $f\in A$. Let us discuss only the former of the two cases, as the other is completely analogous (and even a little simpler). Forgetting the $\mathcal{A}_{U}$-structures on both sides and using Lemma~\ref{basechange} and Lemma~\ref{diffdescofdual}, we are reduced to showing that the map
\begin{equation}
\label{woof}
\mathrm{R}\underline{\mathrm{Hom}}_{R^+}((R^+,\mathbb{Z})_{\blacksquare}[S], M)\overset{L}{\underset{(R^+,\mathbb{Z})_{\blacksquare}}{\otimes}}(R,R^-)_{\blacksquare}\xrightarrow{\sim}\mathrm{R}\underline{\mathrm{Hom}}_{R^+}((R^+,\mathbb{Z})_{\blacksquare}[S], M\overset{L}{\underset{(R^+,\mathbb{Z})_{\blacksquare}}{\otimes}}(R,R^-)_{\blacksquare})~\tag{$*$}
\end{equation}
where $R^+=\mathbb{Z}[T]$, $R^-=\mathbb{Z}[T^{-1}]$, and $R=\mathbb{Z}[T,T^{-1}]$, is an isomorphism. Observe that the following sequence of condensed $R^+$-modules is exact:
\[0\rightarrow M[T]\xrightarrow{T-f} M[T]\rightarrow M \rightarrow 0.\]
Therefore, it suffices to prove the statement for $M[T]$. Using the isomorphism
\[\mathrm{R}\underline{\mathrm{Hom}}_{R^+}((R^+,\mathbb{Z})_{\blacksquare}[S], M[T])\xrightarrow{\sim} \mathrm{R}\underline{\mathrm{Hom}}_{\mathbb{Z}}(\mathbb{Z}_{\blacksquare}[S], M)[T],\]
we rewrite the left hand side of (\ref{woof}) as
\[((\mathrm{R}\underline{\mathrm{Hom}}_{\mathbb{Z}}(\mathbb{Z}_{\blacksquare}[S], M)\overset{L}{\underset{\mathbb{Z}_{\blacksquare}}{\otimes}}(R^-,\mathbb{Z})_{\blacksquare})\overset{L}{\underset{(R^-,\mathbb{Z})_{\blacksquare}}{\otimes}}R^-_{\blacksquare})\overset{L}{\underset{R^-_{\blacksquare}}{\otimes}}(R,R^-)_{\blacksquare},\]
and the right hand side as
\[\mathrm{R}\underline{\mathrm{Hom}}_{\mathbb{Z}}(\mathbb{Z}_{\blacksquare}[S],((M\overset{L}{\underset{\mathbb{Z}_{\blacksquare}}{\otimes}}(R^-,\mathbb{Z})_{\blacksquare})\overset{L}{\underset{(R^-,\mathbb{Z})_{\blacksquare}}{\otimes}}R^-_{\blacksquare})\overset{L}{\underset{R^-_{\blacksquare}}{\otimes}}(R,R^-)_{\blacksquare}).\]
The desired result now follows by combining Lemma~\ref{polynomialization}, Lemma~\ref{localizationT}, and Proposition~\ref{completion derived}.

In the general case, we note that $U$ admits a finite open covering $\{U_i\}_{i=1}^{n}$ by nice open subspaces of $X$ (which are also nice open subspaces of $U$) by virtue of Proposition~\ref{kedlnice}. As the functor
\[\mathcal{D}(\mathcal{A}_U,\mathcal{M}_U)\rightarrow \underset{{i=1}}{\overset{n}{\prod}} \mathcal{D}(\mathcal{A}_{U_{i}},\mathcal{M}_{U_i})\]
is conservative (by Theorem~\ref{maindescent}), the desired statement clearly follows from the above.
\end{proof}

We will also need the following property of the localization functors:

\begin{proposition}
\label{tordimension}
Let $U$ be an affinoid open subspace of $X$. Then the functor $-\underset{(\mathcal{A},\mathcal{M})}{\overset{L}{\otimes}}(\mathcal{A}_{U},\mathcal{M}_{U})$ has finite Tor-dimension. That is, there exists an integer $k\in \mathbb{N}$, such that for every $n\in\mathbb{Z}$ and every $M\in \mathcal{D}^{\geq n}(\mathcal{A},\mathcal{M})$, the module $M\underset{(\mathcal{A},\mathcal{M})}{\overset{L}{\otimes}}(\mathcal{A}_{U},\mathcal{M}_{U})$ lies in $D^{\geq n-k}(\mathcal{A}_{U},\mathcal{M}_{U})$.
\end{proposition}

\begin{proof}
In the case where $U$ is nice rational open subspace, the statement follows directly from Propositions~\ref{completion derived} and \ref{basechange}(ii). In general, $U$ admits a finite cover by nice rational open subspaces by virtue of Proposition~\ref{kedlnice}, and, therefore, the desired statement follows from Theorem~\ref{maindescent}. 
\end{proof}

We now record the following lemma for later use:

\begin{lemma}
\label{cone}
Let $X\rightarrow Y \rightarrow Z\rightarrow X[1]$ be a triangle in $\mathcal{D}(\mathcal{A},\mathcal{M})$. If any two of its three objects are compact (resp. nuclear, resp. dualizable, resp. pseudocoherent), then the remaining third one is compact (resp. nuclear, resp. dualizable, resp. pseudocoherent) as well.
\end{lemma}

\begin{proof}
The claim for pseudocoherent objects is obvious and follows directly from the definition. The compact case readily follows from Lemma~\ref{compact}. Since the functors used in Definition~\ref{definition nuclear} are exact, we obtain the claim for nuclear objects. Finally, the claim for dualizable objects follows by Lemma~\ref{dual}.
\end{proof}

We now prove that the four discussed classes of objects (i.e., nuclear,compact, dualizable, and pseudocoherent) satisfy descent. That is, they form a subsheaf of the sheaf of $\infty$-categories on $X$ introduced in Theorem~\ref{maindescent}.

\begin{theorem}
\label{compact_descent}
Compact objects satisfy descent on $X$.
\end{theorem}

\begin{proof}
Let $M\in \mathcal{D}(\mathcal{A},\mathcal{M})$ be a compact object and $U$ a affinoid open subspace of $X$. To prove that the localization $M\underset{(\mathcal{A},\mathcal{M})}{\overset{L}{\otimes}}(\mathcal{A}_U,\mathcal{M}_U)$ is compact, we compute:
\[\mathrm{R}{\mathrm{Hom}}_{(\mathcal{A}_U,\mathcal{M}_U)}(M\underset{(\mathcal{A},\mathcal{M})}{\overset{L}{\otimes}}(\mathcal{A}_U,\mathcal{M}_U),\underset{i\in I}{\bigoplus} N_i)\cong \mathrm{R}{\mathrm{Hom}}_{(\mathcal{A},\mathcal{M})}(M,\underset{i\in I}{\bigoplus} N_i) \]\[\cong \underset{i\in I}{\bigoplus} \mathrm{R}{\mathrm{Hom}}_{(\mathcal{A},\mathcal{M})}(M,N_i) \cong \underset{i\in I}{\bigoplus}\mathrm{R}{\mathrm{Hom}}_{(\mathcal{A}_U,\mathcal{M}_U)}(M\underset{(\mathcal{A},\mathcal{M})}{\overset{L}{\otimes}}(\mathcal{A}_U,\mathcal{M}_U),N_i).\]

Now let $\mathcal{U}=\{U_k\}_{k=1}^{n}$ be an affinoid open cover of $X$ and $M$ an object of $\mathcal{D}(\mathcal{A},\mathcal{M})$ such that for each $k=1,\dots,n$, the localization $M_k:=M\underset{(\mathcal{A},\mathcal{M})}{\overset{L}{\otimes}}(\mathcal{A}_{U_k},\mathcal{M}_{U_k})$ is compact. Applying Proposition~\ref{kedlnice}, we may assume that the covering is nice. Arguing by induction, we immediately reduce to the case where $n=2$. Denote the localization $M\underset{(\mathcal{A},\mathcal{M})}{\overset{L}{\otimes}}(\mathcal{A}_{U_1\cap U_2},\mathcal{M}_{U_1\cap U_2})$ by $M_{12}$. Suppose we are given a set $I$ and a collection $\{N_i\in \mathcal{D}(\mathcal{A},\mathcal{M})$| $i\in I\}$. For each $i\in I$, denote the localizations  $N_i\underset{(\mathcal{A},\mathcal{M})}{\overset{L}{\otimes}}(\mathcal{A}_{U_1},\mathcal{M}_{U_1})$,  $N_i\underset{(\mathcal{A},\mathcal{M})}{\overset{L}{\otimes}}(\mathcal{A}_{U_2},\mathcal{M}_{U_2})$, and  $N_i\underset{(\mathcal{A},\mathcal{M})}{\overset{L}{\otimes}}(\mathcal{A}_{U_1\cap U_2},\mathcal{M}_{U_1\cap U_2})$ by $N_{i,1}$, $N_{i,2}$, and $N_{i,12}$ respectively. We need to show that the canonical map \[\Phi:\underset{i\in I}{\bigoplus}\mathrm{Hom}_{\mathcal{D}(\mathcal{A},\mathcal{M})}(M, N_i)\rightarrow \mathrm{Hom}_{\mathcal{D}(\mathcal{A},\mathcal{M})}(M,\underset{i\in I}{\bigoplus} N_i)\] is an isomorphism. For every $i\in I$, we have the following exact sequence:
\[  \mathrm{Hom}_{\mathcal{D}(\mathcal{A},\mathcal{M})}(M,N_{i,1}[-1])\oplus \mathrm{Hom}_{\mathcal{D}(\mathcal{A},\mathcal{M})}(M, N_{i,2}[-1]) \rightarrow \mathrm{Hom}_{\mathcal{D}(\mathcal{A},\mathcal{M})}(M,N_{i,12}[-1]) \rightarrow\]\[\rightarrow \mathrm{Hom}_{\mathcal{D}(\mathcal{A},\mathcal{M})}(M, N_i)\rightarrow \mathrm{Hom}_{\mathcal{D}(\mathcal{A},\mathcal{M})}(M,N_{i,1})\oplus \mathrm{Hom}_{\mathcal{D}(\mathcal{A},\mathcal{M})}(M, N_{i,2})\rightarrow \mathrm{Hom}_{\mathcal{D}(\mathcal{A},\mathcal{M})}(M, N_{i,12}) \]
As direct sums are exact, the following sequence is also exact:
\[  \underset{i\in J}{\bigoplus}(\mathrm{Hom}_{\mathcal{D}(\mathcal{A},\mathcal{M})}(M,N_{i,1}[-1])\oplus \mathrm{Hom}_{\mathcal{D}(\mathcal{A},\mathcal{M})}(M, N_{i,2}[-1]))\rightarrow\underset{i\in J}{\bigoplus} \mathrm{Hom}_{\mathcal{D}(\mathcal{A},\mathcal{M})}(M, N_{i,12}[-1]) \rightarrow\]\[\rightarrow \underset{i\in J}{\bigoplus}\mathrm{Hom}_{\mathcal{D}(\mathcal{A},\mathcal{M})}(M, N_i)\rightarrow \underset{i\in J}{\bigoplus}( \mathrm{Hom}_{\mathcal{D}(\mathcal{A},\mathcal{M})}(M,N_{i,1})\oplus \mathrm{Hom}_{\mathcal{D}(\mathcal{A},\mathcal{M})}(M,N_{i,2}))\rightarrow \underset{i\in J}{\bigoplus} \mathrm{Hom}_{\mathcal{D}(\mathcal{A},\mathcal{M})}(M,N_{i,12})\]
Furthermore, the sequence
\[ \mathrm{Hom}_{\mathcal{D}(\mathcal{A},\mathcal{M})}(M,\underset{i\in J}{\bigoplus}N_{i,1}[-1])\oplus \mathrm{Hom}_{\mathcal{D}(\mathcal{A},\mathcal{M})}(M, \underset{i\in J}{\bigoplus}N_{i,2}[-1])\rightarrow \mathrm{Hom}_{\mathcal{D}(\mathcal{A},\mathcal{M})}(M,\underset{i\in J}{\bigoplus} N_{i,12}[-1]) \rightarrow\]\[\rightarrow \mathrm{Hom}_{\mathcal{D}(\mathcal{A},\mathcal{M})}(M,\underset{i\in J}{\bigoplus} N_i)\rightarrow  \mathrm{Hom}_{\mathcal{D}(\mathcal{A},\mathcal{M})}(M,\underset{i\in J}{\bigoplus}N_{i,1})\oplus \mathrm{Hom}_{\mathcal{D}(\mathcal{A},\mathcal{M})}(M,\underset{i\in J}{\bigoplus}N_{i,2})\rightarrow  \mathrm{Hom}_{\mathcal{D}(\mathcal{A},\mathcal{M})}(M,\underset{i\in J}{\bigoplus}N_{i,12})\]
is also exact. As $M$ is locally compact and the functors \[\mathcal{D}(\mathcal{A},\mathcal{M})\rightarrow \mathcal{D}(\mathcal{A}_{U_1},\mathcal{M}_{U_1}),\ \mathcal{D}(\mathcal{A},\mathcal{M})\rightarrow \mathcal{D}(\mathcal{A}_{U_2},\mathcal{M}_{U_2}),\ \mathrm{and}\ \mathcal{D}(\mathcal{A},\mathcal{M})\rightarrow \mathcal{D}(\mathcal{A}_{U_1\cap U_2},\mathcal{M}_{U_1\cap U_2})\] are all fully faithful (by Proposition~\ref{basechange}), we see that $\Phi$ is an isomorphism by virtue of the five lemma.
\end{proof}

\begin{theorem}
\label{nuclear_descent}
Nuclear objects satisfy descent on $X$.
\end{theorem}

\begin{proof}
Observe that localization preserves the class of nuclear objects. Indeed, let $U\subset X$ be a affinoid open subset. By virtue of Proposition~\ref{nuclear colimits} and the fact that localization commutes with colimits, it suffices to show that $-\underset{(\mathcal{A},\mathcal{M})}{\overset{L}{\otimes}}(\mathcal{A}_U,\mathcal{M}_U)$ preserves basic nuclear objects. The desired claim now follows by observing that $-\underset{(\mathcal{A},\mathcal{M})}{\overset{L}{\otimes}}(\mathcal{A}_U,\mathcal{M}_U)$ preserves compact objects and trace class maps.

Now let $\mathcal{U}=\{U_i\}_{i=1}^{n}$ be an affinoid open cover of $X$ and $M$ an object of $\mathcal{D}(\mathcal{A},\mathcal{M})$ such that for each $i=1,\dots,n$, the localization $M_i:=M\underset{(\mathcal{A},\mathcal{M})}{\overset{L}{\otimes}}(\mathcal{A}_{U_i},\mathcal{M}_{U_i})$ is nuclear. Applying Proposition~\ref{basicdualbasechange}, we deduce that for each $i=1,\cdots,n$, the canonical map 
\[(\mathcal{M}[S]^{\vee}\overset{L}{\underset{(\mathcal{A},\mathcal{M})}{\otimes}} M)\overset{L}{\underset{(\mathcal{A},\mathcal{M})}{\otimes}}(\mathcal{A}_{U_i},\mathcal{M}_{U_i})\xrightarrow{\sim}\mathcal{M}_{U}[S]^{\vee}\overset{L}{\underset{(\mathcal{A}_{U_i},\mathcal{M}_{U_i})}{\otimes}} (M\overset{L}{\underset{(\mathcal{A},\mathcal{M})}{\otimes}}(\mathcal{A}_{U_i},\mathcal{M}_{U_i}))\]
is an isomorphism. Using Proposition~\ref{nuclinternal}, we obtain an isomorphism 
\[\mathcal{M}_{U_i}[S]^{\vee}\overset{L}{\underset{(\mathcal{A}_{U_i},\mathcal{M}_{U_i})}{\otimes}} (M\overset{L}{\underset{(\mathcal{A},\mathcal{M})}{\otimes}}(\mathcal{A}_{U_i},\mathcal{M}_{U_i}))\xrightarrow{\sim}\mathrm{R}\underline{\mathrm{Hom}}_{(\mathcal{A}_{U_i},\mathcal{M}_{U_i})}(\mathcal{M}_{U_i}[S], M\overset{L}{\underset{(\mathcal{A},\mathcal{M})}{\otimes}}(\mathcal{A}_{U_i},\mathcal{M}_{U_i})).\] Applying Proposition~\ref{basicdualbasechange} again, we see that the canonical map 
\[\mathrm{R}\underline{\mathrm{Hom}}_{(\mathcal{A},\mathcal{M})}(\mathcal{M}[S], M)\overset{L}{\underset{(\mathcal{A},\mathcal{M})}{\otimes}}(\mathcal{A}_{U_i},\mathcal{M}_{U_i})\xrightarrow{\sim}\mathrm{R}\underline{\mathrm{Hom}}_{(\mathcal{A}_{U_i},\mathcal{M}_{U_i})}(\mathcal{M}_{U_i}[S], M\overset{L}{\underset{(\mathcal{A},\mathcal{M})}{\otimes}}(\mathcal{A}_{U_i},\mathcal{M}_{U_i}))\]
is an isomorphism. Therefore, for each $i=1,\cdots,n$, the restriction of the natural map \[\phi:{\mathcal{M}[S]^{\vee}\overset{L}{\underset{(\mathcal{A},\mathcal{M})}{\otimes}} M\rightarrow \mathrm{R}\underline{\mathrm{Hom}}_{(\mathcal{A},\mathcal{M})}(\mathcal{M}[S],M)}\] to $U_i$ is an isomorphism. As the functor
\[\mathcal{D}(\mathcal{A},\mathcal{M})\rightarrow \underset{i=1}{\overset{n}{\prod}}\mathcal{D}(\mathcal{A}_{U_i},\mathcal{M}_{U_i})\] is conservative, $\phi$ is an isomorphism itself, as desired. 
\end{proof}

Combining Proposition~\ref{dual}, Theorem~\ref{compact_descent}, and Theorem~\ref{nuclear_descent} we obtain the following:

\begin{theorem}\label{dual_descent}
Dualizable objects satisfy descent on $X$.
\end{theorem}

Finally, we prove the following:

\begin{theorem}\label{pseudocoherent objects descent}
Pseudocoherent objects (of $\mathcal{D}^-(\mathcal{A},\mathcal{M})$) satisfy descent on $X$.
\end{theorem}

\begin{proof}
The fact that localization preserves the class of pseudocoherent objects is trivial and follows directly from the definition.

Now let $\mathcal{U}=\{U_k\}_{k=1}^{n}$ be a affinoid open cover of $X$ and $M$ an object of $\mathcal{D}(\mathcal{A},\mathcal{M})$ such that for each $k=1,\dots,n$, the localization $M_k:=M\underset{(\mathcal{A},\mathcal{M})}{\overset{L}{\otimes}}(\mathcal{A}_{U_k},\mathcal{M}_{U_k})$ is pseudocoherent. Applying Proposition~\ref{kedlnice}, we may assume that the covering is nice. Arguing by induction, we immediately reduce to the case where $n=2$. As $M$ is isomorphic to the homotopy pullback of $M_1\rightarrow M_{12}\xleftarrow[]{} M_2$, we see that $M$ is bounded above. Suppose we are given an integer $j$, a directed set $I$, and a collection $\{N_i\in \mathcal{D}^{\geq j}(\mathcal{A},\mathcal{M})$| $i\in I\}$. We need to show that the canonical map \[\underset{i\in I}{\colim}(\mathrm{Hom}_{\mathcal{D}(\mathcal{A},\mathcal{M})}(M, N_i))\rightarrow \mathrm{Hom}_{\mathcal{D}(\mathcal{A},\mathcal{M})}(M,\underset{i\in I}{\colim} N_i)\] is an isomorphism. Invoking Proposition~\ref{tordimension}, we see that there exists an integer $k\in\mathbb{N}$, such that for $U\in \{U_1,U_2,U_1\cap U_2\}$, the localization $N_i\underset{(\mathcal{A},\mathcal{M})}{\overset{L}{\otimes}}(\mathcal{A}_{U},\mathcal{M}_{U})$ is an object of $\mathcal{D}^{\geq j-k}(\mathcal{A}_U,\mathcal{M}_U)$. Arguing verbatim as in the proof of Theorem~\ref{compact_descent}, we obtain the desired isomorphism.
\end{proof}

We conclude this subsection by various remarks regarding the relationship between discrete objects on the one hand and compact, dualizable and nuclear on the other. 

\begin{lemma}
\label{discnucl}
Every discrete object of $\mathcal{D}(\mathcal{A},\mathcal{M})$ is nuclear.
\end{lemma}

\begin{proof}
The claim immediately follows from the fact that $\underline{A}$ is nuclear (by Proposition~\ref{dual}, as it is evidently dualizable) and the class of nuclear objects is stable under colimits (Proposition~\ref{nuclear colimits}).
\end{proof}

\begin{lemma}
\label{perf}
An object $M$ of $\mathcal{D}(A)$ is dualizable if and only if it is compact if and only if it is a perfect complex.
\end{lemma}

\begin{proof}
For the second equivalence, see \cite[\href{https://stacks.math.columbia.edu/tag/0656}{Tag 0656}]{stacks-project}. The implication ``dualizable $\implies$ compact'' is a part of Proposition~\ref{dual}. Therefore, it remains to show that every compact object of $\mathcal{D}(A)$ is dualizable. Using the second equivalence, we may assume that $M$ is a perfect complex; in this case, one can easily construct the required map $i_X$ explicitly.
\end{proof}

\begin{remark}
Observe that the condensification functor sends compact objects (i.e., perfect complexes) to compact ones.
\end{remark}

\begin{lemma}
\label{two_duals}
Let $M$ be a dualizable and discrete object of $\mathcal{D}(\mathcal{A},\mathcal{M})$. Then it is dualizable in $\mathcal{D}(A)$.
\end{lemma}

\begin{proof}
Let $X$ be an object of $\mathcal{D}(A)$ such that $\mathrm{Cond}(X)$ is isomorphic to $M$. As $M$ is compact (by Proposition~\ref{dual}) and the condensification functor is fully faithful and commutes with filtered colimits, we conclude that $X$ is compact in $\mathcal{D}(A)$. The desired claim now follows by Lemma~\ref{perf}.
\end{proof}
\subsection{Proofs}\label{vb4}

We are finally ready to begin the proofs of Theorems~\ref{descent of everything} and \ref{pseudo stack}. We will need the following object:

\begin{lemma}
Denote by $NS_{(\mathcal{A},\mathcal{M})}$ ( ``null sequences'') the module $\mathcal{M}[\mathbb{N}\cup \{\infty\}]/\mathcal{A}$, where the inclusion $\mathcal{A}\hookrightarrow \mathcal{M}[\mathbb{N}\cup \{\infty\}]$ is induced by the map $\{*\}\rightarrow \mathbb{N}\cup \{\infty\},\ *\mapsto \infty$. Then $NS_{(\mathcal{A},\mathcal{M})}$ is compact, projective, and isomorphic to $\colim_M \underset{\mathbb{N}}{\prod} \underline{M}$, where the colimit is taken over all finitely generated subrings $B\subset A^+$ and all quasi-finitely generated $B$-submodules $M$ of $A$.
\end{lemma}
\begin{proof}
Using the argument of Lemma~\ref{freeguyshuber}, the claim is easily derived from Lemma~\ref{powerseries} and its proof.
\end{proof}

We also need the following technical result:

\begin{lemma}
\label{contfact}
Let $S$ and $T$ be profinite sets and $f$ a map $\mathcal{M}[T]\rightarrow \underline{C(S,A)}$. Then $f$ can be factored through $NS_{(\mathcal{A},\mathcal{M})}$. 
\end{lemma}

\begin{proof}
Fix a pair of definition $(A_0,I)$ of $A$ and denote by $\phi$ the map $T\rightarrow C(S,A)$ which gives rise to $f:\mathcal{M}[T]\rightarrow \underline{C(S,A)}$. Note that $C(S,A)$ is homeomorphic to $\varprojlim C(S,A/I^n)$. As $S$ is compact and $A/I^n$ is discrete, each space $C(S,A/I^n)$ is discrete as well. Therefore, each induced map $\phi_n: T\rightarrow C(S,A/I^n)$ factors as
\begin{center}
    \begin{tikzcd}
 T \arrow[r, "\phi_n"] \arrow[d,swap, two heads,"\pi_n"]
    & C(S,A/I^n) \\
   T_n \arrow[ru, hook, "\varphi_n"]
    \end{tikzcd}
\end{center}
where $T_n$ is a finite discrete set, $\pi_n$ is surjective, and $\varphi_n$ is injective. Each element $t\in T_n$ can be identified with its image in $C(S,A/I^n)$. Pick arbitrary lifts $f_{0,1},\dots, f_{0,|T_1|}\in C(S,A)$ of the elements of $T_1\subset C(S,A/I)$. Arguing inductively, we construct for each $n>0$, elements $f_{n,1},\dots, f_{n,|T_n|}\in C(A,I^{n-1})$ such that each $t\in T_n$ may be written as $\underset{j=1}{\overset{n}{\sum}}f_{j,i_j}$ modulo $I^n$ for certain $i_j\in \{1,\dots, |T_j|\}$. Now consider the map $\tilde{f}:NS_{(\mathcal{A},\mathcal{M})} \rightarrow C(S,A)$ induced by the following sequence $\mathbb{N}\cup \{\infty\}\rightarrow C(S,A)$:
\[(x_i)_{i\in \mathbb{N}}=(f_{1,1},f_{1,2},\dots,f_{1,|T_1|}, f_{2,1},\dots, f_{2,|T_2|}, \dots).\] We claim that $f$ can be factored as $\tilde{f}\circ g$ for some $g: \mathcal{M}[T]\rightarrow NS_{(\mathcal{A},\mathcal{M})}$. Write $NS_{(\mathcal{A},\mathcal{M})}$ as $\colim_M \underset{\mathbb{N}}{\prod} \underline{M}$; for each $t\in T$, let $(t_i)_{i\in \mathbb{N}}$ be an element of $NS_{(\mathcal{A},\mathcal{M})}$ such that every $t_i$ is either $0$ or $1$ and $\phi(t)=\underset{i\in\mathbb{N}}{\sum}t_ix_i$ holds. It is clear that if the map $g':T\rightarrow \underset{\mathbb{N}}{\prod} A,\ t\mapsto (t_i)_{i\in \mathbb{N}}$ is continuous, then it provides the required factorization. To prove continuity, we are required to show that for each $k\geq 0$, the induced map \[T\rightarrow \underset{i=1}{\overset{A_K}{\prod}} A,\ A_k:=|T_1|+\cdots|T_k|\] is continuous. In fact, it can be easily shown that every such map factors through $\pi_k:T\rightarrow T_k$, which implies that it is continuous, as desired.
\end{proof}

\begin{theorem}
\label{d_is_d}
Let $M$ be a nuclear pseudocoherent object of $\mathcal{D}^-(\mathcal{A},\mathcal{M})$. Then $M$ is discrete. In particular, any dualizable object is discrete.
\end{theorem}
\begin{proof}
Without loss of generality, we assume that $M$ can be written as  
\[\ldots \rightarrow \mathcal{M}[S_{2}]\rightarrow \mathcal{M}[S_{1}]\rightarrow \mathcal{M}[S_0]\rightarrow 0, \]
where $\{S_i\}_{i\geq 0}$ are profinite sets. Since $M$ is nuclear, the canonical map $f':\mathcal{M}[S_0] \rightarrow M$ can be factored as the composition \[\mathcal{M}[S_0] \xrightarrow{1\otimes g'} \mathcal{M}[S_0] \overset{L}{\underset{(\mathcal{A},\mathcal{M})}{\otimes}}  \mathcal{M}[S_0]^{\vee} \overset{L}{\underset{(\mathcal{A},\mathcal{M})}{\otimes}} M  \xrightarrow{ev_{\mathcal{M}[S_0]}\otimes 1} M\]
for some $g': \mathcal{A}\rightarrow \mathcal{M}[S_0]^{\vee} \overset{L}{\underset{(\mathcal{A},\mathcal{M})}{\otimes}} M$. As $\mathcal{A}$ is projective, $g'$ can in turn be factored as the composition \[\mathcal{A}\xrightarrow{g} \mathcal{M}[S_0]^{\vee} \overset{L}{\underset{(\mathcal{A},\mathcal{M})}{\otimes}} \mathcal{M}[S_0] \xrightarrow{1\otimes f'} \mathcal{M}[S_0]^{\vee} \overset{L}{\underset{(\mathcal{A},\mathcal{M})}{\otimes}}M.\]
Thus, we have constructed a trace class map $f:\mathcal{M}[S_0]\rightarrow \mathcal{M}[S_0]$ such that the composition $f'\circ f$ is equal to $f'$. Therefore, the composition $f'\circ (1-f)$ is trivial, and, hence, we obtain the following morphism of triangles:
\begin{center}
\begin{tikzcd}[]
\mathcal{M}[S_0] \arrow{r}{1-f}\arrow{d} & \mathcal{M}[S_0] \arrow{r}\arrow{d}{f'} & \cone(1-f) \arrow{r}\arrow{d}{\phi}& \mathcal{M}[S_0][1] \arrow{d} \\
0 \arrow{r}{} & M \arrow{r}{\mathrm{id}_{M}} & M \arrow{r}& 0
\end{tikzcd}
\end{center}
It can be easily checked that the induced map $H^{0}(\cone(1-f))\xrightarrow{} H^0(M)$ is surjective. Applying Lemma~\ref{trace-class}, we conclude that $\cone(1-f_n)$ is a perfect complex over $A$. Consequently, $\cone\phi$ is pseudocoherent and nuclear (by Lemma~\ref{cone}). It is straightforward to check that the cohomology group $H^{0}(\cone(\phi))$ is trivial. Thus, we have the triangle
\begin{center}
\begin{tikzcd}[]
\cone\phi [-1] \arrow{r}{}& \cone(1-f_n) \arrow{r}{\phi} & M \arrow{r} & \cone\phi,
\end{tikzcd}
\end{center}
whose left hand side is pseudocoherent and nuclear. Arguing by induction, one shows that $M$ can be written as a colimit of discrete objects, and, hence, is discrete itself.
\end{proof}

\begin{lemma}
\label{trace-class}
Let $S$ be a profinite set and $f:\mathcal{M}[S]\xrightarrow{} \mathcal{M}[S]$ a trace-class morphism. Then $\cone(f-1)$ is quasi-isomorphic to a complex of the form $\underline{A}^n[-1]\rightarrow \underline{A}^{n}[0]$.
\end{lemma}

\begin{proof}
By Lemma~\ref{diffdescofdual}, we may identify $\mathcal{M}[S]^{\vee}$ with $\underline{C(S,A)}$. Let $g:\underline{A}\xrightarrow{} \underline{C(S,A)}\overset{L}{\underset{(\mathcal{A},\mathcal{M})}{\otimes}} \mathcal{M} [S]$ be a morphism in $\mathcal{D}(\mathcal{A},\mathcal{M})$ realizing $f$ as a trace-class map. It is easily verified that the set of maps $\{\underline{A}\xrightarrow{} \underline{C(S,A)}\overset{L}{\underset{(\mathcal{A},\mathcal{M})}{\otimes}} \mathcal{M} [S]\}$ in $\mathcal{D}(\mathcal{A},\mathcal{M})$ is in bijection with the set of maps $\{\underline{A}\rightarrow\underline{C(S,A)}\overset{}{\underset{(\mathcal{A},\mathcal{M})}{\otimes}} \mathcal{M} [S]\}$ in $\mathrm{Mod}_{(\mathcal{A},\mathcal{M})}^{\mathrm{cond}}$. It implies the map $f: \mathcal{M}[S]\rightarrow \mathcal{M}[S]$ can be factored (in $\mathrm{Mod}_{(\mathcal{A},\mathcal{M})}^{\mathrm{cond}}$) as 
\[\mathcal{M}[S]\xrightarrow{1\otimes g} \mathcal{M}[S] \overset{}{\underset{(\mathcal{A},\mathcal{M})}{\otimes}} \underline{C(S,A)}\overset{}{\underset{(\mathcal{A},\mathcal{M})}{\otimes}} \mathcal{M} [S] \xrightarrow{\mathrm{ev}_{\mathcal{M}[S]}\otimes 1} \mathcal{M}[S].\]
Pick be a surjective map $\tau':\underset{i\in I}{\bigoplus}\mathcal{M}[Q_i]\twoheadrightarrow \underline{C(S,A)}$, where $\{Q_i\}_{i\in I}$ are profinite sets. As $\underline{A}$ is compact, we obtain the following factorization of the map $g$: 
\[ \underline{A}\rightarrow \mathcal{M}[T\times S]\xrightarrow{\tau \otimes 1}\underline{C(S,A)}\overset{}{\underset{(\mathcal{A},\mathcal{M})}{\otimes}} \mathcal{M}[S],\]
where $T=\underset{j\in J}{\coprod} Q_j$ for some finite subset $J\subset I$ and $\tau$ is the restriction of $\tau'$ to $\mathcal{M}[T]$.

Applying Lemma~\ref{contfact}, we obtain a factorization $\mathcal{M}[T]\xrightarrow[]{\tilde{\tau}}NS_{(\mathcal{A},\mathcal{M})}\xrightarrow[]{\sigma} \underline{C(S,A)}$ of $\tau$; it gives rise to the following factorization:
\[\phi:\underline{A}\rightarrow NS_{(\mathcal{A},\mathcal{M})}\otimes \mathcal{M}[S],\ \psi:\mathcal{M}[S]\otimes NS_{(\mathcal{A},\mathcal{M})}\rightarrow \underline{A},\]
\[f: \mathcal{M}[S]\xrightarrow{\mathrm{id}\otimes\phi} \mathcal{M}[S]\otimes NS_{(\mathcal{A},\mathcal{M})}\otimes \mathcal{M}[S]\xrightarrow{\psi\otimes \mathrm{id}}\mathcal{M}[S].\]

Let us describe the image of $\phi$. Since $\mathcal{M}[\mathbb{N}\cup \{\infty\}]\underset{(\mathcal{A},\mathcal{M})}{\otimes}\mathcal{M}[S]$ is isomorphic to $\mathcal{M}[(\mathbb{N}\cup \{\infty\})\times S]$, we have $NS_{(\mathcal{A},\mathcal{M})}\underset{(\mathcal{A},\mathcal{M})}{\otimes} \mathcal{M}[S]=\colim_M \underset{\mathbb{N}\times J}{\prod} \underline{M}$, where the colimit is taken over all finitely generated subrings $B\subset A^+$ and all quasi-finitely generated $B$-submodules $M$ of $A$ and $J$ is a set such that $\mathbb{Z}_{\blacksquare}[S]\cong \underset{J}{\prod}\mathbb{Z}$. The map $\phi$ is uniquely determined by the image of $1\in A$, that is, by an element of $\underset{\mathbb{N}\times J}{\prod} \underline{M}$ for some $M$. Every such element may be written as $\underset{i\in\mathbb{N}}{\sum} x_i\otimes y_i$, where $x_i=(0,0,\dots,0,1,0,\dots)\in \underset{\mathbb{N}}{\prod} B^+\subset  NS_{(\mathcal{A},\mathcal{M})}(*)$ with $1$ in the $i$-th position, $y_i\in \underset{J}{\prod} M$, and decomposable tensors are identified with their image under the canonical map $(\underset{\mathbb{N}}{\prod} B^+)\underset{{A}}{\otimes}(\underset{J}{\prod} M)\rightarrow \underset{\mathbb{N}\times J}{\prod} M$; this sum decomposition is just the obvious decomposition "into lines". Observe that this sum does converge in $\underset{\mathbb{N}\times J}{\prod} M$, which is endowed with the product topology, and $\{y_i\}_{i\in\mathbb{N}}$ is bounded by construction, since it is a subset of $M$. It follows directly from the construction used in the proof of Lemma~\ref{contfact} that the mapping $\sigma$ maps $(x_i)_{i\in\mathbb{N}}$ into a sequence $(f_i)_{i\in \mathbb{N}}\subset C(S,A)$ converging to $0$. Putting all this together, we see that the mapping $f:\mathcal{M}[S]\rightarrow\mathcal{M}[S]$ is induced by the map of the underlying modules and can be given by the formula $ m\mapsto \underset{i\in \mathbb{N}}{\sum} f_i(m)y_i,\ m\in \mathcal{M}[S]$, where we identify $C(S,A)$ and $\mathrm{Hom}_{(\mathcal{A},\mathcal{M})}(\mathcal{M}[S],\mathcal{A})$. 

Now let us analyze the sequence $\{f_i\}_{i\in \mathbb{N}}$. Identify $\mathcal{M}[S]$ with $\underset{J}{\prod} \mathbb{Z} \underset{\mathbb{Z}_{\blacksquare}}{\otimes} (\mathcal{A},\mathcal{M})$ and $\mathrm{Hom}_{ (\mathcal{A},\mathcal{M})}(\mathcal{M}[S],\mathcal{A})$ with $\mathrm{Hom}^{\mathrm{cont}}_{\mathbb{Z}}(\underset{J}{\prod}{\mathbb{Z}},{A})$. We claim that for each $i\geq 0$, there exists an at most countable set $K_i\subset J$ such that $f_i(\underset{J\setminus K_i}{\prod} \mathbb{Z})\equiv 0$. Indeed, fix $i$ and pick a pair of definition $(A_0,I)$ of $A$. Writing $A$ as $\varprojlim A/I^n$, we obtain $\mathrm{Hom}^{\mathrm{cont}}_{\mathbb{Z}}(\underset{J}{\prod}{\mathbb{Z}},{A})\cong \varprojlim\mathrm{Hom}^{\mathrm{cont}}_{\mathbb{Z}}(\underset{J}{\prod}{\mathbb{Z}},A/I^n)$. As each $A/I^n$ is discrete, there exists a finite set $K_{i,n}$ such that $f_i(\underset{J\setminus K_{i,n}}{\prod} \mathbb{Z})\equiv 0$ modulo $I^n$. It is evident that the set $K_i:=\underset{n\geq 0}{\cup}K_{i,n}$ satisfies the required property. Moreover, the set $K:=\underset{n\geq 0}{\cup}K_{i}$, which is at most countable, satisfies that property for each $i\geq 0$.

Writing $\underset{J}{\prod} \mathbb{Z}$ as $\underset{K}{\prod} \mathbb{Z}\oplus \underset{J\setminus K}{\prod} \mathbb{Z}$, we obtain a decomposition $\mathcal{M}[S]\cong Y \oplus Y'$ such that for each $i\geq 0$, $f_i(Y')\equiv 0$. Denote the inclusion map $Y\rightarrow \mathcal{M}[S]$ by $\iota$ and the projection map $\mathcal{M}[S]\rightarrow Y$ by $\pi$. As the restriction of $1_{\mathcal{M}[S]}-f$ to $Y'$ is just the identity map, we obtain the following commutative diagram:

\begin{center}
\begin{tikzcd}
\mathcal{M}[S] \arrow{r}{1_{\mathcal{M}[S]}-f}\arrow{d}{\pi} & \mathcal{M}[S]\arrow{d}{\pi} \\
Y \arrow[r, dashed]{} & Y
\end{tikzcd}
\end{center}
where the bottom map is supplied by the universal property of the cokernel of $Y'\hookrightarrow \mathcal{M}[S]$. Observe that it is also evidently equal to $1_Y-\pi f\iota$. It follows immediately from the construction of $K$ that $\pi f \iota$ is trace class as well; indeed, it is given by the formula $m\mapsto \underset{i\in \mathbb{N}}{\sum} \bar{f}_i(m)\bar{y}_i$, where $\bar{f}_i$ denotes the restriction of $f_i$ to $Y$ and $\bar{y}_i=\pi(y_i)$. It is straightforward to check that the kernel and cokernel (and, hence, the cone) of the bottom map are isomorphic to those of the top map. Thus, the upshot of the discussion so far is that we may assume that $J$ is at most countable. Of course, the case where $J$ is finite, is trivial.

Therefore, let us assume that $J=\mathbb{N}$. Denote by $e_l$ the element $(0,0,\dots,0,1,0,\dots)\in\underset{\mathbb{N}}{\prod} A$ with $1$ in the $l$-th position. Denote also the $k$-th coordinate of $f(e_l)$ by $d_{kl}$ and consider the matrix $D=(d_{kl})_{k,l\in\mathbb{N}}$. As the set $\{y_i\}_{i\in \mathbb{N}}$ is bounded and the sequence $(f_i)_{i\in\mathbb{N}}$ converges to zero, we easily see that for every $r\in\mathbb{N}$, there exists $R_r\in\mathbb{N}$ such that $d_{kl}\in I^r$ for every $k\in\mathbb{N},\ l\geq R_r$. Now decompose $Y$ as $(\underset{i=1}{\overset{R_1-1}{\prod}} \mathbb{Z} \oplus\underset{i=R_1}{\overset{\infty}{\prod}} \mathbb{Z}) \underset{\mathbb{Z}_{\blacksquare}}{\otimes} (\mathcal{A},\mathcal{M})$ and denote the first factor by $Y_1$ and the second by $Y_2$; additionally, denote the matrix $(d_{kl})_{k,l\geq R_1}$ by $\tilde{D}$ and the map $\pi f\iota$ by $\tilde{f}$, where $\pi$ is the canonical projection $Y\twoheadrightarrow Y_2$ and $\iota$ is the canonical inclusion $Y_2\hookrightarrow Y$. We claim that the sum $\underset{n\geq 0}{\sum}\tilde{f}^n$ is well-defined. To prove it, let us analyze each individual $\tilde{f}^n$. View it as a map $\tilde{f}^n:\underset{i=R_1}{\overset{\infty}{\prod}}\mathbb{Z}\rightarrow Y_2$; it is clearly given by the matrix $\tilde{D}^n$, all of whose coefficients lie in $I^n$. Write $Y_2$ as $\colim_M \underset{i=R_1}{\overset{\infty}{\prod}} \underline{M}$; since $\underset{i=R_1}{\overset{\infty}{\prod}}\mathbb{Z}$ is a compact $\mathbb{Z}_{\blacksquare}$-module, we see that there exists a finitely generated subring $B_n\subset A^+$ and a quasi-finitely generated submodule $M_n\subset A$ over $B_n$ such that the map $\tilde{f}^n$ factors through $ \underset{i=R_1}{\overset{\infty}{\prod}} \underline{M}_n$. Moreover, we may assume that all $B_n$'s are the same. Indeed, there exists a finitely generated subring $\tilde{B}\subset A^+$ and quasi-finitely generated submodules $N_1,N_2\subset A$ over $\tilde{B}$ such that for each $i\geq 0$, the image of $f_i$ lies inside $N_1$, and $\{y_i\}_{i\in \mathbb{N}}\subset N_2$; one easily verifies that we may assume that $B_n$ is equal to $\tilde{B}$ (and $M_n$ is equal to $N_1^n N_2$). Note that the image of $\tilde{f}^n$ clearly lies in $ \underset{i=R_1}{\overset{\infty}{\prod}} \underline{I}^n$; replacing $M_n$ by $M_n\cap I^n$, we may assume that the module $M_n$ comes from a tower $(M_{nk})_{k\in \mathbb{N}}$ such that $M_{nk}$ is trivial for each $k=0,\dots,n$. It therefore follows that the submodule $M:=\underset{i\geq 0}{+}M_n$ is quasi-finitely generated over $\tilde{B}$ and $\underset{n\geq 0}{\sum}\tilde{f}^n$ is well-defined. Thus, the mapping $1_{Y_2}-\tilde{f}$ is bijective and, hence, we may assume that it is just the identity map; in other words, we may assume that $1-f$ is given by a matrix of the form $\begin{pmatrix} (a_{kl})_{k,l=0}^{R_1-1} & *\\ * & E \end{pmatrix}$. Finally, it is easily verified that there exists an invertible endomorphism $F:Y\rightarrow Y$ such that $(1-f)\circ F$ is given by a matrix of the form $\begin{pmatrix} (a_{kl})_{k,l=0}^{R_1-1} & *\\ 0 & E \end{pmatrix}$. The desired statement now follows from the fact that the cone of any such map is isomorphic to the cone of the composition $Y_1\hookrightarrow Y\xrightarrow[]{1-f} Y\twoheadrightarrow Y_1$, which is clearly discrete.
\end{proof}

Combining Lemma~\ref{perf}, Lemma~\ref{two_duals}, and Theorem~\ref{d_is_d}, we obtain the following:

\begin{corollary}
Every dualizable object of $\mathcal{D}(\mathcal{A},\mathcal{M})$ comes from a perfect complex in $\mathcal{D}(A)$.
\end{corollary}

We conclude with the promised proofs of Theorem~\ref{descent of everything} and Theorem~\ref{pseudo stack}.

\begin{proof}[Proof of Theorems~\ref{descent of everything} and \ref{pseudo stack}]

Combing Theorems \ref{nuclear_descent}, \ref{dual_descent}, \ref{pseudocoherent objects descent}, and \ref{d_is_d}, we see that pseudocoherent and perfect complexes over $A$ satisfy descent. We now show that for each $m\in \mathbb{Z}$, pseudocoherent complexes in $\mathcal{D}^{\leq m}(A)$ also satisfy descent. Without loss of generality, we assume that $m=0$. Let $\{U_i\}_{i=1}^{n}$ be an affinoid open covering of $X$ and $M$ an a pseudocoherent complex in $\mathcal{D}^{-}(A)$ such that for each $i=1,\dots, n$, the localization $M\underset{A}{\overset{L}{\otimes}}A_{U_i}$ is a pseudocoherent complex in $\mathcal{D}^{\leq 0}(A)$. Applying Proposition~\ref{kedlnice}, we may assume that the covering is nice. Arguing by induction, we immediately reduce the proof to the case where $n=2$. Denote the localization $M\underset{A}{\overset{L}{\otimes}}A_{U_1}$ by $M_1$, $M\underset{A}{\overset{L}{\otimes}}A_{U_2}$ by $M_2$ and $M\underset{A}{\overset{L}{\otimes}}A_{U_1\cap U_2}$ by $M_{12}$. As $M$ is isomorphic to the homotopy pullback of $M_1\rightarrow M_{12}\xleftarrow[]{} M_{2}$, we see that its non-trivial cohomology groups may sit only in cohomological degrees $\leq 1$. It follows from Lemma~\ref{no first coh} that $H^1(M)$ is trivial; hence, $M$ is a pseudocoherent complex in $\mathcal{D}^{\leq 0}(A)$.

We now handle the case of (topological) stably pseudocoherent modules. Let $U$ be an affinoid open subspace of $X$ and $M$ a (topological) stably pseudocoherent module over $A$. It follows directly from the definitions that $M\underset{A}{\overset{L}{\otimes}}A_U$ is (topological) stably pseudocoherent over $A_U$. Conversely, let $\{U_1,U_2\}$ be a nice open covering of $X$ and $M$ a pseudocoherent complex over $A$ such that for $i=1,2$, the localization $M\underset{A}{\overset{L}{\otimes}}A_{U_i}\in \mathrm{Mod}_{A_{U_i}}$ is a (topological) stably pseudocoherent module over $A_{U_i}$. We need to show that $M$ itself is a (topological) stably pseudocoherent module over $A$. We note that its non-trivial cohomology groups may sit only in positive homological degrees as $M$ is isomorphic to the homotopy pullback of $M_1\rightarrow M_{12}\xleftarrow[]{} M_{2}$ (with notation as above); hence, $M$ is a pseudocoherent module over $A$. To show that it is stably pseudocoherent, we need to check that for every affinoid open subspace $U\subset X$, the localization $M\underset{A}{\overset{L}{\otimes}}A_U$ is concentrated in degree $0$. To prove it, we apply the above argument to $\{U\cap U_1,U\cap U_2\}$. We now note that if $M\underset{A}{\overset{L}{\otimes}}A_{U_i}\in \mathrm{Mod}_{A_{U_i}}$ is complete for $i=1,2$, then $M$ is complete as well. Indeed, it follows directly from the fact that for two-element covers, $M$ is isomorphic to the kernel of $M_1 \oplus M_2 \xrightarrow{f_1-f_2} M_{12}$. Finally, suppose that for $i=1,2$, the localization $M\underset{A}{\overset{L}{\otimes}}A_{U_i}$ is topological stably pseudocoherent and let $U$ be an affinoid open subspace. Applying the same argument as above to $\{U\cap U_1,U\cap U_2\}$, we see that $M\underset{A}{\otimes} A_U$ is complete as well.

Finally, we prove that perfect complexes with amplitude in $[a,b]$ satisfy descent (recall that $\mathrm{FP}_A$ is equal to $\mathrm{Perf}_A^{[0,0]}$). Without loss of generality, we assume that $b=0$. Let $\{U_1,U_2\}$ be a nice open covering of $X$ and $P$ a perfect complex in $\mathcal{D}^{\leq 0}(A)$ such that for $i=1,2$, the complex $P\underset{A}{\overset{L}{\otimes}}A_{U_i}$ lies in $\mathrm{Perf}_{A_{U_i}}^{[a,0]}$. We note that the complex $P^{\vee}$ is a perfect complex as well. As localization commutes with taking duals in the present situation (by Proposition~\ref{basicdualbasechange}) and for each $n\in \mathbb{Z}$, pseudocoherent complexes in $\mathcal{D}^{\leq n}(A)$ satisfy descent, we see that the tor-amplitude of $P^{\vee}$ lies in $[-\infty,-a]$. As $P$ is isomorphic to $(P^{\vee})^{\vee}$, we see that $P$ has tor-amplitude in $[a;+\infty)$. Therefore, $P$ has tor amplitude in $[a,0]$, as desired.
\end{proof}

\begin{lemma}
\label{no first coh}
Let $\{U_1,U_2\}$ be an affinoid open covering of $X$ and $f:A^n\rightarrow A^m$ a map such that the induced maps $f_{U_1}:A_{U_1}^n\rightarrow A_{U_1}^m$ and $f_{U_2}: A_{U_2}^n\rightarrow A_{U_2}^m$ are surjective. Then $f$ is surjective.
\end{lemma}

\begin{proof}
As the global section functor $\Gamma:\mathrm{QCoh}_{\spec A}\rightarrow \mathrm{Mod}_A$ is exact, it suffices to prove that for every closed point $\mathfrak{m}\in \spec A$, the induced map $f_{\mathfrak{m}}:A_{\mathfrak{m}}^n\rightarrow A_{\mathfrak{m}}^m$ is surjective. Represent the map $f$ by an $m\times n$-matrix $M$. By assumption, for every point $x\in X$, the map $f_x:\mathcal{O}_{X,x}^n\rightarrow \mathcal{O}_{X,x}^m$ is surjective. Therefore, for every point $x\in X$, there exists a minor of $M$ whose determinant is invertible in $\mathcal{O}_{X,x}$ (by Nakayama's lemma). Recall that for every maximal ideal $\mathfrak{m}\subset A$, there exists a valuation $\nu\in X$ such that $\nu(\mathfrak{m})=0$ (see, for example, the proof of \cite[Lemma 1.7.5]{Kdl}). Thus, for every closed point $\mathfrak{m}\subset A$, there exists a minor of $M$ whose determinant does not lie in $\mathfrak{m}$, and, hence, is invertible in $A_{\mathfrak{m}}$, as desired.
\end{proof}

\clearpage
\printbibliography

\textsc{Max-Planck-Institut für Mathematik, Vivatsgasse 7, 53111 Bonn, Germany}

\textsc{Universität Bonn, Endenicher Allee 60, 53115 Bonn, Germany}

\textit{E-mail address}: \texttt{griga@math.uni-bonn.de}
\end{document}